\documentclass[10pt,a4paper]{article}
\usepackage{amstext}
\usepackage{array}
\usepackage{amsfonts}
\usepackage{amssymb}
\usepackage{graphicx}
\usepackage{epstopdf}
\usepackage{algorithm}     %serve per scrivere lo pseudo codice
\usepackage{algpseudocode} %serve per scrivere lo pseudo codice
\usepackage{url}
\usepackage{caption}
\usepackage{subfig}
\usepackage{csquotes}
\usepackage{tikz}
\usetikzlibrary{arrows}
\usepackage{mathtools}
\usepackage{bbm}
\usepackage{amsthm}
\newcommand{\vertiii}[1]{{\left\vert\kern-0.25ex\left\vert\kern-0.25ex\left\vert #1 
   \right\vert\kern-0.25ex\right\vert\kern-0.25ex\right\vert}}

\newtheorem{theorem2}{Theorem}[section]

\newtheorem{remark}[theorem2]{Remark}

\usepackage{geometry}
\makeatletter
\def\blfootnote{\xdef\@thefnmark{$\star$}\@footnotetext}
\makeatother
\newenvironment{Authors}%
  {\begin{center}\begin{bfseries}}%
  {\end{bfseries}\end{center}}
\newenvironment{Addresses}%
  {\begin{flushleft}\begin{itshape}}%
  {\end{itshape}\end{flushleft}}
  \newcommand{\email}[1]{\hspace*{\stretch{1}}\emph{\texttt{#1}}}

 \usepackage{fancyhdr}  
  \fancypagestyle{plain}{
\fancyhead{}
\fancyhead[C]{\hfill Submitted to Elsevier, August 2023}
 }
\begin{document}

\thispagestyle{plain}

\title{Registration-based model reduction of parameterized PDEs with spatio-parameter adaptivity}
 \date{}
 \maketitle

 \maketitle
\vspace{-50pt} 
 
\begin{Authors}
Nicolas Barral$^{1}$, Tommaso Taddei$^{2}$, Ishak Tifouti$^{1}$
\end{Authors}

\begin{Addresses}
$^1$
Univ. Bordeaux, CNRS, Bordeaux INP, IMB, UMR 5251, F-33400 Talence, France\\ Inria Bordeaux Sud-Ouest, Team CARDAMOM, 33400 Talence, France, \email{nicolas.barral@inria.fr,ishak.tifouti@inria.fr} \\
$^2$
Univ. Bordeaux, CNRS, Bordeaux INP, IMB, UMR 5251, F-33400 Talence, France\\ Inria Bordeaux Sud-Ouest, Team MEMPHIS, 33400 Talence, France, \email{tommaso.taddei@inria.fr} \\
\end{Addresses}

%65N30 = Finite elements, Rayleigh-Ritz and Galerkin methods, finite methods
% 41A45 = Approximation by arbitrary linear expressions
% 35L02 = First-order hyperbolic equations
% 35J15 = Second-order elliptic equations
%90C26 = mathematical programming  nonconvex problems

%-----------------------------
%      your text
%-----------------------------

\begin{abstract}
We propose an automated nonlinear model reduction and mesh adaptation framework for rapid and reliable solution of parameterized advection-dominated problems,  with emphasis on compressible flows.
The key features of our approach are threefold:
(i) a metric-based mesh adaptation technique to generate an accurate mesh for a range of parameters,
(ii) a general (i.e., independent of the underlying equations) registration procedure for the computation of a mapping $\Phi$ that tracks moving features of the solution field, and
(iii) an hyper-reduced least-square Petrov-Galerkin  reduced-order model  for the rapid and reliable estimation of the mapped solution.
We discuss a general paradigm --- which mimics the refinement loop considered in mesh adaptation --- to simultaneously construct the high-fidelity and the reduced-order approximations, and we discuss actionable strategies to accelerate the offline phase.
We present extensive numerical investigations for a quasi-1D nozzle problem and for a two-dimensional inviscid flow past a Gaussian bump to display the many features of the methodology and to assess the performance for problems with discontinuous solutions.
\end{abstract}

\medskip

%
%\noindent
%\emph{Highlights}
%\begin{itemize}
%\item
%Application of model reduction to parametric problems with discontinuous solutions.
%\item
%Mesh adaptation reduces memory and CPU costs of parametric problems.
%\item
%Registration improves compression of the solution manifold.
%\item
%Registration facilitates the task of mesh adaptation.
%\end{itemize}

\section{Introduction}
\label{sec:intro}
\subsection{Lagrangian model   reduction of  steady conservation laws} 
In the past few decades, there has been an increasing demand for rapid and reliable reduced-order models (ROMs) for many-query and real-time applications  such as design optimization, uncertainty quantification, real-time control and monitoring.
Despite the many contributions to the field, model order  reduction of advection-dominated partial differential equations (PDEs) remains a formidable task that requires major improvements of state-of-the-art procedures.
The goal of this paper is to devise an integrated model order reduction (MOR) mesh adaptation (MA) procedure for nonlinear advection-dominated PDEs: our approach combines projection-based MOR, mesh adaptation and registration  techniques to simultaneously build a parsimonious yet accurate high-fidelity (HF) discretization, a low-rank representation of the solution field that depends on a modest number of generalized coordinates,  and a ROM that can be rapidly solved for new values of the parameters.
 
 {We consider PDE problems that depend on a vector of $P$ parameters.}
We denote by $\mu$ the vector of model parameters in the region $\mathcal{P}\subset \mathbb{R}^P$; we denote by $\Omega\subset \mathbb{R}^d$ the open computational domain\footnote{To simplify the
presentation, in the introduction we assume that the domain does not depend on the parameters; however, in the numerical examples, we shall consider the case of parameterized geometries.}; given the parametric field $w:\Omega\times \mathcal{P}\to \mathbb{R}$, we also introduce notation $w_{\mu}:=w(\cdot; \mu): \Omega\to \mathbb{R}$.
Given $\mu\in \mathcal{P}$, we denote by $q_{\mu}^{\rm true}:\Omega\to \mathbb{R}^D$ the vector of $D$ state variables that satisfies the hyperbolic  conservation law:
\begin{equation}
\label{eq:conservation_law}
\nabla \cdot {F}_{\mu}( q_{\mu}^{\rm true} ) = {S}_{\mu}(q_{\mu}^{\rm true}) \quad
{\rm in} \; \Omega,
\end{equation}
where $F: \mathbb{R}^D \times \mathcal{P} \to \mathbb{R}^{D \times d}$ is the physical flux and 
${S}: \mathbb{R}^D \times \mathcal{P} \to \mathbb{R}^{D}$ is the source term.
We further introduce the Hilbert space 
$\mathcal{X}:= [L^2(\Omega)]^D$   endowed with the inner product $(\cdot, \cdot)$ and the induced norm
$\| \cdot \|:=\sqrt{(\cdot,\cdot)}$, such that
$({w}, {v}) = \int_{\Omega} {w} \cdot  {v} \,  d{x}$ for all ${w}, {v} \in \mathcal{X}$;
 we define   the  solution manifold
$\mathcal{M} = \{q_{\mu}^{\rm true}: \mu\in \mathcal{P} \} \subset \mathcal{X}$ that collects the solutions to \eqref{eq:conservation_law} for all parameter values in the prescribed parameter range. 
We denote by 
$\mathcal{T}_{\rm hf} = \left(   \{   {x}_j^{\rm hf}\}_{j=1}^{N_{\rm nd}}, \texttt{T}   \right)$ a mesh of the domain $\Omega$ 
with nodes  $\{   {x}_j^{\rm hf}\}_j$ and  connectivity matrix $\texttt{T}$ (see section \ref{sec:formulation}); given the bijection $\Phi:\Omega\to \mathbb{R}^d$, we use notation 
$\Phi(\mathcal{T}_{\rm hf})$ to refer to the mesh with deformed nodes 
$\{  \Phi( {x}_j^{\rm hf} ) \}_j$ and the same connectivity $\texttt{T}$ as $\mathcal{T}_{\rm hf}$.

As extensively discussed in the MOR literature, 
effective model reduction of advection-dominated PDEs  is extremely challenging for state-of-the-art   procedures.
First, the vast majority of MOR methods rely on 
linear or affine approximations, that is 
\begin{equation}
\label{eq:linear_ansatz}
q_{\mu}^{\rm true} \approx \widehat{q}_{\mu}^{\rm lin} =  \texttt{Z} \widehat{\boldsymbol{\alpha}}_{\mu},
\end{equation}
where  $\texttt{Z}:\mathbb{R}^n \to \mathcal{X}$ is a linear or affine operator,
 and $\widehat{\boldsymbol{\alpha}}:\mathcal{P} \to \mathbb{R}^n$ is a  function of the parameter --- in the MOR literature, 
$\texttt{Z}$ is typically dubbed   \emph{reduced-order basis} (ROB), while  
 $\widehat{\boldsymbol{\alpha}}_{\mu}$ are referred to as \emph{generalized coordinates}.
As shown in several studies (e.g.,
\cite{ohlberger2016reduced}), linear methods  are fundamentally ill-suited to deal with parameter-dependent sharp gradients that naturally  arise in the  solutions to conservation laws of the form \eqref{eq:conservation_law}.
Second, MOR methods typically rely on a single HF mesh to describe all elements of the solution manifold $\mathcal{M}$.
For advection-dominated problems, MA  is of paramount importance for computational tractability. However, if parametric variations strongly affect the location of sharp-gradient regions, 
we
are forced to refine  the mesh over a vast portion of the domain $\Omega$, which leads 
to HF discretizations of intractable size. Effective MOR procedures for conservation laws should thus embed an effective parametric MA strategy to track moving structures. 

The provable inadequacy of linear ansatzs \eqref{eq:linear_ansatz} 
for  conservation laws has  motivated the development
of several nonlinear approximation methods
\cite{amsallem2008interpolation,barnett2022quadratic,lee2020model,peherstorfer2020model};
a  promising class of nonlinear approximations is given by 
\emph{Lagrangian  methods}
\cite{ching2022model,iollo2014advection,mirhoseini2023model,mojgani2017arbitrary,ohlberger2013nonlinear,taddei2015reduced,taddei2020registration} based on the ansatz
\begin{equation}
\label{eq:lagrangian}
\widehat{q}_{\mu} = 
\widetilde{q}_{\mu}
\circ \Phi_{\mu}^{-1},
\quad
{\rm where} \;\;
\widetilde{q}_{\mu} =   \texttt{Z} \widehat{\boldsymbol{\alpha}}_{\mu}, \;\;
\Phi_{\mu} = \texttt{N} \left(  \widehat{\mathbf{a}}_{\mu} \right).
 \end{equation}
As in \eqref{eq:linear_ansatz},
$\texttt{Z} :\mathbb{R}^n \to \mathcal{X}$ is a linear (or affine) operator, and
$\widehat{\boldsymbol{\alpha}}:\mathcal{P} \to \mathbb{R}^n$ is a  vector-valued function of  generalized coordinates;
on the other hand,   $\texttt{N} :\mathbb{R}^m \to {\rm Lip}(\Omega; \mathbb{R}^d)$ is a suitable, possibly nonlinear, operator that is informed by  the domain $\Omega$ and
$\widehat{\mathbf{a}} :\mathcal{P} \to \mathbb{R}^m$ is a  vector-valued function of  generalized coordinates
for the mapping.

Lagrangian approaches are motivated by the observation (see, e.g., \cite{iollo2022mapping} and \cite{taddei2020registration}) that for many problems in continuum mechanics coherent structures that are troublesome for linear approximations --- such as shear layers, wakes, shocks and cracks --- vary smoothly with the parameter.
The mapping $\Phi: \Omega\times \mathcal{P} \to \Omega$ in \eqref{eq:lagrangian} should hence be 
 designed to track moving features of the solution field and ultimately improve the compressibility of the mapped solution manifold
$\widetilde{\mathcal{M}} = \{ 
\widetilde{q}_{\mu}^{\rm true}:=
q_{\mu}^{\rm true}\circ \Phi_{\mu}  : \mu\in \mathcal{P}\}$.
The task of finding the mapping $\Phi$ based on approximate snapshots of the solution manifold $\mathcal{M}$ is referred to as \emph{registration problem} \cite{taddei2020registration}.
Note that, by tracking sharp features of the solution field, registration  facilitates also the task of building a common mesh for all elements of the (mapped) solution manifold: the mapping $\Phi$ hence provides a systematic way to perform parameter-dependent r-adaptivity \cite{budd2009adaptivity,mcrae2018optimal}.

\subsection{Adaptive construction of Lagrangian reduced-order models}
In this paper, we propose a general paradigm for the simultaneous construction of the HF and reduced-order approximations, which mimics the refinement loop considered in MA.
The general procedure is sketched in Algorithm \ref{alg:offline_training}. 
Given an initial mesh $\mathcal{T}_{\rm hf}^{(0)}$ of $\Omega$ and  the training set 
$\mathcal{P}_{\rm train} = \{\mu^k \}_{k=1}^{n_{\rm train}} \subset \mathcal{P}$, our method returns an HF mesh $\mathcal{T}_{\rm hf}$, a low-rank mapping $\Phi$, a ROB $\texttt{Z}$ and a ROM   for the generalized coordinates $\widehat{\boldsymbol{\alpha}}$ (cf. \eqref{eq:lagrangian}) based on an iterative procedure that comprises four distinct steps.

\begin{enumerate}
\item
\underline{Snapshot generation:} 
$
\left( \mathcal{T}_{\rm hf}  , \Phi ,  \mathcal{P}_{\rm train}  \right)
 \rightarrow 
\left\{ q_{\mu}^{{\rm hf}} \, : \,  \mu\in \mathcal{P}_{\rm train} \right\}$.
We generate snapshots of the solution field for all values of the parameter $\mu$ in  
$ \mathcal{P}_{\rm train}$ based on the parametric mesh $\mu \mapsto \Phi_{\mu}( \mathcal{T}_{\rm hf} )$.
\item
\underline{Mesh adaptation:}
$\left\{  \widetilde{q}_{\mu}^{{\rm hf}} :=   q_{\mu}^{{\rm hf}} \circ \Phi_{\mu} :  \mu\in \mathcal{P}_{\rm train}
\right\} \rightarrow  \mathcal{T}_{\rm hf}$. We exploit the available set of snapshots  to generate an accurate yet parsimonious mesh for the elements of the mapped manifold  $\widetilde{\mathcal{M}}$.
\item
\underline{Registration:}
$\left( 
\left\{  q_{\mu}^{{\rm hf}} \, : \, \mu\in \mathcal{P}_{\rm train}
\right\}
, \;\; 
\mathcal{T}_{\rm hf} 
\right) \rightarrow  \Phi 
$. We exploit the  available set of snapshots 
to find a parametric mapping $\Phi$ that tracks coherent, parameter-dependent structures of the solution field.
The method should ensure that the deformed mesh $\Phi_{\mu}(  \mathcal{T}_{\rm hf}  )$ is a proper mesh of $\Omega$ for all $\mu \in \mathcal{P}$.
\item
\underline{Linear-subspace model reduction:}
$\left( \mathcal{T}_{\rm hf}, \, \Phi,  \,
\mathcal{P}_{\rm train}
\right)
 \rightarrow \left(   \texttt{Z} ,  {\rm ROM} \right)
$. We apply linear-subspace MOR to determine the low-rank expansion 
$\mu \mapsto  \widetilde{q}_{\mu} $, that is we build the ROB $ \texttt{Z}$ and the ROM for 
$\mu \mapsto   \widehat{\boldsymbol{\alpha}}_{\mu}$.
\end{enumerate}

\begin{algorithm}[H]                      
\caption{: adaptive training procedure.}     
\label{alg:offline_training}     
 \normalsize 

\begin{algorithmic}[1]
\State
Initialization: define the mesh
$\mathcal{T}_{\rm hf}^{(0)} =\mathcal{T}_{\rm hf}^{(1)}$,   the mapping
$\Phi^{(0)} = \texttt{id}$ (identity map), and the training set
$\mathcal{P}_{\rm train} = \{\mu^k \}_{k=1}^{n_{\rm train}} \subset \mathcal{P}$.
\vspace{3pt}

\For {$k=1, \ldots, N_{\rm it}$ }

\State
Snapshot generation
\hfill 
$
\left( \mathcal{T}_{\rm hf}^{(k-1)} , \Phi^{(k-1)},  \mathcal{P}_{\rm train}  \right)
 \rightarrow 
\left\{ q_{\mu}^{{\rm hf}, (k)} \, : \,  \mu\in \mathcal{P}_{\rm train} \right\}.$ 
\vspace{3pt}

\If {$k>1$}
\State
Mesh adaptation  (cf. section \ref{sec:mesh_adaptation})
\hfill
$
\left\{  q_{\mu}^{{\rm hf}, (k)} \circ \Phi_{\mu}^{(k-1)} :  \mu\in \mathcal{P}_{\rm train}
\right\}
 \rightarrow  \mathcal{T}_{\rm hf}^{(k)}
$.

\EndIf

\State 
Registration (cf. section \ref{sec:registration})
\hfill
$\left( 
\left\{  q_{\mu}^{{\rm hf}, (k)} \, : \, \mu\in \mathcal{P}_{\rm train}
\right\}
, \;\; 
\mathcal{T}_{\rm hf}^{(k)}
\right) \rightarrow  \Phi^{(k)}
$.
\smallskip

\State
 Linear-subspace model reduction (cf. section \ref{sec:pMOR})
\hfill
$\left(\mathcal{T}_{\rm hf}^{(k)}, \, \Phi^{(k)},  \, \mathcal{P}_{\rm train} \right)
\;\; 
 \rightarrow 
\left(\texttt{Z}^{(k)},     {\rm ROM}^{(k)} \right)$.
 
\EndFor
\end{algorithmic}
\end{algorithm}

Similarly to the standard MA loop, our method relies on multiple iterations to address the inaccuracy of the HF estimates at early iterations.
We show that the iterative procedure in Algorithm \ref{alg:offline_training} can be significantly accelerated using information from previous iterations (cf. section \ref{sec:adaptive_training}).

The outline of the paper is as follows.
In section \ref{sec:formulation}, we introduce relevant notation and the two model problems considered for  numerical assessment. In sections  \ref{sec:mesh_adaptation}, 
\ref{sec:registration},
and \ref{sec:pMOR}, we discuss the problems of mesh adaptation,
registration, 
 and model reduction; 
 in section \ref{sec:adaptive_training} we discuss how to accelerate the training procedure by exploiting information from previous iterations; in section \ref{sec:numerics}, we present extensive numerical investigations to illustrate the effectiveness of our approach.
Section \ref{sec:conclusions} concludes the paper.

\subsection{Contributions and relation to previous works}
This paper extends  the work in \cite{ferrero2022registration}
in several ways:  
first, we propose an adaptive, iterative procedure for the simultaneous construction of the HF mesh, the mapping $\Phi$, and the reduced-order approximation for the mapped field; 
second, we incorporate an automated parametric mesh adaptation strategy that is directly informed by the estimated solution fields; 
third, we discuss viable new strategies to accelerate the training procedure.
As in \cite{ferrero2022registration}, 
the building blocks of Algorithm \ref{alg:offline_training} exploit methodologies from previous works.
%As in \cite{ferrero2022registration}, 
The  registration procedure was first proposed in \cite{taddei2020registration} and then extended in  \cite{taddei2021registration,taddei2021space}, while 
 we rely on a projection-based least-squares Petrov-Galerkin (LSPG, \cite{carlberg2017galerkin,carlberg2013gnat}) ROM
 with empirical test space chosen as in \cite{taddei2021space}, hyper-reduction based on a variant of the  mesh-sampling/empirical quadrature  procedures first proposed in  \cite{farhat2015structure,yano2019discontinuous}, and 
 \emph{discretize-then-map} treatment of 
 geometry variations 
(cf. \cite{dal2019hyper,taddei2021discretize,washabaugh2016use}).
We also observe that the idea of using ROMs and/or HF models of variable fidelity at training stage to reduce training costs has been explored in
\cite{feng2023accelerating,kast2020non}.

We rely on a metric-based  approach
(\cite{loseille2011continuous1,loseille2011continuous2})
to mesh adaptation.
Given the HF field $q_{\mu}^{\rm hf}$, we compute the Hessian of the Mach number in the reference configuration to determine the metric $\mathfrak{M}_{\mu}$ for the parameter $\mu\in \mathcal{P}_{\rm train}$; then, we resort to metric intersection
\cite{barral2015thesis,barral2017time}
to devise a common metric for the entire snapshot set.
In this work, we rely on the open-source mesh adaptation toolkit
\texttt{mmg2d} (\cite{arpaia2022h,dapogny2014three}) to generate adapted meshes from a (possibly anisotropic) metric  $\mathfrak{M}$.

As discussed in section \ref{sec:registration}, the mapping coefficients
$\widehat{\mathbf{a}}_{\mu}$ in \eqref{eq:lagrangian} are computed using a non-intrusive  (regression) approach, while  the solution generalized coordinates 
$\widehat{\boldsymbol{\alpha}}_{\mu}$ are computed using projection;
on the other hand, Mirhoseini and Zahr in \cite{mirhoseini2023model} have recently proposed a coupled Lagrangian MOR approach to simultaneously  learn solution and mapping generalized coordinates.
As opposed to
\cite{mirhoseini2023model}, our choice enables the use of standard projection-based linear-subspace MOR methods for parameterized geometries: it hence has the potential to achieve faster online predictions and much easier integration with existing HF and MOR routines, possibly at the price of larger offline training costs.

Several authors have proposed to  include
(parametric) mesh adaptation procedures in the MOR framework.
Simultaneous adaptivity in space   and in parameter  
{--- in effect, spatio-parameter adaptivity ---}
 was proposed by Yano in \cite{yano2018reduced} and further developed  in \cite{sleeman2022goal}.
The approaches in \cite{sleeman2022goal,yano2018reduced}   exploit 
 $h$-MA and weak-greedy sampling of the parameter space, and 
  rely on the explicit  instantiation of a \emph{super-mesh} over the entire parameter domain; the size of the  super-mesh might hence be prohibitively large for  advection-dominated problems.   
 %Simultaneous adaptivity in space --- via $h$-MA --- and in parameter --- via greedy sampling ---
% was proposed by Yano in \cite{yano2018reduced} and further developed  in \cite{sleeman2022goal};  {\color{red}and dubbed \emph{spatio-parameter adaptivity}}.
To address this issue, 
Little and Farhat have proposed 
in \cite{little2023nonlinear}
to combine $h$-MA with clustering in parameter domain --- more precisely in state space --- to avoid the explicit definition of a super-mesh that is valid for all parameters.
The present work represents the first attempt to systematically combine parametric $r$-MA (induced by the mapping) with parameter-independent $h$-MA, in the  MOR framework.

We finally  remark that several authors have  considered ansatzs of the form 
\begin{equation}
\label{eq:lagrangian_modified}
\widehat{q}_{\mu} = \widetilde{q}_{\mu} \circ \Phi_{\mu},
\quad
{\rm where} \;\;\; 
\widetilde{q}_{\mu}: \mathbb{R}^d \to \mathbb{R}^D,
\quad
\Phi_{\mu}: \Omega \to \mathbb{R}^d,
\end{equation}
with  $\widetilde{q}_{\mu},\Phi_{\mu}$
possibly nonlinear    low-rank operators (see, e.g.,
\cite{black2021efficient,krah2023front}).
Note that, 
unlike in \eqref{eq:lagrangian},  
$\Phi_{\mu}$ is not necessarily a bijection from $\Omega$ in itself; note also that the field
$\widetilde{q}_{\mu}$ needs to be defined over $\mathbb{R}^d$.
As shown in \cite{krah2023front}, approximations of the form \eqref{eq:lagrangian_modified}
can potentially handle  shock topology changes;
on the other hand,  
we remark that \eqref{eq:lagrangian_modified} is inherently nonlinear and hence requires the development of specialized projection techniques.
%
%\textbf{MISSING:} mention that mappings have also been used in an Eulerian setting.
%\cite{krah2023front}, paper by Benner \cite{sarna2021data}. 

\section{Problem statement and finite element discretization}
\label{sec:formulation}
We consider the problem of approximating the solution to the parameterized Euler equations; we consider the equations in non-dimensional form.
We refer to \cite{toro2013riemann} for a thorough introduction to the mathematical model and to its physical interpretation. We denote by $\rho$ the fluid density, by $u$ the velocity  field, by $E$ the total energy and by $p$ the (static) pressure; we consider 
the following relationship between pressure and conserved
variables:
\begin{subequations}
\label{eq:constitutive_law}
\begin{equation}
p(q) = (\gamma - 1) \left( E  - \frac{1}{2} \rho \| u  \|_2^2\right),
\end{equation}
where $\gamma$ is the ratio of specific heats, which is set equal to $1.4$. We further introduce the speed of sound $a$, the Mach number ${\rm Ma}$, the total temperature $T_{\rm tot}$, the total pressure $p_{\rm tot}$ and the total enthalpy $H_{\rm tot}$ such that
\begin{equation}
\begin{array}{l}
\displaystyle{
a = \sqrt{ \gamma \frac{p}{\rho} }, 
\;\;
{\rm Ma} =   \frac{\|u \|_2}{a} , 
\;\;
T = \frac{p}{R \rho},
\;\;
T_{\rm tot} = T \left( 1+ \frac{\gamma-1}{2}  {\rm Ma} ^2\right),
} \\[3mm]
\displaystyle{
p_{\rm tot} = p \left( 1+ \frac{\gamma-1}{2}  {\rm Ma} ^2\right)^{\frac{\gamma-1}{\gamma}},
H_{\rm tot} = \frac{E+p}{\rho},
\;\;
R=\gamma-1.
}\\
\end{array}
\end{equation}
\end{subequations}

We introduce the finite element (FE) mesh
 $\mathcal{T}_{\rm hf} = \left(   \{   {x}_j^{\rm hf}\}_{j=1}^{N_{\rm hf}}, \texttt{T}   \right)$ of the domain $\Omega\subset \mathbb{R}^d$:
 the points 
 $\{   {x}_j^{\rm hf}\}_{j=1}^{N_{\rm nd}} \subset \overline{\Omega}$ are the nodes of the mesh,
 the matrix 
 $\texttt{T} \in \mathbb{N}^{n_{\rm lp}, N_{\rm e}}$ is 
the   connectivity matrix where $n_{\rm lp}$ is the number of degrees of freedom in each element and $N_{\rm e}$ is the number of elements.
We define the reference element 
$\widehat{\texttt{D}} = \{ {x}\in (0,1)^d:  \sum_{i=1}^d x_i <  1  \}$,
the   space $\mathbb{P}_{\texttt{p}}(\widehat{\texttt{D}})$ of  polynomials of degree lower or equal  to $\texttt{p}$,
the Lagrangian basis $\{ \ell_i \}_{i=1}^{n_{\rm lp}}$  of the polynomial space $\mathbb{P}_{\texttt{p}}(\widehat{\texttt{D}})$ associated with the nodes 
 $\{  \tilde{x}_i  \}_{i=1}^{n_{\rm lp}} \subset \overline{\widehat{\texttt{D}}}$; then, we define 
the elements
$\{ \texttt{D}_k  \}_{k=1}^{N_{\rm e}}$ as the images of the FE maps 
$ {\Psi}_k^{\rm hf}:\widehat{ \texttt{D} } \to \texttt{D}_k$ such that
\begin{equation}
\label{eq:psi_mapping}
 {\Psi}_k^{\rm hf}( {\tilde{x}})
 =
 \sum_{i=1}^{n_{\rm lp}} \;
 {x}_{ \texttt{T}_{i,k}   }^{\rm hf}  \; \ell_i(\tilde{x}),
 \quad
 k=1,\ldots,N_{\rm e}.
\end{equation}
We also define the FE space associated with the mesh
$\mathcal{T}_{\rm hf}$,
\begin{equation}
\label{eq:DG_space}
\mathcal{X}_{\rm hf} = \left\{
v\in [ L^2(\Omega) ]^D\, : \, v\circ  {\Psi}_k^{\rm hf} \in [ \mathbb{P}_{\texttt{p}}
(\widehat{\texttt{D}})]^D,
\;\;k=1,\ldots,N_{\rm e}
\right\},
\end{equation}
where $D=d+2$ corresponds to the number of state variables.
If $u\in \mathcal{X}_{\rm hf}$, we denote by $\mathbf{u}\in \mathbb{R}^{N_{\rm hf}}$ the corresponding FE vector such that
\begin{equation}
\label{eq:vector2field}
u \Big|_{  \texttt{D}_{k}  }  \, = \,  
 \sum_{i=1}^{n_{\rm lp}} \; 
  \sum_{\ell=1}^{D} 
 \left( \mathbf{u}  \right)_{\texttt{I}_{i,k,\ell}} \; \,  
 \ell_{i,k} \;\; e_{\ell},
 \qquad
 k=1,\ldots,N_{\rm e},
 \;\;
 \texttt{I}_{i,k,\ell} = i+(k-1) n_{\rm lp} + (\ell - 1)  n_{\rm lp} N_{\rm e},
 \end{equation}
 where $e_1,\ldots,e_d$ are the vectors of the canonical basis of $\mathbb{R}^D$, and 
 $N_{\rm hf} = D \cdot n_{\rm lp} \cdot N_{\rm e}$.
 
In view of the FE formulation, we introduce the facets 
$\mathcal{F}_{\rm hf} = \{   \texttt{F}_j \}_{j=1}^{N_{\rm f}}$ of the mesh: for each facet, we denote by
$\mathbf{n}^+$ the positive normal\footnote{The choice of the positive normal is arbitrary for internal facets and coincides with the outward normal to $\Omega$ for boundary facets.} to the facet and we define the element
 $\texttt{D}_j^+$  
 that contains $\texttt{F}_j$ and whose normal on $\texttt{F}_j$  is equal to 
 $\mathbf{n}^+$; for internal facets, we also define the element  
$\texttt{D}_j^-$ such that 
$\texttt{D}_j^+ \cap \texttt{D}_j^- = \texttt{F}_j$.
We also define the restriction operators
$E_k: \mathcal{X}_{\rm hf} \to [ L^2( \texttt{D}_k  ) ]^D$ and
$E_j^{\pm}: \mathcal{X}_{\rm hf} \to [ L^2( \texttt{D}_j^{\pm}  ) ]^D$ such that
$E_k u = u|_{\texttt{D}_k}$ and
$E_j^{\pm} u = u|_{\texttt{D}_j^{\pm}}$ for 
$k=1,\ldots,N_{\rm e}$ and $j=1,\ldots,N_{\rm f}$.

\begin{remark}
\label{remark:lagrangian_ansatz}
Exploiting \eqref{eq:vector2field}, we find that any FE field $u\in \mathcal{X}_{\rm hf}$ is uniquely characterized by the pair $(\mathcal{T}_{\rm hf}, \mathbf{u})$. Given the bijection $\Phi:\Omega\to \mathbb{R}^d$, we introduce the deformed mesh
$\Phi( \mathcal{T}_{\rm hf}   ) = \left(
\{   \Phi(x_j^{\rm hf}) \}_{j=1}^{N_{\rm hf}},
\texttt{T} \right)$, and the corresponding  FE maps
$\{ \Psi_{k,\Phi}^{\rm hf} \}_k$ and FE space
$\mathcal{X}_{\rm hf, \Phi}$: it is easy to verify that if 
$(\Phi( \mathcal{T}_{\rm hf} ), \mathbf{u})$  interpolates the field $u$ in the nodes of 
$\Phi( \mathcal{T}_{\rm hf} )$, then 
$( \mathcal{T}_{\rm hf} , \mathbf{u})$ interpolates $u\circ \Phi$ in the nodes of $\mathcal{T}_{\rm hf} $. This implies that the ansatz \eqref{eq:lagrangian} can be stored as
\begin{equation}
\label{eq:lagrangian_discrete}
\mu\in\mathcal{P} \mapsto 
\left(
\Phi_{\mu}(  \mathcal{T}_{\rm hf}  ), \mathbf{Z} \widehat{\boldsymbol{\alpha}}_{\mu}
\right)
\end{equation}
where $\mathbf{Z}\in \mathbb{R}^{N_{\rm hf} \times n}$ is a parameter-independent matrix.
\end{remark}

\subsection{Finite element formulation}
\label{sec:FEM_notation}
We consider a discontinuous Galerkin (DG) FE formulation of the compressible Euler equations. We rely on a Laplacian artifical viscosity model (see, e.g., \cite{persson2006sub}) based on the piecewise-constant dilation-based viscosity:
\begin{equation} 
\label{eq:viscosity_dilationbased}
\nu |_{\texttt{D}_k}
\, = \,
c_{\rm \nu}
\left(
\frac{h_k}{\texttt{p}}
\right)^2
\int_{\texttt{D}_k} (-\nabla \cdot u)_+ \, dx,
\quad
{\rm with} \;\;
h_k = | \texttt{D}_k   |^{1/d}, 
\;\;
c_{\rm \nu}>0.
\end{equation}
We refer to \cite{yu2020study} for a thorough review of artificial viscosity models for DG formulations.
We consider the local Lax-Friedrichs (Rusanov) convective flux, and symmetric interior penalty diffusive flux. 
If we denote by $q$ the vector of state variables, the DG formulation
of the conservation law \eqref{eq:conservation_law}  can be stated as follows:
find $q^{\rm hf} \in \mathcal{X}_{\rm hf}$ such that
\begin{subequations}
\label{eq:abstract_HFM}
\begin{equation}
\mathfrak{R}^{\rm hf}(q^{\rm hf}, v)
=
\sum_{k=1}^{N_{\rm e}} r_k^{\rm e}(E_k q^{\rm hf}, E_k v  ) \; + \;
\sum_{j=1}^{N_{\rm f}} r_j^{\rm f}(
E_j^{+} q^{\rm hf}, E_j^{+} q^{\rm hf},
 E_j^{+} v,  E_j^- v  ),
 \quad
 \forall \, v\in  \mathcal{X}_{\rm hf};
\end{equation}
where the elemental residuals $r_k^{\rm e}$ are given by 
\begin{equation}
r_k^{\rm e}(q,v )
=
\int_{\texttt{D}_k} \left(
- F(q) + \nu(q) \nabla q 
\right) : \nabla v \, - \, S(q) \cdot v 
\, dx,
\quad
k=1,\ldots,N_{\rm e},
\end{equation}
while the facet residuals $r_j^{\rm f}$ are given by 
\begin{equation}
r_j^{\rm f}(q,v )
=
\int_{\texttt{F}_j} 
\mathfrak{H}(q, \mathbf{n^+}) \cdot J (v)  \, dx
-
\int_{\texttt{F}_j \setminus \partial \Omega} 
\left\{  \nu(q) \nabla q \mathbf{n}^+ \right\} \cdot J (v) 
+
\left\{  \nu(q) \nabla v \mathbf{n}^+ \right\} \cdot J (q) 
-
\frac{\gamma_{\rm ip}}{|\texttt{F}_j|} J(q) \cdot J(v) 
\, dx,
\end{equation}
for $j=1,\ldots,N_{\rm f}$.
Here,
 $q^{\pm}(x) = \lim_{\epsilon\to 0^+} q(x  \mp  \epsilon \mathbf{n}^+)$,  $J(v) = v^+-v^-$ and 
$\{ v \} = \frac{1}{2} (v^++v^-)$ if $x\notin \partial \Omega$, while
$J(v) = \{ v \} = v$  if $x\in \partial \Omega$; finally, 
  $\mathfrak{H}(q, \mathbf{n^+})$ is the convective numerical flux which embeds the definition of the boundary conditions and for which we omit the explicit expression
(see, e.g., \cite[Appendix B]{fidkowski2004high} for  the details).
\end{subequations}

In the numerical experiments, we consider 
  polynomials of degree $\texttt{p}=2$; we set 
$c_{\rm \nu}=0.1$ for the model problem of  section 
\ref{sec:model_problems_nozzle} and 
$c_{\rm \nu}=10$ for the model problem of section   
\ref{sec:model_problems_bump}; on the other hand, we consider $\gamma_{\rm ip}=10 \texttt{p}^2$.
Finally, we solve
the discrete problem \eqref{eq:abstract_HFM}  using the  pseudo-transient continuation (PTC) strategy 
discussed in \cite{yano2011importance}; in the absence of prior information about the solution field, we initialize the iterative procedure with the free-stream solution.

\subsection{Model problems}
\label{sec:model_problems}

\subsubsection{Inviscid flow through a nozzle}
\label{sec:model_problems_nozzle}
\begin{figure}[H]
\centering
\subfloat[] 
{\includegraphics[width=0.4\textwidth]
 {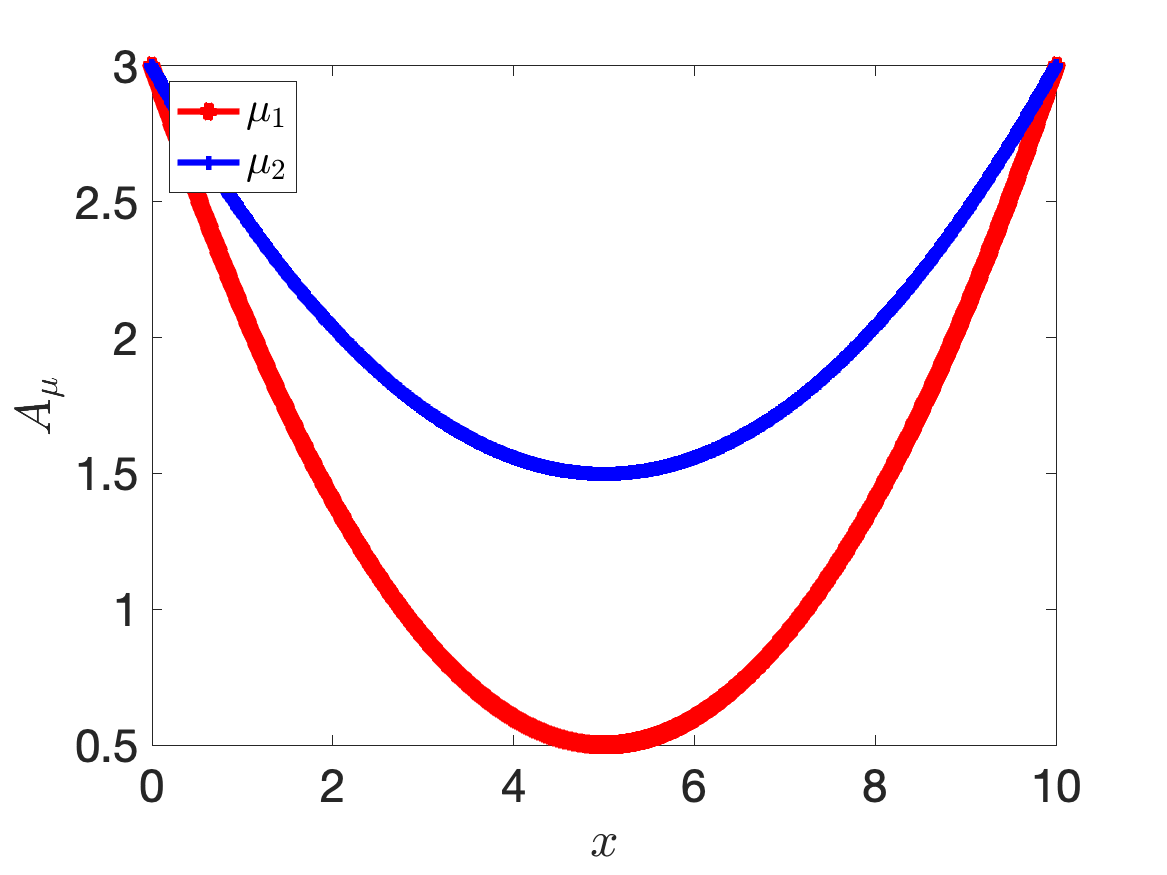}}
 ~~
   \subfloat[] 
{\includegraphics[width=0.4\textwidth]
 {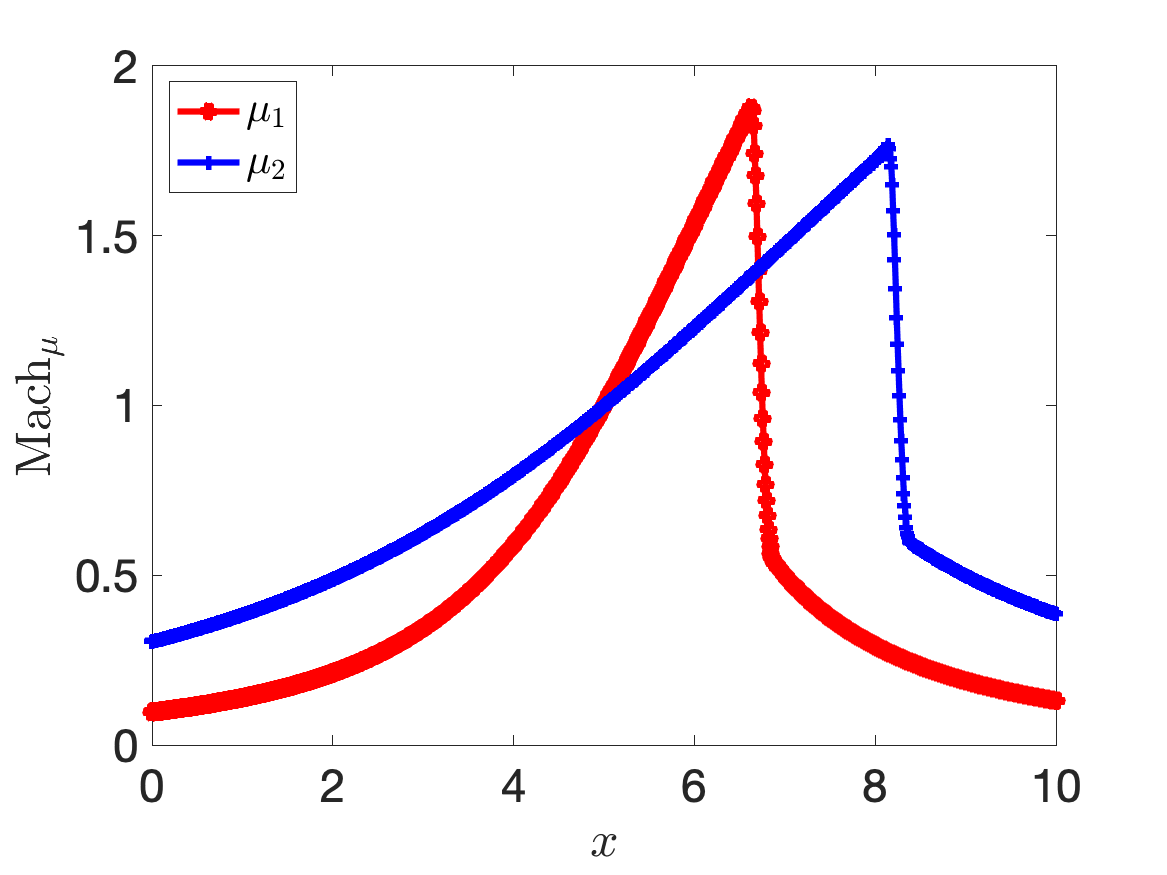}}
 
 \caption{inviscid flow through a nozzle. 
 (a)-(b) throat area and Mach field for  $\mu_1=[0.5,0.7]$ and  $\mu_2=[1.5,0.7]$.
}
 \label{fig:nozzle_vis}
 \end{figure}

We study the inviscid transonic flow of an ideal gas through a converging-diverging duct.
We define the domain $\Omega=(0,L)$, the area $A:\Omega \to \mathbb{R}_+$, the state $q$, the flux $F$ and the source term $S$ such that
\begin{subequations}
\label{eq:nozzle_flow}
\begin{equation}
q = \left[
\begin{array}{l}
A \rho \\
A \rho u \\
A E \\
\end{array}
\right],
\quad
F(q) = \left[
\begin{array}{l}
A \rho u  \\
A (\rho u^2 + p) \\
A u(E+p) \\
\end{array}
\right],
\quad
S(q) = \left[
\begin{array}{l}
0 \\
p   \partial_x A  \\
0 \\
\end{array}
\right],
\quad
A(x) = 3 + 4(A_0 - 3) \frac{x}{L} \left (   1 - \frac{x}{L}  \right).
\end{equation}
Then, we consider the conservation law:
\begin{equation}
\partial_x F(q^{\rm true}) = S(q^{\rm true}) \quad x\in \Omega,
\end{equation}
completed with a subsonic inlet condition where we prescribe total pressure $p_{\rm tot} = 0.95$ and
total temperature $T_{\rm tot} = 0.95$, and a subsonic outlet condition where we prescribe the static pressure  $p_0$.  
Note that the free-stream field is uniquely determined by the data
$p_{\rm tot}, T_{\rm tot}, p_0$ through \eqref{eq:constitutive_law}.
We set $L=10$; furthermore, we consider the parameter vector $\mu=[A_0,p_0]$ in the region 
$\mathcal{P}  = [0.5,1.5] \times [0.7,0.85]$. Figures
\ref{fig:nozzle_vis}(a) and (b) show the behavior of the area throat and the Mach number for two parameter values.
\end{subequations}

\subsubsection{Inviscid flow over a Gaussian bump}
\label{sec:model_problems_bump}
\begin{figure}[H]
\centering
 \subfloat[] 
{
\begin{tikzpicture}[scale=3]
\linethickness{0.3 mm}
\linethickness{0.3 mm}

\draw[ultra thick]  (-1.5,0)--(-1.5,0.8)--(1.5,0.8)--(1.5,0);

\draw[smooth, domain = 0:1.5, color=blue, very thick] plot (\x,{0.0625*exp(-25*\x^2)});

\draw[smooth, domain = 0:1.5, color=blue, very thick] plot (-\x,{0.0625*exp(-25*\x^2)});

\coordinate [label={above:  {\Huge {$\Omega$}}}] (E) at (1, 0.4) ;
\end{tikzpicture}
 }

\subfloat[] 
{\includegraphics[width=0.4\textwidth]
 {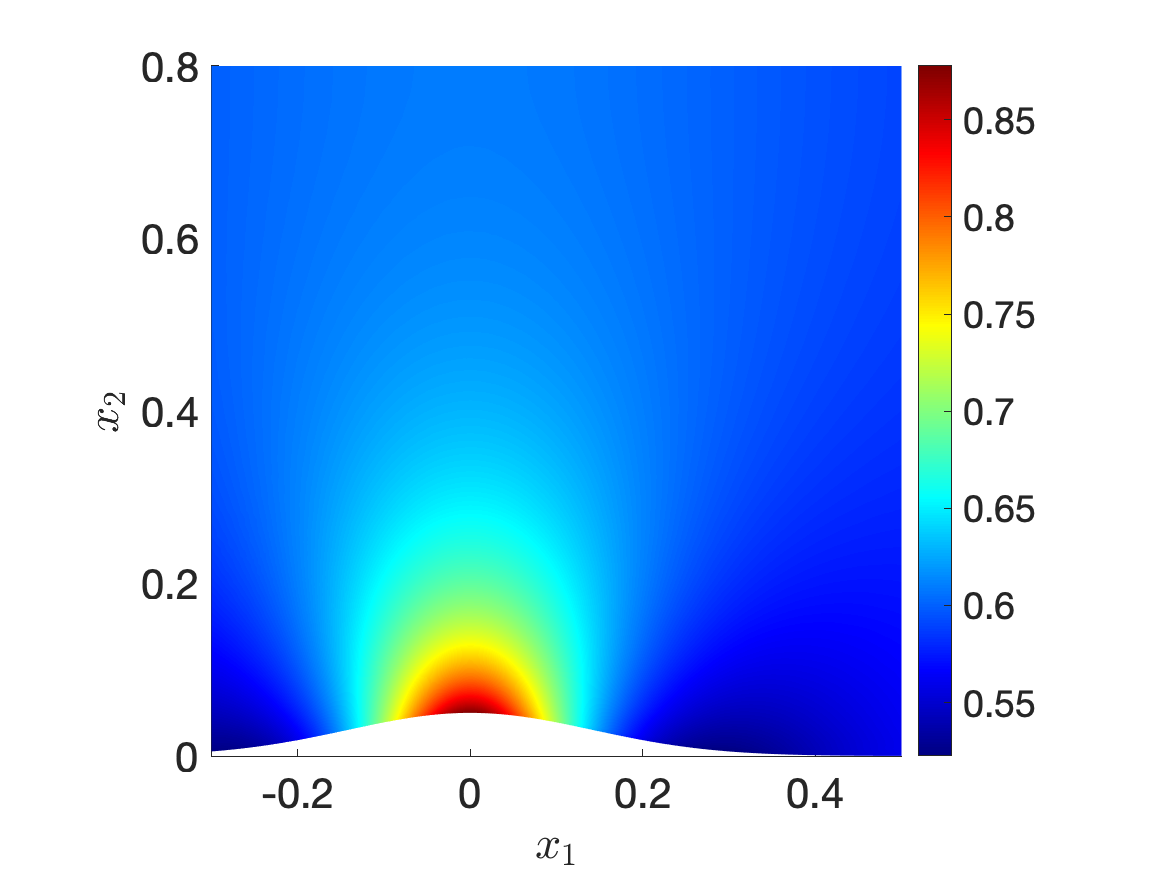}}
 ~~
   \subfloat[] 
{\includegraphics[width=0.4\textwidth]
 {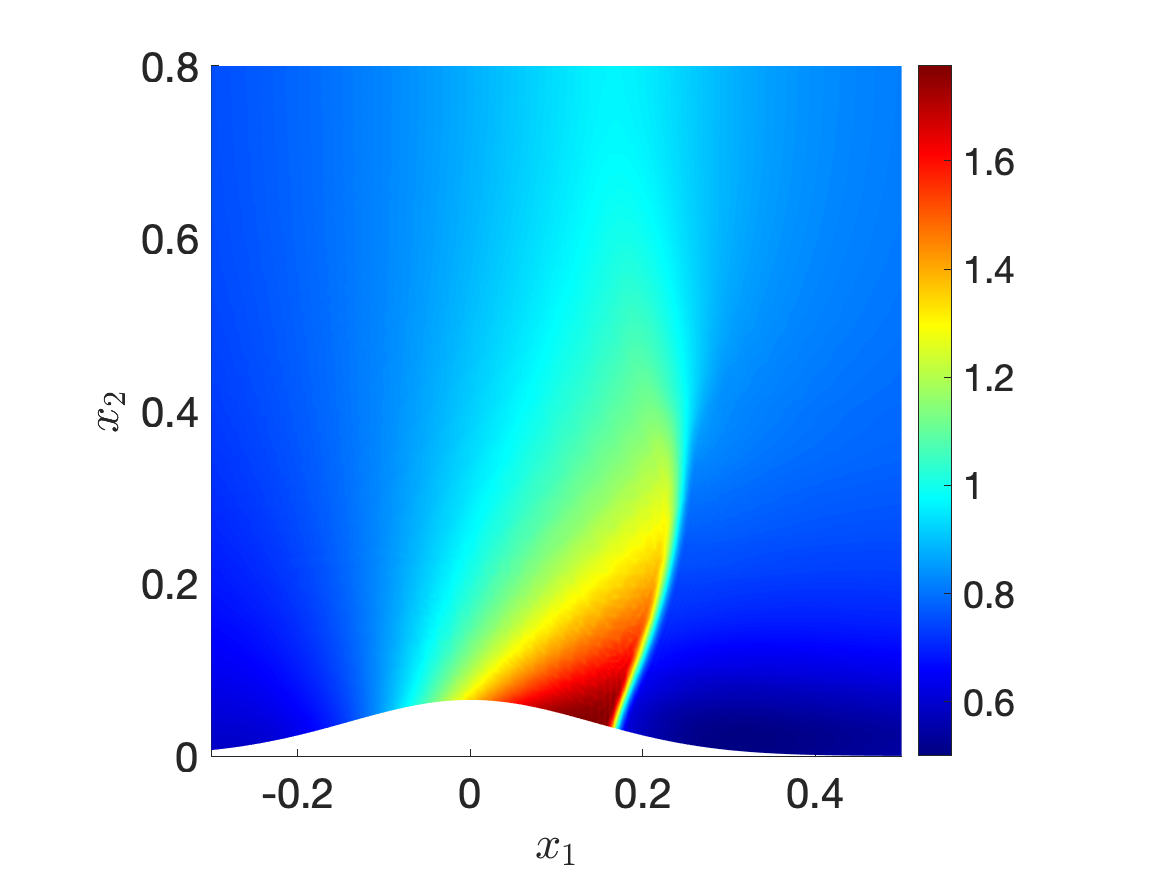}}
 
 \caption{inviscid flow over a Gaussian bump. 
 (a) computational domain.
 (b)-(c) visualization of the Mach field for $\mu=[0.05,0.58]$ and
 $\mu=[0.065,0.78]$.
}
 \label{fig:transbump_vis}
 \end{figure}
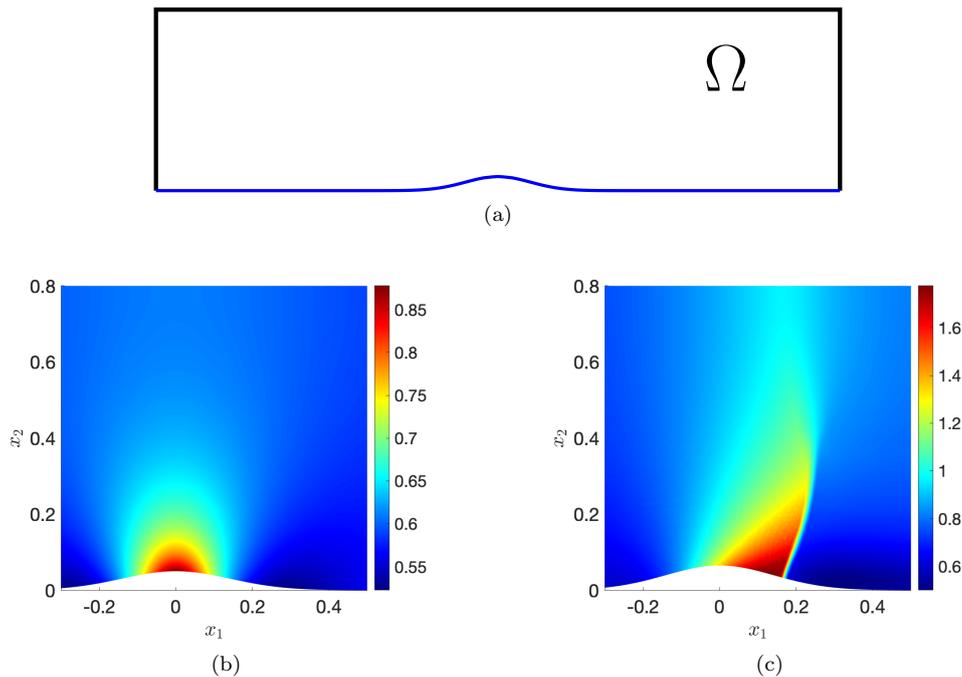

We also consider the two-dimensional inviscid flow past a Gaussian bump. We introduce  the domain $\Omega = \{
x\in (-1.5,1.5)\times (0,0.8): x_2 > h e^{-25 x_1^2}
\}$ where $h>0$ is a given parameter (cf. Figure \ref{fig:transbump_vis}(a)).
We consider the conservation law:
\begin{equation}
\label{eq:transbump_flow}
\nabla \cdot F(q^{\rm true}) = 0,
\quad {\rm where} \;\;
q^{\rm true} = \left[
\begin{array}{l}
 \rho^{\rm true} \\
 \rho^{\rm true} u^{\rm true} \\
 E^{\rm true} \\
\end{array}
\right],
\quad
F(q) = \left[
\begin{array}{l}
\rho u^{\top}  \\
 \rho u u^{\top}  \\
(E+p) u^{\top}  \\
\end{array}
\right],
\end{equation}
completed with wall boundary conditions on top and bottom boundaries, subsonic inlet condition (total temperature, total pressure and flow direction) at the left boundary and  subsonic outlet condition (static pressure) at the right boundary --- the symbol $(\cdot)^{\top} $ denotes the transposition operator. We express the free-stream field ${q}_{\infty}$ in terms of the Mach number ${\rm Ma}_{\infty}$,
$$
T_{\infty}=1, \;\;
p_{\infty}=\frac{1}{\gamma}, \;\;
\rho_{\infty}= 1, \;\;
u_{\infty}= {\rm Ma}_{\infty} e_1. \;\;
$$
Finally, we introduce the parameter vector $\mu=[h, {\rm Ma}_{\infty}]$ and the parameter region $\mathcal{P}=[0.05,0.065] \times [0.58,0.78]$.
Note that the computational  domain $\Omega$ depends on the geometric parameter $h$; therefore,  we should introduce a geometric mapping to recast the problem over a parameter-independent configuration. We here resort to a Gordon-Hall transformation; we refer to \cite[section 2]{ferrero2022registration} for the details.

Figures \ref{fig:transbump_vis}(b) and (c) show the behavior of the Mach field for two values of the parameters in $\mathcal{P}$: 
we observe that the flow is completely subsonic for moderate values of ${\rm Ma}_{\infty}$ and develops a normal shock over the bump   for   ${\rm Ma}_{\infty} \gtrsim 0.65-0.7$.

\section{Parametric mesh adaptation}
\label{sec:mesh_adaptation}
We consider the problem of determining an adapted mesh $\mathcal{T}_{\rm hf}$ of the domain $\Omega$ based on a set of snapshots $\widetilde{\mathcal{S}}_{\rm hf} = \{  \widetilde{q}_{\mu}^{\rm hf}: \mu\in \mathcal{P}_{\rm train} \}$ defined over a mesh $\mathcal{T}_{\rm hf}^0$.
We here pursue a metric-based MA
approach: first, we resort to the snapshots in $\widetilde{\mathcal{S}}_{\rm hf}$ to define a metric $\mathfrak{M}$; then, we resort to a MA toolkit  to devise the P1 mesh and subsequently we use standard FE routines to devise the high-order mesh.
In section \ref{sec:fundamentals_MA}, 
we review the key elements   of anisotropic mesh adaptation;
then,
in section \ref{sec:1D_meshadapt}
 we briefly comment on mesh adaptation for one-dimensional problems;
in section \ref{sec:isotropic_markrefine}, we 
discuss an isotropic mark-then-refine strategy employed for the two-dimensional problem of section \ref{sec:model_problems_bump};
finally, in section 
\ref{sec:anisotricMA},
 we review Hessian-based
 anisotropic mesh adaptation, and we discuss metric intersection for parametric problems.
  We refer to \cite{loseille2011continuous1,loseille2011continuous2} for a thorough introduction to metric-based MA.
 
\subsection{Fundamentals of anisotropic metric-based mesh adaptation}
\label{sec:fundamentals_MA}
A Riemannian metric  field is derived from the error estimate and prescribes the size and shape of the mesh elements.

\paragraph{Euclidean space.}

A scalar product is a symmetric positive definite (SDP) form, which can be represented by an SPD matrix $\mathfrak{M}$, which is dubbed \emph{metric tensor} or simply \emph{metric}. The scalar product is then written:
\begin{equation}
     \left(x,\, y\right)_{\mathfrak{M}}~=~ x^\top\,\mathfrak{M}\,y \quad\textrm{    where     }\quad x, y \in \mathbb{R}^d.
\end{equation}
%Note that when $\mathfrak{M}$ is the identity matrix, the scalar product is the canonical Euclidian dot product. 
A vector space with a scalar product is called an \emph{Euclidean space}.
The scalar product is associated with a distance that can be used to compute lengths in the Euclidean space :
\begin{equation}
  \ell_{\mathfrak{M}}(x,y) = \sqrt{(y-x)^\top\, \mathfrak{M}\,(y-x)}\,,
\end{equation}
from which we deduce classic geometrical quantities such as angles or volumes.
%the volume of an element $K$:
%\begin{equation}
% |K|_{\mathfrak{M}} = \sqrt{\det{\mathfrak{M}}}\,|K|_{\mathfrak{I}}\, .
%\end{equation}
	
The metric
$\mathfrak{M}$ is diagonalizable in an orthonormal basis:
\begin{equation}
  \mathfrak{M} = \mathcal{R}^\top\, \Lambda \, \mathcal{R}\,,
\end{equation}
where
  $\Lambda = \textrm{diag}\left(\lambda_1,\, \ldots,\, \lambda_d\right)$ is the {diagonal matrix of eigenvalues} 
   of $\mathfrak{M}$ and
  $\mathcal{R} = \left(\,\mathbf{r}_1 \,|\,\mathbf{r}_2\,|\, \ldots \, | \mathbf{r}_d\,\right)$ is the unitary matrix (\textit{i.e. } $\mathcal{R}^\top\, \mathcal{R} = \mathcal{I}_d )$ of eigenvectors  of $\mathfrak{M}$.
A metric tensor has an intuitive geometric representation:
the set of points that are at constant distance from  a point  $x\in \mathbb{R}^d$ is an ellipsoid centered in $x$ whose axes are aligned with the eigenvectors $\mathbf{r}_1,\ldots,\mathbf{r}_d$ of $\mathfrak{M}$. The set of points at   distance one  from a point (the unit ball of $\mathfrak{M}$), is an ellipsoid for which the sizes of the axes are $h_i = \lambda_i^{-\frac{1}{2}}$. 
In other words, in the context of anisotropic mesh adaptation, the eigenvectors of the metric tensor prescribe the orientation of the elements, while the eigenvalues prescribe the sizes in these directions. In an Euclidean metric space, these sizes and orientations are constant over all the domain; {for mesh adaptation, we want them to vary in space depending on the solution features.
This observation motivates the introduction of Riemannian metric spaces.
}

\paragraph{Riemannian metric space.}

We now  define a metric tensor field $\mathfrak{M}:\Omega\to \mathbb{R}^{d\times d}$ such that $\mathfrak{M}(x)$ is symmetric positive definite for all $x\in \Omega$. {There is no notion of global scalar product}; however,  we can extend the notion of distance. Given $x,y \in \Omega$, we define the distance
\begin{equation}
\label{eq:MA_distance}
\ell_{\mathfrak{M}}(x,y)
=\int_0^1 
\sqrt{ (y-x)^{\top}  \mathfrak{M}( (1-t) x + t y  ) (y-x)} \, dt\,,
\end{equation}
and, given a set $A\subset \Omega$, we define the volume:
\begin{equation}
\label{eq:MA_volume}
\big| A \big|_{\mathfrak{M}}
=
\int_A \sqrt{ {\rm det} (\mathfrak{M} (x))     } \, dx.
\end{equation}
Locally, the 
eigenvalues and the  eigenvectors of the 
 metric tensor $\mathfrak{M}(x)$ define  size  and orientations, respectively.

\paragraph{Unit mesh.} 

Given the mesh $\mathcal{T}_{\rm hf}$ of $\Omega$ with elements $\{\texttt{D}_k \}_{k=1}^{N_{\rm e}}$, we say that an element 
$\texttt{D}_k$ is a quasi-unit element with respect to  $\mathfrak{M}$ if the lengths \eqref{eq:MA_distance} of all its edges are approximately equal to one and its volume is approximately equal to $\frac{\sqrt{3}}{4}$ for $d=2$ and $\frac{\sqrt{2}}{12}$ for $d=3$. Similarly, we say that the mesh $\mathcal{T}_{\rm hf}$ is unit if all its elements are quasi-unit. Adapting a mesh with respect to $\mathfrak{M}$ comes to generating a mesh that is unit in that metric field.

\begin{remark}
\label{remark:mesh2metric}
Exploiting the geometric interpretation of the metric tensor, we can devise  a practical strategy to identify the  metric tensor $\mathfrak{M}$ associated to the triangle $\texttt{D}$:
{we shall use this strategy in section \ref{sec:isotropic_markrefine}.
}
We define the vertices $\{x_1^{\rm v}, x_2^{\rm v}, x_3^{\rm v}\}$ so that
 the longest (in Euclidean norm) edge is the one that connects  $x_1^{\rm v}$ and $x_2^{\rm v}$.
First, we set  $\lambda_1$ equal to the square of the inverse of ${\|  x_2^{\rm v} - x_1^{\rm v}\|_2}$,
and  $\lambda_2$ equal to the square of the inverse of the distance between the  vertex $x_3^{\rm v}$ and the  edge $\overline{ x_1^{\rm v}   x_2^{\rm v}}$; 
second, we set $\mathbf{n}_1 = \frac{x_2^{\rm v} - x_1^{\rm v}}{\|  x_2^{\rm v} - x_1^{\rm v}\|_2}$ and 
$\mathbf{n}_2 = [(\mathbf{n}_1)_2, - (\mathbf{n}_1)_1]^\top$;
finally, we define $\mathcal{R}  = [\mathbf{n}_1, \mathbf{n}_2 ]$, the diagonal matrix
$\Lambda = {\rm diag} \left(\lambda_1, \lambda_2 \right)$ and the metric $\mathfrak{M}  =\mathcal{R} ^{\top} \Lambda   \mathcal{R}$.
\end{remark}

\subsection{Mesh adaptation for one-dimensional problems}
 \label{sec:1D_meshadapt}
 For one-dimensional problems, we resort to the standard de Boor's algorithm (see, e.g., \cite[Chapter 2]{huang2010adaptive}): 
for consistency with   section \ref{sec:fundamentals_MA}, we present  
 the method in a slightly different formalism than the one of \cite{huang2010adaptive}.
 Given the  
metric $\mathfrak{M}:\Omega \to \mathbb{R}^+$, we define the  
 mesh density  function  $\mathfrak{d}: \Omega \to \mathbb{R}^+$:
\begin{equation}
   \label{eq:mesh_adapt_1D_density}
   \mathfrak{d}(x) = \sqrt{\mathfrak{M}(x)}\,, \quad x\in\Omega.
\end{equation} 
Our goal is to construct a (quasi-)unit mesh
with $N$ nodes   $\{x_i^{\rm hf}  \}_{i=1}^N$
 with respect to the metric $\mathfrak{M}$, that is
\begin{equation}
   \label{eq:mesh_adapt_1D_density_unit_mesh}
   \int_{x_1^{\rm hf}  }^{x_N^{\rm hf}   } \mathfrak{d}(x) \, dx = N
   \textrm{\quad and\quad}
   \int_{x_i^{\rm hf}}^{x_{i+1}^{\rm hf}} \mathfrak{d}(x) \, dx = 1\,,
   \quad
   i=1,\ldots,N-1\,.
\end{equation}
 
De Boor's algorithm constructs a unit mesh for an approximate metric.
 First, we introduce  an initial grid $ x_1^{\rm hf,0} \leq \ldots \leq x_{N_0}^{\rm hf,0}$ and the piecewise-constant approximation $\widehat{\mathfrak{d}}$ of $\mathfrak{d}$ such that,
$$
\widehat{\mathfrak{d}}(x) = 
\frac{1}{x_{i+1}^{\rm hf,0} - x_{i}^{\rm hf,0}}
\int_{ x_{i}^{\rm hf,0}  }^{ x_{i+1}^{\rm hf,0}   }
\mathfrak{d}(x) \, dx,
\quad
x\in \left(x_{i}^{\rm hf,0},  x_{i+1}^{\rm hf,0} \right),
\;\;
i=1,\ldots,N-1.
$$
Then, we find the unique set of points $\{x_i  \}_{i=1}^N$ that satisfies  \eqref{eq:mesh_adapt_1D_density_unit_mesh} for the mesh density function  $\widehat{\mathfrak{d}}$:
since $\widehat{\mathfrak{d}}$ is piecewise-constant, we can obtain an explicit expression for $\{x_i\}_{i=1}^N$. We refer to \cite{huang2010adaptive} for the explicit formula.
Note that several iterations of this algorithm are typically  required to obtain a good approximation of the equidistributed mesh for $\mathfrak{d}$.

 For the nozzle flow problem, we define the mesh density function
 $\mathfrak{d}$ based on the second-order derivative of the Mach number. Given the mapped fields
 $\{ \widetilde{q}_{\mu}^{\rm hf} : \mu\in \mathcal{P}_{\rm train} \}$, we define the corresponding Mach fields
  $\{ \widetilde{\rm Ma}_{\mu}^{\rm hf} : \mu\in \mathcal{P}_{\rm train} \}$ and the non-normalized density $\rho$
\begin{equation}
\label{eq:mesh_sensor_nozzle}
\rho(x)
=
\max_{\mu \in \mathcal{P}_{\rm train}}
\max \left\{
 \big| \partial_{xx}  \widetilde{\rm Ma}_{\mu}^{\rm hf} (x)\big|,
\, C_{\mu}
\right\}  
  \quad
  {\rm where}
  \quad
  C_{\mu}  = 10^{-2}
\sup_{x\in \Omega} \big| \partial_{xx}  \widetilde{\rm Ma}_{\mu}^{\rm hf} (x)\big|.
\end{equation}
Finally, we define the normalized density:
$$
 \mathfrak{d}(x) = \frac{N}{\int_\Omega \rho(x)\, dx} \rho(x).
$$
 The function $\rho$ (and thus $\mathfrak{d}$)  is well-defined in the interior of each element of the mesh $\mathcal{T}_{\rm hf}^0$; in our implementation, we rely on 
 the evaluation of the sensor \eqref{eq:mesh_sensor_nozzle} in the elements' quadrature points to define the piecewise-constant function that is used by the de Boor's algorithm.

\subsection{Isotropic mark-then-refine mesh adaptation}
\label{sec:isotropic_markrefine}
For the Euler equations, the total enthalpy $H_{\rm tot}$
(cf. \eqref{eq:conservation_law})
is constant and can be computed exactly from the boundary conditions for any $\mu\in \mathcal{P}$.
Given the set of snapshots $\widetilde{\mathcal{S}}_{\rm hf}$ defined over the mesh $\mathcal{T}_{\rm hf}^0$, we first compute the average error in total enthalpy
\begin{subequations}
\label{eq:total_enthalpy_error_indicator}
\begin{equation}
\label{eq:total_enthalpy_error_indicator_a}
\eta_{\mu,k} = \frac{1}{| \texttt{D}_k^0  |} \int_{ \texttt{D}_k^0  } \left( H_{\rm tot, \mu}^{\rm hf} -  H_{\rm tot, \mu}^{\rm true} \right)^2 \, dx,
\quad
k=1,\ldots,N_{\rm e}^0,
\end{equation}
with $H_{\rm tot, \mu}^{\rm true}:=H_{\rm tot}(q_{\mu}^{\rm true})$ and
$H_{\rm tot, \mu}^{\rm hf}:=H_{\rm tot}(q_{\mu}^{\rm hf})$, 
and the maximum over the training set of parameters
\begin{equation}
\label{eq:total_enthalpy_error_indicator_b}
\eta_{k}^{\rm max}  =
\max_{\mu\in \mathcal{P}_{\rm train}}
\eta_{\mu,k} ,
\quad
k=1,\ldots,N_{\rm e}^0.
\end{equation}
\end{subequations}
Second, 
given $\gamma_{\rm ref}\in (0,1)$,
we mark the $\gamma_{\rm ref} \cdot 
N_{\rm e}^0$ elements that maximize $\{ \eta_{k}^{\rm max}    \}_k$; we denote by 
$\texttt{I}_{\rm max}\subset \{1,\ldots,N_{\rm e}^0\}$ the indices of the marked elements.
Third, we extract the metric $\{  \mathfrak{M}_k^0 \}_{k=1}^{N_{\rm e}^0}$ from the mesh 
$\mathcal{T}_{\rm hf}^0$ using the strategy in Remark \ref{remark:mesh2metric} and we define the new metric as follows:
\begin{equation}
\label{eq:new_metric_iso}
\mathfrak{M}_k = \mathfrak{M}_k^0 \;\;{\rm if} \; k\notin
\texttt{I}_{\rm max} , \quad
  \mathfrak{M}_k = 4 \mathfrak{M}_k^0 \;\;{\rm if} \; k\in
\texttt{I}_{\rm max}.
\end{equation}
Multiplication by four in \eqref{eq:new_metric_iso} corresponds to an isotropic reduction of the local mesh density by a factor of two.
Fourth, we define the metric $\mathfrak{M}$ in the vertices of the mesh $\mathcal{T}_{\rm hf}^0$ using the simple average:
$$
\mathfrak{M}_i = \frac{1}{\# {\rm Neigh}_i } \sum_{k\in {\rm Neigh}_i} \mathfrak{M}_k,
\quad
i=1,\ldots,N_{\rm v}^0,
$$
where ${\rm Neigh}_i \subset \{1,\ldots,N_{\rm e}^0\}$ contains the indices of the elements that contain the $i$-th vertex of $\mathcal{T}_{\rm hf}^0$. Finally, we apply a  mesh adaptation toolkit %\texttt{mmg2d} 
to generate the new mesh.

The choice of $\gamma_{\rm ref}$ regulates how quickly we  increase the size of the HF mesh. 
Since the mesh $\mathcal{T}_{\rm hf}$ should be accurate for the mapped manifold $\{ q_{\mu}^{\rm true} \circ \Phi_{\mu} : \mu\in \mathcal{P} \}$ where the mapping $\Phi$ changes at each outer-loop iteration, it is not worth to refine the mesh if the mapping $\Phi$ is excessively inaccurate. On the other hand, modest values of 
$\gamma_{\rm ref}$ might require a large number of   iterations 
in Algorithm \ref{alg:offline_training}. A thorough investigation of the  choice of $\gamma_{\rm ref}$ on performance is beyond the scope of the present work;  in the numerical experiments, we mark $10\%$ of the elements
(i.e., $\gamma_{\rm ref}=10\%$)
of the mesh at each iteration.

We observe that the computation of \eqref{eq:total_enthalpy_error_indicator} for all $\mu\in \mathcal{P}_{\rm train}$ might be expensive: in the numerical experiments, we hence run a strong greedy algorithm (cf. Appendix \ref{sec:appendix_greedy}) to identify the most relevant parameters in $\mathcal{P}_{\rm train}$.
We also notice that the choice of the error indicator
\eqref{eq:total_enthalpy_error_indicator}
 is specific to the Euler equations: several alternatives have been considered in the literature such that the $\texttt{p}+1$ residual or goal-oriented error estimates (e.g., \cite{yano2012optimization}). 
Finally, we observe that at each iteration of   
Algorithm \ref{alg:offline_training} we generate a new mesh that is independent of the previous meshes: the advantage of this choice is that we allow the mesh adaptation toolkit to perform smoothing operations that ensure well-behaved meshes; clearly, the disadvantage is that we need to perform   mesh interpolation 
between unstructured meshes
at several steps of our training phase.

\subsection{Anisotropic mesh adaptation}
\label{sec:anisotricMA}

Anisotropic mesh adaptation refers to a class of methods where both the elements sizes and orientations are optimized with respect to an error estimate. As explained in section \ref{sec:fundamentals_MA}, this can be achieved by generating a unit mesh with respect to a prescribed metric field. In this work, we use the multiscale metric defined in \cite{loseille2011continuous2}:

\begin{equation}
\label{eq:hessian_based_metric}
\mathfrak{M}_{s}(x) 
=
\left(
\frac{N}{\int_{\Omega}  
\left( {\rm det}( | H_s(\bar{x}) | )\mathrm{d}\bar{x} 
\right)^\frac{p}{2p + d}  }
\right)^{2/d}
\,
{\rm det}( | H_s(x) | )  ^{\frac{-1}{(2p + d)}}
\,
| H_s(x) |,
\quad
\forall x \in \Omega\,,
\end{equation}
where $s: \Omega\to \mathbb{R}$ is a scalar field driving the adaptation, $H_s:\Omega\to \mathbb{R}^{d\times d}$ is the Hessian matrix of $s$. Since  $H_s$ is real-valued and  symmetric, it can  be diagonalized:
$H_s =  \mathcal{R}_s^{\top} \Lambda_s \mathcal{R}_s$ 
where $\mathcal{R}_s$ is the orthogonal matrix of eigenvectors and 
$\Lambda_s = {\rm diag}(\lambda_{s,1},\lambda_{s,2})$ is the diagonal matrix of eigenvalues. 
We further define $|H_s|$  by taking the absolute value of the eigenvalues:
$| H_s | := 
\mathcal{R}_s^{\top}  {\rm diag}(| \lambda_{s,1} |,| \lambda_{s,2} |) \mathcal{R}_s$.
The constant $p$ is associated with the $L^p$ norm 
 that is used to derive 
the error estimate and ultimately  the metric.
On the other hand, $N$ is the \emph{continuous complexity} and is directly related to the target number of vertices of the mesh.
 
In this work, we do not only want to adapt to one field, but to as many fields as we have snapshots. 
Instead of generating several meshes and taking the supermesh of several adapted meshes, we choose to construct one metric based on the available fields and generate only one adapted mesh. 
A fairly common strategy is used to construct that unique metric: the multiscale metric (from Eq.~\eqref{eq:hessian_based_metric}) is computed for each snapshot, and those metrics are then combined vertex-wise using metric intersection :
\begin{equation}
\label{eq:parametric-multiscale-metric}
	\mathfrak{M}(x) = \bigcap_{k =1}^{n_\textrm{train}} \mathfrak{M}_{s^k}(x)\,,
\end{equation}
where $n_\textrm{train}$ is the number of snapshots in $\mathcal{P}_\textrm{train}$ and $\mathfrak{M}_{s^k}$ is the multiscale metric field computed for the chosen scalar field of the $k$-th snapshot of the dataset.

A procedure to intersect two metrics is proposed in~\cite{barral2015thesis}. Geometrically, we have seen that the unit ball of a metric is an ellipsoid. Intersecting two metrics is equivalent to finding the largest ellipsoid included in the two corresponding ellipsoids. Note that this method has limitations, which we aim to solve in future work. First, metric intersection is not associative, which means that depending on the order in which the intersections are carried out, the final result will slightly change, although it is not known to change significantly. Second, if each metric field associated with each snapshot would result in a mesh with a certain prescribed number of vertices, once the metrics are intersected, we do not control the new prescribed number of vertices. To avoid this issue, we could intersect the Hessians and then apply the multiscale normalization. However, the geometric interpretation of intersecting quantities with potentially very different orders of magnitudes is  unclear.

\subsection{Practical considerations}
Good approximations of the   gradient  and   the Hessian 
of the solution field
are key in the implementation of the metric. If the order of the solver is high enough (strictly greater than one), those quantities are computed exactly in each element. Otherwise, we perform $L^2$ projections based  on  Clement's interpolation~\cite{clement1975approximation} to reconstruct a $P^1$ gradient (resp. Hessian) from a $P^0$ gradient (resp. Hessian).
Metric~\eqref{eq:hessian_based_metric} depends on the choice of scalar field $s$, and parameter $p$.
We here use the Mach number field as sensor to drive the adaptation.
In the numerical simulations, we set $p=1$, which is found to better capture small-scale features of the solution field.

The task of generating the adapted meshes is left to the  anisotropic remesher \texttt{mmg2d} (\cite{arpaia2022h,dapogny2014three}). It takes 
as input a mesh and a metric field defined on the mesh, and returns 
a  unit  mesh   for the given metric. To do so, it performs iteratively a series of local mesh operations : vertex addition, removal, smoothing and topology changes. In this work, the default parameters of \texttt{mmg2d} are used. The eigenvalues of   metric~\eqref{eq:hessian_based_metric} are previously truncated to avoid excessively small edge sizes.

Mesh adaptation is a non-linear problem, where the convergence of the mesh/solution couple has to be considered: a better mesh gives a better solution, which in turns gives a better mesh, etc. 
This is addressed in Algorithm~\ref{alg:offline_training}: 
during each outer loop iteration, we indeed generate snapshots (cf. Line 3) associated to the current mesh and we use them to adapt the mesh (cf. Line 4).
We hence expect to converge to a final mesh --- and a reduced-order approximation --- that is accurate for all parameters in $\mathcal{P}$.
%Initially, for each high-fidelity solution, several iterations of solution and adaptation are performed. Then, the quality of the common adapted mesh is improved by adding new snapshots in the metric.

\section{Registration}
\label{sec:registration}
The {second} ingredient of our method is a registration algorithm that is designed to track coherent structures of the solution field, to facilitate the tasks of mesh adaptation and linear-subspace model reduction; the algorithm takes as input
(i)  a set of snapshots $\{ q_{\mu}^{\rm hf} : \mu\in \mathcal{P}_{\rm train} \}$ and
(ii) a mesh $\mathcal{T}_{\rm hf}$ of the domain $\Omega$, 
and returns a parameterized map $\Phi: \Omega \times \mathcal{P} \to \Omega$ such that
$\Phi_{\mu}(\mathcal{T}_{\rm hf})$ is a proper mesh of $\Omega$ for all $\mu\in \mathcal{P}$.
The development and the analysis of registration methods for MOR remains a challenging task that requires many advances; in this work, we briefly summarize the procedure employed in the numerical experiments
 and we refer to a future work for a thorough discussion on registration methods.

\subsection{Spectral maps}
\label{sec:spectral_maps}
Given the family of domains $\{ \Omega_{\mu} : \mu\in \mathcal{P} \}$, we define the ``reference''  domain ${\Omega}_{\rm p}$ and the geometric map 
$\Psi^{\rm geo}:  {\Omega}_{\rm p} \times \mathcal{P} \to \mathbb{R}^d$ such that 
$\Psi_{\mu}^{\rm geo}(  {\Omega}_{\rm p} ) = \Omega_{\mu}$ for all $\mu\in \mathcal{P}$. 
We denote by ${\mathbf{n}}_{\rm p}$ the outward normal to $\partial \Omega_{\rm p}$ and we define the space of tensorized polynomials $\mathbb{Q}_J$ of degree at most $J$ in each variable.
For the nozzle problem, we consider 
$\Psi^{\rm geo} = \texttt{id}$ --- where $\texttt{id}(x)=x$ is the identity map ---
and ${\Omega}_{\rm p} = \Omega$;
for the transonic bump problem, we consider a Gordon-Hall map
(cf. \cite[section 2]{ferrero2022registration})  and 
 ${\Omega}_{\rm p} = (0,1)^2$;
in the latter, we introduce the reference parameter $\bar{\mu}\in \mathcal{P}$ 
and we define
$\Omega:=\Omega_{\bar{\mu}}$. 
 Then, we consider mappings of the form
\begin{subequations}
\label{eq:map_gen}
\begin{equation}
\label{eq:map_gen_a}
\texttt{N}_{\mu}(\mathbf{a} )
\, = \,
\Psi_{\mu}^{\rm geo} \circ 
\texttt{N}_{\rm p}(\mathbf{a})
 \circ  \Lambda_{\bar{\mu}}^{\rm geo},
\quad
\texttt{N}_{\rm p}(\mathbf{a})= \texttt{id} +  \sum_{i=1}^m (\mathbf{a})_i \varphi_i,
\end{equation}
where $\Lambda_{\bar{\mu}}^{\rm geo} := \left( \Psi_{\bar{\mu}}^{\rm geo} \right)^{-1}$,  $\{  \varphi_i \}_{i=1}^m$ spans the space $\mathcal{U}_{\rm p}$ of tensorized polynomials such that
\begin{equation}
\label{eq:map_gen_b}
\mathcal{U}_{{\rm p}} = {\rm span}   \{  \varphi_i \}_{i=1}^m 
\subset 
\mathcal{U}_{{\rm hf,p}} =
\left\{
\varphi \in [\mathbb{Q}_J]^d: \;
\varphi \cdot \mathbf{n}_{\rm p} |_{\partial  {\Omega}_{\rm p} } = 0
\right\}.
\end{equation}
\end{subequations}
We equip $\mathcal{U}_{{\rm hf,p}}$ with the inner product
\begin{equation}
\label{eq:norm_mapping_space}
( \varphi, \psi )_{H^2({\Omega}_{\rm p})}^2 : =
 \int_{{\Omega}_{\rm p}} \; 
 \left(  \sum_{i,j,k=1}^d  \partial_{j,k} (\varphi)_i  \cdot \partial_{j,k} (\psi)_i  \; + \;
 \varphi \cdot \psi \right) \, d x,
\end{equation} 
and we assume that  $\{  \varphi_i \}_{i=1}^m$ is an orthonormal basis of $\mathcal{U}_{{\rm p}} $.
We observe that, if $\texttt{N}_{\rm p}(\mathbf{a})$ is a bijection from $\Omega_{\rm p}$ in itself, $\texttt{N}_{\mu}(\mathbf{a} )$ is a bijection from $\Omega$ to $\Omega_{\mu}$, for all $\mu\in \mathcal{P}$. 

We denote by $\mathcal{B}_{\mu}$ the space of diffeomorphisms from $\Omega$ to $\Omega_{\mu}$; exploiting the analysis in \cite{taddei2020registration,taddei2021registration}, we can prove that 
(i)
for any $\mu\in \mathcal{P}$ the model class $\texttt{N}$ \eqref{eq:map_gen} 
is dense in a meaningful subspace  of $\mathcal{B}_{\mu}$ and 
(ii)
the set of admissible maps
$A_{\rm bj, \mu} =\{ \mathbf{a} \in \mathbb{R}^m : \texttt{N}_{\mu}(\mathbf{a} ) \in     
\mathcal{B}_{\mu} 
\}$ has a non-empty interior for any choice of $\mathcal{U}_{{\rm p}}$ in \eqref{eq:map_gen_b}.  The latter is extremely important for the numerical robustness of registration methods. 
We remark that the two results are currently restricted to domains that are diffeomorphic to the unit hyper-cube $\Omega_{\rm p} = (0,1)^d$. The extension of these results to a broader class of domains is the subject of ongoing research.

\subsection{Optimization-based registration}
\label{sec:optreg}

Given the training set of snapshots
$\{q_{\mu}^{\rm hf} : \mu\in \mathcal{P}_{\rm train} \}$,
we determine the mapping coefficients $\widehat{\mathbf{a}}_{\mu}$ (cf. \eqref{eq:lagrangian}) by solving the optimization problem:
\begin{equation}
\label{eq:optreg_statement}
\min_{\mathbf{a}\in \mathbb{R}^m}
\mathfrak{f}_{\mu}^{\rm obj}(\mathbf{a}) := 
\mathfrak{f}_{\mu}^{\rm tg}(\mathbf{a}) \, + \,
\xi \left(
\big|
\texttt{N}_{\rm p}(\mathbf{a})
\big|_{H^2(\Omega_{\rm p})}^2
\, + \,
\mathfrak{f}_{\rm msh} (   
\texttt{N}_{\mu}(\mathbf{a})
)
\,+ \,
\mathfrak{f}_{\rm jac}(  
\texttt{N}_{\rm p}(\mathbf{a})
)
\right).
\end{equation}
Here, $\mathfrak{f}_{\mu}^{\rm tg}$ denotes the target (or proximity) function that measures the degree of similarity between the available estimate of the solution field $q_{\mu}^{\rm true}$ and a suitable template solution or template reduced space, while the terms multiplied by the weighting parameter $\xi>0$ are regularization terms that promote the smoothness of the map and ensure bijectivity.
In more detail, 
$\big| \cdot  \big|_{H^2(\Omega_{\rm p})}^2 = 
\int_{\Omega_{\rm p}} \; \left( 
\sum_{i,j,k=1}^d
\partial_{j,k} (\cdot) _i^2 \right) \, dx$ is the $H^2$ seminorm; 
$\mathfrak{f}_{\rm msh}$ controls the quality of the deformed mesh
%$\texttt{N}_{\mu}(\widehat{\mathbf{a}}_{\mu})(  \mathcal{T}_{\rm hf} )$
(cf. \cite{zahr2020implicit}),  
\begin{equation}
\label{eq:f_msh}
\mathfrak{f}_{\rm msh}({\Phi}   )
=
\frac{1}{|\Omega|}
\sum_{k=1}^{N_{\rm e}^{\rm pb}  } \; \; 
|  \texttt{D}_k^{\rm pb}  |
\int_{  \widehat{\texttt{D}}  } 
{\rm exp} 
\left(
\frac{q_{k}^{\rm msh}(\Phi)}{q_{k}^{\rm msh}(\texttt{id})}
\, - \, 
\kappa_{\rm msh}
\right)
 \; dx,
 \quad
 q_{k}^{\rm msh}(\Phi) : = 
 \frac{1}{d^2}
\left(
\frac{ \|  \nabla \Psi_{\Phi,k}^{\rm hf} \|_{\rm F}^2  }{( {\rm det} ( \nabla  \Psi_{\Phi,k}^{\rm hf}   )    )_+^{2/d}}\right)^2,
\end{equation}
where
$\{  \Psi_{\Phi,k}^{\rm hf} \}_k$ are the elemental maps
\eqref{eq:psi_mapping} associated to the deformed mesh 
${\Phi}  (  \mathcal{T}_{\rm hf} )$; and
$\mathfrak{f}_{\rm jac}$ is designed to ensure that the selected map is non-singular,
\begin{equation}
\label{eq:f_jac}
\mathfrak{f}_{\rm jac}( \Phi_{\rm p} )
=
\frac{1}{|\Omega_{\rm p}|}
\int_{\Omega_{\rm p}}
{\rm exp} \left(
\frac{\epsilon - {\rm det} (\nabla \Phi_{\rm p})}{C_{\rm exp}} 
\right)
\, dx,
\quad
{\rm with} \;\;
\epsilon \in (0,1), \;\;
C_{\rm exp} \ll \epsilon.
\end{equation}
Note that \eqref{eq:f_msh} and \eqref{eq:f_jac} depend on several hyper-parameters: in the numerical experiments, we consider
$$
\epsilon = 0.1, \quad
C_{\rm exp} = 0.025 \epsilon, \quad
\kappa_{\rm msh} = 10, \quad
\xi=10^{-3}.
$$

We observe that $ q_{k}^{\rm msh}(\Phi) \equiv  1$ for $d=1$ dimensional problems that are discretized using linear elements: we hence omit the mesh regularization term for the nozzle problem.
Furthermore, we empirically found that the regularization \eqref{eq:f_jac} is not strictly needed for two-dimensional problems based on discretize-then-map treatment of geometry parameterizations (cf. section \ref{sec:pMOR}): in the numerical experiments we hence omit the regularization \eqref{eq:f_jac}  for the transonic bump test case.
In the remainder of this section, we discuss the choice of the target function for the two model problems considered in the numerical section.

\subsubsection*{Target function for the nozzle problem}
Given the snapshot $q_{\mu}^{\rm hf}$, we compute the Mach field 
${\rm Ma}_{\mu}^{\rm hf}$ and we estimate the maximum of its derivative $x_{\mu}^{\star}$; then, we consider the target 
\begin{equation}
\label{eq:target_function_nozzle}
\mathfrak{f}_{\mu}^{\rm tg}(\mathbf{a} )
=
\big|
\texttt{N}  (  x_{\bar{\mu}}^{\star} ; \mathbf{a}) 
 \, - \,
x_{{\mu}}^{\star}
\big|^2.
\end{equation}
In the numerical experiments, we estimate $x_{{\mu}}^{\star}$ using the formula
\begin{equation}
\label{eq:max_estimate}
x_{\mu}^{\star}
=
\frac{1}{\# \texttt{I}_+} \sum_{i\in \texttt{I}_+}
x_i^{\rm hf,qd},
\quad
\texttt{I}_+ = \left\{
i \in \{1,\ldots,N_{\rm hf,q}\} :
|
\partial_x  {\rm Ma}_{\mu}^{\rm hf}(  x_i^{\rm hf,qd}   )
|
> \delta \max_j | \partial_x   {\rm Ma}_{\mu}^{\rm hf}(  x_i^{\rm hf,qd}   )  |
\right\},
\end{equation}
where $\{ x_i^{\rm hf,qd} \}_{i=1}^{N_{\rm hf,q}}$ are the quadrature points of the FE mesh and $\delta>0$ is a threshold that is set equal to $0.5$.
We observe that the definitions of \eqref{eq:target_function_nozzle} and \eqref{eq:max_estimate} exploit the knowledge that the solution exhibits a single discontinuity in $\Omega$; we refer to 
\cite{iollo2022mapping} and \cite{taddei2022optimization} for a generalization to a more general setting.

\subsubsection*{Target function for the transonic bump problem}
We consider the target 
\begin{equation}
\label{eq:target_transonic_bump}
\mathfrak{f}_{\mu}^{\rm tg}( \mathbf{a}  )
=
\min_{\nu\in \mathcal{S}_n}
\frac{1}{|\Omega_{\rm p}|}
\,
 \int_{\Omega_{\rm p}} \;
\big| s_{\mu}^{\rm hf} \circ \texttt{N}_{\rm p}(\mathbf{a}) - \nu \big|^2 \,  {\rm det}(\nabla \Psi_{\bar{\mu}})  dx
\, + \,
\|   \texttt{N}_{\mu} (  x_{\bar{\mu}}^{\star} ; \mathbf{a}) \, - \,
x_{{\mu}}^{\star}  \|_2^2,
\end{equation}
where
$s_{\mu}^{\rm hf} = {\rm Ma}_{\mu}^{\rm hf} \circ  \Psi_{\mu}^{\rm geo}$,
and  $x_{{\mu}}^{\star}$  is equal to the maximum of the Mach number over the bump if the flow is subsonic, and equal to the maximum of the tangential derivative of the Mach number --- which is  practically  estimated using \eqref{eq:max_estimate} --- if the flow is transonic.
We observe that the evaluation of the target \eqref{eq:target_transonic_bump} requires the evaluation of the field $s_{\mu}^{\rm hf}$ in arbitrary points of $\Omega_{\rm p}$: it is thus important to define $s_{\mu}^{\rm hf}$  over a structured grid.

For completeness, we comment on the choice of the first term in \eqref{eq:target_transonic_bump}. Given the reduced space $\widetilde{\mathcal{S}}_n \subset L^2(\Omega_{\mu})$, the goal of registration is to find a mapping $\Phi$ such that
$$
\min_{\nu\in  \widetilde{\mathcal{S}}_n}  \;
\int_{ \Omega_{\mu}  } 
\, \big| {\rm Ma}_{\mu}^{\rm hf} \circ \Phi- \nu \big|^2 \,  dx,
$$
where the choice to consider the Mach number as registration sensor is justified by the observation that it is a scalar quantity that exhibits relevant features (shocks, contact discontinuities) of the full field $q_{\mu}^{\rm true}$.
Exploiting the expression of $\Phi$, 
$\Phi = \Psi_{\mu}^{\rm geo} \circ  \Phi_{\rm p}  \circ  \Lambda_{\bar{\mu}}^{\rm geo}$,
and the change of variable $x=\Psi_{\bar{\mu}}^{\rm geo}(\xi)$, we find
$$
\min_{\nu\in  \widetilde{\mathcal{S}}_n}  \;
\int_{ \Omega_{\mu}  } 
\, \big| s_{\mu}^{\rm hf} \circ \Phi- \nu \big|^2 \, 
 {\rm det}(\nabla \Psi_{\bar{\mu}})  dx
\quad
{\rm where} \;\; 
{\mathcal{S}}_n
=
\left\{
\nu \circ \Psi_{\mu}^{\rm geo} \, : \, 
\nu \in \widetilde{\mathcal{S}}_n
\right\},
\quad
s_{\mu}^{\rm hf}  = {\rm Ma}_{\mu}^{\rm hf} \circ \Psi_{\mu}^{\rm geo}.
$$
The space $\mathcal{S}_n \subset L^2(\Omega_{\rm p})$ in \eqref{eq:target_transonic_bump} is dubbed \emph{template space} and is built using the greedy procedure proposed in
 \cite{taddei2021space}.

\subsection{Parametric registration}
\label{sec:param_reg}
We combine the optimization statement discussed in the previous section with the greedy algorithm proposed in \cite{taddei2021space} for the adaptive construction of the template space $\mathcal{S}_n$ in \eqref{eq:target_transonic_bump}, and a standard regression procedure to obtain the parametric mapping $\Phi$
--- for completeness, we report the greedy method in
Appendix \ref{sec:appendix_greedy}. For the nozzle problem, the greedy procedure is not necessary: in this case we simply rely on \cite[Algorithm 1]{taddei2020registration}.
In both cases, the cost of the procedure is dominated by the solution to the optimization statement \eqref{eq:optreg_statement} for all $\mu\in \mathcal{P}_{\rm train}$,
\begin{equation}
\label{eq:optimization_all_params}
\widehat{\mathbf{a}}_{\mu} \in {\rm arg} \min_{\mathbf{a}\in \mathbb{R}^m} \mathfrak{f}_{\mu}^{\rm obj}(\mathbf{a}),
\quad
\mu\in \mathcal{P}_{\rm train}
\end{equation} 
for the first iteration of the algorithm --- which 
corresponds to the choice $\mathcal{S}_1 ={\rm span}  \{ s_{\bar{\mu}}^{\rm hf}  \}$ for the transonic bump test case.

We rely on the Matlab function \texttt{fminunc} which implements a quasi-Newton method; since the problem is non-convex, 
the choice of the initial condition for the optimizer is critical to achieve accurate performance.
Towards this end, following \cite{taddei2020registration},
we first reorder the parameters in $\mathcal{P}_{\rm train}$ so that
$\mu^{(1)}
=
{\rm arg}\min_{\mu \in \mathcal{P}_{\rm train}}
\| \mu - \bar{\mu} \|_2$ and
$$
\mu^{(k)} = {\rm arg} \min_{\mu \in \mathcal{P}_{\rm train} 
\setminus \{  \mu^{(i)} \}_{i=1}^{k-1}
 }
 \left(
 \min_{\mu' \in   \{  \mu^{(i)} \}_{i=1}^{k-1}}
 \| \mu - \mu'  \|_2
 \right),
 \quad
k=2,\ldots,n_{\rm train};
$$
then, we choose the initial condition as follows:
$$
\mathbf{a}_{\mu^{(1)}}^0 = 0,
\quad
\mathbf{a}_{\mu^{(k)}}^0
=
\widehat{\mathbf{a}}_{\mu^{({\rm ne}_k)}},
\;\;
{\rm with} \;\;
{\rm ne}_k
={\rm arg}\, \min_{j=1,\ldots,k-1} 
 \| \mu^{(j)} - \mu^{(k)} \|_2,
 \quad k=2,\ldots,n_{\rm train}.
$$
We observe that this choice of the initial condition prevents the parallelization of the registration procedure.

\begin{remark}
\label{remark:implementation_registration}
In the numerical experiments for the two-dimensional test case, we consider polynomials of degree $J=10$ and we rely on a P1 $61\times 21$ Cartesian FE grid of the unit square to represent the sensors
$\mu\mapsto s_{\mu}^{\rm hf}$. If we denote by 
$\{  x_j^{\rm hf,r} \}_{j=1}^{N_{\rm nd,r}}$ the nodes of the mesh on $\Omega_{\rm p}$, computation of $s_{\mu}^{\rm hf}$ requires the interpolation of the FE field ${\rm Ma}_{\mu}^{\rm hf}$ in the points
$\{ \Psi_{\mu}^{\rm geo} (  x_j^{\rm hf,r} ) \}_{j=1}^{N_{\rm nd,r}}$. To ensure that the objective function is sufficiently smooth for gradient-based optimization, we post-treat the sensor by applying  a low-pass filter (moving average) in each spatial direction.
\end{remark}

\section{Linear-subspace projection-based model order reduction}
\label{sec:pMOR}
In this section, we  present the projection-based MOR procedure employed to estimate the mapped field $\widetilde{q}_{\mu}^{\rm true} : = {q}_{\mu}^{\rm true} \circ \Phi_{\mu}$. As anticipated in the introduction, we seek approximations of the form
\begin{subequations}
\label{eq:MOR_intro}
\begin{equation}
\label{eq:MOR_intro_a}
\widetilde{q}_{\mu} = \texttt{Z} \widehat{\boldsymbol{\alpha}}_{\mu}
\quad
{\rm with} \;\;
\widehat{\boldsymbol{\alpha}}_{\mu} \in  {\rm arg} \min_{  \boldsymbol{\alpha} \in \mathbb{R}^n  }
\max_{\psi \in \widehat{\mathcal{Y}}} \;
\frac{\mathfrak{R}_{\mu}^{\rm eq}(
\texttt{Z}  \boldsymbol{\alpha}, \psi)}{\vertiii{\psi}}
\end{equation}
where $\texttt{Z}:\mathbb{R}^n \to \mathcal{X}_{\rm hf}$ is a suitable linear operator,
\begin{equation} 
\label{eq:MOR_intro_b}
\mathfrak{R}_{\mu}^{\rm eq}(
q, v)
=
\sum_{k=1}^{N_{\rm e}} 
\rho_{k}^{\rm eq,e}
\,
r_{k,\mu}^{\rm e}(E_k q, E_k v  ) \; + \;
\sum_{j=1}^{N_{\rm f}} 
\rho_{k}^{\rm eq,f}
r_{j,\mu}^{\rm f}(
E_j^{+} q, E_j^{-} q,
 E_j^{+} v,  E_j^- v  ),
 \quad
 \forall \, q, v\in  \mathcal{X}_{\rm hf},
\end{equation}
\end{subequations}
is a weighted residual that depends on the sparse weights
$\boldsymbol{\rho}^{\rm eq,e}\in \mathbb{R}^{N_{\rm e}}$ and
$\boldsymbol{\rho}^{\rm eq,f}\in \mathbb{R}^{N_{\rm f}}$, 
$\widehat{\mathcal{Y}} \subset \mathcal{X}_{\rm hf}$ is a $m$-dimensional linear space with $m\geq n$, and $\vertiii{\cdot} = \sqrt{((\cdot, \cdot ))}$ is  the norm associated to the test space.
As in \cite{ferrero2022registration}, we consider a discrete $L^2$ inner product for the trial space and a discrete $H^1$ inner product for the test space such that
\begin{subequations} 
\label{eq:discrete_inner_products}
\begin{equation}
\left\{
\begin{array}{l}
\displaystyle{
(q,v)
=
\sum_{k=1}^{N_{\rm e}} \, 
\int_{\texttt{D}_k} 
  q \cdot  v 
\, dx
}
\\[3mm]
\displaystyle{
((q, v ))
\; =  \;
\sum_{k=1}^{N_{\rm e}} \, 
\int_{\texttt{D}_k} 
\left(
\nabla q   : \nabla v \, + \,   q \cdot  v 
\right)
\, dx
\; - \;
\sum_{j=1}^{N_{\rm f}} \, 
\int_{\texttt{F}_j} 
\left\{    \nabla q \mathbf{n}^+ \right\} \cdot J (v)
+
\left\{   \nabla v \mathbf{n}^+ \right\} \cdot J (q)
-
\eta \{  \texttt{r}_{j}( J (q)  )  \} \cdot J(v)
\, dx
}
\\
\end{array}
\right.
\end{equation}
where $\texttt{r}_{j}:[ L^2(\texttt{F}_j)]^D \to \mathcal{X}_{\rm hf}$ is the BR2  lifting operator (cf. \cite{bassi1997high}) given by
\begin{equation} 
\label{eq:BR2}
( \texttt{r}_{j}(w), v )
=
- \int_{\texttt{F}_j} w \cdot \{ v \} \, dx
\quad
\forall \,  
w \in  [ L^2(\texttt{F}_j)]^D, \;\; 
v\in \mathcal{X}_{\rm hf} 
\quad
j=1,\ldots,N_{\rm f},
\end{equation}
and $\eta>0$ is a stabilization parameter that is here set equal to $d+1$.
\end{subequations}
In the remainder of this section, we discuss the  construction of the various pieces of the formulation.

\subsection{Online solution method}
\label{sec:online_solution}
We denote by $\{   \psi_i \}_{i=1}^m$ an orthonormal basis of $\widehat{\mathcal{Y}}$;
we introduce the set of indices
$\texttt{I}_{\rm eq,e} = \{k\in \{1, \ldots,N_{\rm e} \} \, : \, \rho_{k}^{\rm eq,e} \neq  0 \}$ and
$\texttt{I}_{\rm eq,f} = \{j\in \{1, \ldots,N_{\rm f} \} \, : \, \rho_{j}^{\rm eq,f} \neq  0 \}$.
Then, we rewrite 
the minimization statement in \eqref{eq:MOR_intro_a} as the nonlinear least-square problem 
\begin{equation}
\label{eq:LSPG_algebraic}
\min_{  \boldsymbol{\alpha} \in \mathbb{R}^n  }
\left\|
\boldsymbol{\mathfrak{R}}_{\mu}^{\rm eq}(
 \boldsymbol{\alpha} )
\right\|_2,
\quad
{\rm with} 
\;\;
\left(
\boldsymbol{\mathfrak{R}}_{\mu}^{\rm eq}(
 \boldsymbol{\alpha} )
\right)_i =
\mathfrak{R}_{\mu}^{\rm eq}
(\texttt{Z}  \boldsymbol{\alpha}, \psi_i),
\;\; i=1,\ldots,m,
\end{equation}
which can be solved using the Gauss-Newton method
(GNM). Note that the computation of the entries of 
$\boldsymbol{\mathfrak{R}}_{\mu}^{\rm eq}(
 \boldsymbol{\alpha} )$ for any 
$\boldsymbol{\alpha} \in \mathbb{R}^n$  requires to compute the local elemental residuals 
$\{ r_{k,\mu}^{\rm e} \}_k$ for all $k\in \texttt{I}_{\rm eq,e}$ and the facet residuals
$\{ r_{j,\mu}^{\rm f} \}_j$ for all $j\in \texttt{I}_{\rm eq,f}$; towards this end, we should store the trial and test ROBs in the sampled elements
\begin{equation}
\label{eq:reduced_mesh}
\Omega_{\rm eq} :=
\left(
\bigcup_{k \in \texttt{I}_{\rm eq,e}}
\texttt{D}_k 
\right)
\cup
\left(
\bigcup_{j \in \texttt{I}_{\rm eq,f}}
\texttt{D}_j^+ \cup \texttt{D}_j^-\right).
\end{equation}
We conclude that online storage and computational costs scale linearly with the cardinality of 
$|\texttt{I}_{\rm eq,e}|$ and $|\texttt{I}_{\rm eq,f}|$.

Our formulation enables a straightforward
\emph{discretize-then-map} 
 treatment of geometry variations: 
  the elemental residual 
 $r_{k,\mu}^{\rm e}(\cdot,\cdot)$ depends on the nodes  $\{ x_{\texttt{T}_{i,k}}^{\rm hf} \}_{i=1}^{n_{\rm lp}}$ of the $k$-th element of the mesh;
 given a new value of the parameter $\mu$, it hence suffices to deform  the nodes of the sampled elements through the mapping $\Phi_{\mu}$  \emph{before} starting the GNM iterations. Similar reasoning applies to the facet integrals. As discussed in \cite{yano2019discontinuous} this approach enables the use of the routines of the DG HF code and is thus simple  to implement.
 
Several variants of the present approach are available in the literature.
 In \cite{yano2019discontinuous}, Yano considered an element-wise EQ procedure that guarantees relevant conservation properties, while in \cite{du2022efficient} Du and Yano proposed a pointwise EQ procedure that generates sparse quadrature rules within each element and facet. Our approach enables slightly larger reductions than the approach in \cite{yano2019discontinuous} and, unlike the approach in \cite{du2022efficient} can cope with elementwise  terms such as the BR2 lifting operator
(see \eqref{eq:BR2}) or elementwise artificial viscosities of the form \eqref{eq:viscosity_dilationbased}. 
A thorough comparison of our method with other EQ formulations is beyond the scope of this work.
 
 We finally comment on the choice of the initial condition for   GNM. We here rely on  nearest-neighbor regression:
 given the training set of simulations 
 $\{ \widetilde{q}_{\mu}^{\rm hf} : \mu \in \mathcal{P}_{\rm train}    \}$,  we define the corresponding best-fit generalized coordinates
 $\{  \boldsymbol{\alpha}_{\mu}^{\rm bf} : \mu \in \mathcal{P}_{\rm train}   \}$ obtained by projecting the available snapshots on the ROB $\texttt{Z}$; then, for any $\mu \in \mathcal{P}$, we initialize GNM with 
$ \boldsymbol{\alpha}_{\mu_{\rm nn}}^{\rm bf}$ 
with  $\mu_{\rm nn} = {\rm arg} \min_{\mu' \in \mathcal{P}_{\rm train}} \| \mu - \mu' \|_2$.
We observe that the present approach might be highly suboptimal if the cardinality of $\mathcal{P}_{\rm train}$ is modest: in section \ref{sec:adaptive_training}, we discuss how to improve the initialization of GNM using information from the previous iterations of Algorithm  
\ref{alg:offline_training}.

\subsection{Construction of the empirical test space}
\label{sec:test_space}
As in  \cite{ferrero2022registration}, 
we here resort to the sampling strategy 
based on
proper orthogonal decomposition
(POD, \cite{sirovich1987turbulence,volkwein2011model})
proposed in \cite{taddei2021space} to construct the test space $\widehat{\mathcal{Y}}$  
in \eqref{eq:MOR_intro}.
Given the training set $\mathcal{P}_{\rm train}=\{\mu^k \}_{k=1}^{n_{\rm train}} \subset \mathcal{P}$,
the associated snapshots $\{ \widetilde{q}_{\mu}^{\rm hf} : \mu\in \mathcal{P}_{\rm train}\}$, and the trial ROB $\{  \zeta_i \}_{i=1}^n$, we compute the test snapshot set 
$$
(({\psi}_{k,i}, v  ))
\, = \,
\mathfrak{J}_{\mu^k}^{\rm hf}[  \widetilde{q}_{\mu^k}^{\rm hf}  ]
(  {\zeta}_i, v ),
\quad
\forall \; v \in  \mathcal{X}_{\rm hf},
$$
for $i=1,\ldots,n$ and $k=1,\ldots,n_{\rm train}$,
where
$\mathfrak{J}_{\mu}^{\rm hf}[q]: \mathcal{X}_{\rm hf}\times  \mathcal{X}_{\rm hf} \to \mathbb{R}$ denotes 
 the Fr{\'e}chet derivative of the HF residual at $q$.
 Then, we perform POD on the 
 test snapshot set 
$\{ {\psi}_{k,i} \}_{k,i}$ based on the $((\cdot,\cdot))$ inner product \eqref{eq:discrete_inner_products} to obtain $\widehat{\mathcal{Y}}$.
In all the numerical experiments, we consider test spaces of size $j_{\rm es} =
{\rm dim}( \widehat{\mathcal{Y}} ) =  2n$; alternatively, we might   choose the dimension of $\widehat{\mathcal{Y}}$ using  an  energy criterion.
We refer to \cite[Appendix C]{taddei2021space} for a rigorous justification of our method for 
 linear inf-sup stable problems.

\subsection{Hyper-reduction}
\label{sec:hyper_reduction}
We seek 
$\boldsymbol{\rho}^{\rm eq,e} \in \mathbb{R}_+^{N_{\rm e}}$ and
$\boldsymbol{\rho}^{\rm eq,f} \in \mathbb{R}_+^{N_{\rm f}}$ in
\eqref{eq:MOR_intro_b}  
 such that
 \begin{enumerate}
 \item[(i)]
 (\emph{efficiency constraint})
the number of nonzero entries in 
 $\boldsymbol{\rho}^{\rm eq,e},\boldsymbol{\rho}^{\rm eq,f}$,  
 $\| \boldsymbol{\rho}^{\rm eq,e}   \|_{\ell^0}$ and
 $\| \boldsymbol{\rho}^{\rm eq,f}   \|_{\ell^0}$,   is as small as possible;
 \item[(ii)]
 (\emph{constant function constraint})
  the constant function is approximated correctly in $\Omega$ (for  ${\Phi} = \texttt{id}$), 
 \begin{equation}
 \label{eq:constant_function_constraint}
\Big|
\sum_{k=1}^{N_{\rm e}} \rho_k^{\rm eq,e} | \texttt{D}_k  |
\,-\,
| \Omega | 
 \Big|
 \ll 1,
 \quad
\Big|
\sum_{j=1}^{N_{\rm f}} \rho_j^{\rm eq,f} | \texttt{F}_j  |
\,-\,
\sum_{j=1}^{N_{\rm f}}  | \texttt{F}_j  |
 \Big|
 \ll 1; 
\end{equation}
\item[(iii)]
(\emph{manifold accuracy constraint})
for all $\mu \in \mathcal{P}_{\rm train,eq} = \{  \mu^k \}_{k=1}^{n_{\rm train}+n_{\rm train,eq}}$, the  empirical residual satisfies
\begin{subequations} 
\label{eq:accuracy_constraint}
\begin{equation}
\Big \|
   {\boldsymbol{\mathfrak{R}}}_{\mu}^{\rm hf}
( \boldsymbol{\alpha}_{\mu}^{\rm train}   )   
\, - \,
   {\boldsymbol{\mathfrak{R}}}_{\mu}^{\rm eq}
( \boldsymbol{\alpha}_{\mu}^{\rm train}   )   
 \Big \|_2
 \ll 1.
\end{equation}
where 
${\boldsymbol{\mathfrak{R}}}_{\mu}^{\rm hf}$  corresponds to substitute
$\rho_1^{\rm eq,e} = \ldots = \rho_{N_{\rm e}}^{\rm eq,e} = 
\rho_1^{\rm eq,f} = \ldots = \rho_{N_{\rm f}}^{\rm eq,f} = 1$ in
\eqref{eq:MOR_intro_b} 
 and $\boldsymbol{\alpha}_{\mu}^{\rm train}$ satisfies
\begin{equation}
\boldsymbol{\alpha}_{\mu}^{\rm train} = 
\left\{
\begin{array}{ll}
\displaystyle{ {\rm arg} \min_{\boldsymbol{\alpha} \in \mathbb{R}^n} \; 
\|   \texttt{Z}  \boldsymbol{\alpha}    - \widetilde{q}_{\mu}^{\rm hf} \|_2,}
 & {\rm if} \; \mu \in \mathcal{P}_{\rm train} ; \\[3mm]
\displaystyle{ {\rm arg} \min_{\boldsymbol{\alpha} \in \mathbb{R}^n} \; 
\| {\boldsymbol{\mathfrak{R}}}_{\mu}^{\rm hf}
( \boldsymbol{\alpha}   )   
\|_2,}
  & {\rm if} \; \mu \notin \mathcal{P}_{\rm train} ; \\
\end{array}
\right.
\end{equation}
and 
$\mathcal{P}_{\rm train} = \{  \mu^k \}_{k=1}^{n_{\rm train}}$ is the set of parameters for which the HF solution is available.
\end{subequations}
 \end{enumerate}

We refer to the above-mentioned literature for a thorough motivation of  the previous constraints.
We remark that several authors 
 (see \cite[Algorithm 1]{yano2019discontinuous})
have observed that considering an augmented training set  
$\mathcal{P}_{\rm train,eq}$
in \eqref{eq:accuracy_constraint} might  improve  performance of the hyper-reduced ROM, particularly for small values of $n_{\rm train}$.
However, for the numerical experiments of this work, we empirically observed  that the choice $\mathcal{P}_{\rm train}=\mathcal{P}_{\rm train,eq}$ leads to accurate results.

It is possible to show  (see, e.g.,  \cite{taddei2021space})   that (i)-(ii)-(iii) lead to a  sparse representation problem of the form
\begin{equation}
\label{eq:sparse_representation}
\min_{  \boldsymbol{\rho} \in \mathbb{R}^{N_{\rm e}   + N_{\rm f}}}
\;
\| \boldsymbol{\rho}   \|_{\ell^0},
\quad
{\rm s.t} \quad
\left\{
\begin{array}{l}
\|\mathbf{G} \boldsymbol{\rho} - \mathbf{b}  \|_2 \leq \delta; \\[3mm]
\boldsymbol{\rho} \geq \mathbf{0}; \\
\end{array}
\right.
\end{equation}
for a suitable threshold  $\delta>0$, and for  a suitable  choice of 
$\mathbf{G}, \mathbf{b}$. Following \cite{farhat2015structure}, we here resort to the non-negative  least-squares  method  to find approximate solutions to \eqref{eq:sparse_representation}. In particular, we use the Matlab  function \texttt{lssnonneq}, which takes as input the pair $(\mathbf{G}, \mathbf{b})$ and a tolerance $tol_{\rm eq}>0$ and returns the sparse vectors $\boldsymbol{\rho}^{\rm eq,e},\boldsymbol{\rho}^{\rm eq,f}$,
\begin{equation}
\label{eq:lsqnonneg}
[\boldsymbol{\rho}^{\rm eq,e}, \boldsymbol{\rho}^{\rm eq,f}] =  \texttt{lsqnonneg} \left( \mathbf{G}, \mathbf{b}, tol_{\rm eq} \right).
\end{equation}
We refer to \cite{chapman2017accelerated} for an efficient implementation of the non-negative least-squares method for large-scale problems.

\subsection{Construction of the trial space via greedy sampling}
\label{sec:greedy_sampling}

We resort to the weak-greedy algorithm 
(cf. \cite{veroy2003posteriori}) to build the ROM and the trial ROB $\texttt{Z}$; the weak-greedy method relies on the repeated maximization of an error indicator to adaptively sample the parameter domain; Algorithm \ref{alg:weak_greedy} summarizes the overall procedure, while Algorithm \ref{alg:construction_rom} summarizes the construction of the ROM.
In this work, we consider the residual-based error indicator
(cf. \cite[section 3.2.3]{ferrero2022registration}),
\begin{equation}
\label{eq:error_indicator}
\Delta: \mu\in \mathcal{P} \mapsto
\sup_{v\in \mathcal{X}_{\rm hf}} 
\frac{\mathfrak{R}_{\mu}^{\rm hf} (\widetilde{q}_{\mu}, v)}{\vertiii{v}}.
\end{equation}
Note that the evaluation of \eqref{eq:error_indicator} requires the solution to a linear system of size $N_{\rm hf}$: it is hence ill-suited for real-time online computations; nevertheless, in our experience the offline cost associated with the evaluation of \eqref{eq:error_indicator} is a fraction of the cost to perform  hyper-reduction and to build the test space $\widehat{\mathcal{Y}}$.
  We refer to 
\cite{ferrero2022registration} and to the references therein for a thorough discussion on the construction of an inexpensive surrogate of \eqref{eq:error_indicator}.
Even if we empirically observe that our residual-based error indicator is   highly-correlated with the true error, it does not provide a  rigorous bound; for this reason, after having computed the new HF solution (cf. Line 5, Algorithm \ref{alg:weak_greedy}) we check if the relative error is below a given threshold for the parameter that maximizes the error indicator.

\begin{algorithm}[H]                      
\caption{: weak-greedy algorithm. }     
\label{alg:weak_greedy}     

 \small
\begin{flushleft}
\emph{Inputs:}  $\mathcal{P}_{\rm train,gr} : =   \{ \mu_{\rm gr}^k \}_{k=1}^{n_{\rm train}} $ training parameter set, 
$\Phi:\Omega\times \mathcal{P} \to \Omega$ mapping;
$\mathcal{T}_{\rm hf}$ mesh.
\smallskip

\emph{Outputs:} 
$\texttt{Z}$ trial ROB;
$\mu \in \mathcal{P} \mapsto \widehat{\boldsymbol{\alpha}}_{\mu}$
ROM for the solution coefficients.

\end{flushleft}                      

 \normalsize 

\begin{algorithmic}[1]
\State
Choose 
$\mathcal{P}_{\star} =  \{ \mu^{\star, i}  \}_{i=1}^{n_0}$ and compute the HF solutions
$\mathcal{S}_{\star} =  \{ \widetilde{q}_{\mu}^{\rm hf}  : \mu \in  \mathcal{P}_{\star}   \} $.
\vspace{3pt}

\For {$n=n_0+1,\ldots,n_{\rm max}$}

\State
Update the ROB $\texttt{Z}$ and the ROM (cf. Algorithm  \ref{alg:construction_rom}).
\vspace{3pt}

\State
Estimate the solution   $\widetilde{q}_{\mu}^{\rm hf}$   and compute the indicator $\Delta_{\mu}$ in  \eqref{eq:error_indicator}
for all $\mu \in \mathcal{P}_{\rm train,gr}$.
\vspace{3pt}

\State
Compute $\widetilde{q}_{\mu^{\star,n}}^{\rm hf}$ for 
 $ \mu^{\star,n} = {\rm arg} \max_{ \mu\in \mathcal{P}_{\rm train,gr}  } \Delta_{\mu}$;
update $\mathcal{P}_{\star}$ and 
$\mathcal{S}_{\star}$. 
\vspace{3pt}

\If{$\| \widetilde{q}_{\mu^{\star,n}}^{\rm hf} - \widetilde{q}_{\mu^{\star,n}}   \| < {\texttt{tol}} 
\| \widetilde{q}_{\mu^{\star,n}}^{\rm hf}   \|$}

\State
Update the ROB $\texttt{Z}$ and the ROM.

\State
\texttt{break}
\EndIf

\EndFor
\end{algorithmic}
\bigskip

\end{algorithm}

\begin{algorithm}[H]                      
\caption{: construction of the ROM. }     
\label{alg:construction_rom}     

\small
\begin{flushleft}
\emph{Inputs:}  snapshot set 
$\mathcal{S}_{\star} : =   \{ (\mu, \widetilde{q}_{\mu}^{\rm hf}) : \mu\in \mathcal{P}_{\star} \}$.
\smallskip

\emph{Outputs:} 
$\texttt{Z}$ trial ROB;
$\mu \in \mathcal{P} \mapsto \widehat{\boldsymbol{\alpha}}_{\mu}$
ROM for the solution coefficients.

\end{flushleft}                      

 \normalsize 

\begin{algorithmic}[1]
\State
Define the test space $\widehat{\mathcal{Y}}$ (cf.
section \ref{sec:test_space})
\medskip

\State
Define the EQ weights $\boldsymbol{\rho}^{\rm eq,e}, \boldsymbol{\rho}^{\rm eq,f}$ (cf.
section \ref{sec:hyper_reduction}).
\medskip

\State
Store trial and test ROBs, and   grid points in the reduced mesh (cf.  \eqref{eq:reduced_mesh}).
\end{algorithmic}

\end{algorithm}

We observe that the weak-greedy algorithm requires multiple definitions of the ROM, which imply multiple constructions of the test space $\widehat{\mathcal{Y}}$, the quadrature weights
$\boldsymbol{\rho}^{\rm eq,e},\boldsymbol{\rho}^{\rm eq,f}$ and multiple greedy searches over the training set $\mathcal{P}_{\rm train, gr}$ (cf. Line 4, Algorithm \ref{alg:weak_greedy}). As reported in Table \ref{tab:offline_costs_transbump}, the overhead costs of the greedy procedure --- that is, the total  cost of the procedure minus the cost of the HF solves --- might be significant.
This observation motivates the development of more sophisticated training strategies to reduce offline costs. We address this issue in section \ref{sec:adaptive_training}.

\section{Adaptive procedure}
\label{sec:adaptive_training}
Each iteration of Algorithm \ref{alg:offline_training} generates a large amount of data about the parametric problem, which can be used to speed up offline computations.
In the remainder of this section, we illustrate computational bottlenecks of the training phase and we discuss  actionable strategies to reduce the computational burden; in the numerical investigations, we assess the impact of these choices.
\begin{itemize}
\item
The construction of the snapshot set for registration (cf. Line 3, Algorithm \ref{alg:offline_training}) based on the HF model is prohibitively expensive.
Instead, we propose to rely on the ROM built at the previous iteration; for the first iteration, we first execute the weak-greedy algorithm and then we use the ROM to generate the dataset of simulations.
\item
The GNM for \eqref{eq:LSPG_algebraic} is sensitive to the choice of the initial condition.
In our implementation, we initialize GNM based on nearest-neighbor regression which is clearly highly inaccurate for modest values of $n$. To face this issue, we propose to rely on a large dataset of initial conditions defined as follows:
\begin{equation}
\label{eq:initial_condition_ROM}
\left\{
\widehat{\boldsymbol{\alpha}}_{\mu}^0 : \mu\in \mathcal{P}_{\rm train}
\right\}
\quad
{\rm where}
\;\;
\widehat{\boldsymbol{\alpha}}_{\mu}^0
={\rm arg} \min_{ \boldsymbol{\alpha} \in \mathbb{R}^n  }
\;
\|  \texttt{Z}  \boldsymbol{\alpha}  - \widehat{q}_{\mu}^{\rm old} \circ \Phi_{\mu}^{-1} \|.
\end{equation} 
Note that the fields $\{ \widehat{q}_{\mu}^{\rm old} : \mu \in   \mathcal{P}_{\rm train} \} $ are generated for registration (cf. Line 3, Algorithm \ref{alg:offline_training}); nevertheless, cost of \eqref{eq:initial_condition_ROM} is significant due to the need 
to compute  the composition of $\widehat{q}_{\mu}^{\rm old} $ with $\Phi_{\mu}^{-1}$ --- which requires
 mesh interpolation.
In practice, we estimate the $L^2$ norm $\| \cdot    \|$ in \eqref{eq:initial_condition_ROM} using $10^3$  randomly-sampled points in $\Omega$ to reduce offline costs.
\item
The PTC strategy employed to solve the HF problem (cf. section \ref{sec:formulation})  might require many iterations to reach convergence. To reduce the computational burden, we initialize the PTC solver with the reduced-order solution $\widehat{q}_{\mu}$ from the previous iteration, as opposed to the free-stream solution.
Thanks to this choice, we can consider a much larger initial CFL number\footnote{In the numerical experiments, we set ${\rm CFL}_0=100$ instead of ${\rm CFL}_0=1$; see
\cite[section II.B]{yano2011importance}.} without experiencing any stability issue.
%{\color{red}Similarly, we might   exploit information from the previous iteration to initialize the GNM for \eqref{eq:LSPG_algebraic}:
%since the ROBs, the FE mesh and the mapping vary  at each iteration, 
%the conversion of the  generalized coordinates is computationally expensive; furthermore,
%for the model problems considered in the present work, 
% we did not find any benefit compared to nearest neighbor regression.
% }
\item
The registration procedure discussed in section \ref{sec:registration} relies on multiple solutions to a nonlinear non-convex optimization problem of size  $m=\mathcal{O}(10^2)$ that is sensitive to the initial condition. In our experience, the initialization strategy reviewed in section \ref{sec:param_reg}  leads to accurate performance; however, it requires a sufficiently dense discretization of $\mathcal{P}$ and is not parallelizable. To address this issue, we propose to store the mapping coefficients $\{  \widehat{\mathbf{a}}_{\mu} : \mu\in \mathcal{P}_{\rm train}  \}$ obtained during the first iteration of the registration method and then use them as  initial conditions for the subsequent iteration: note that
for this choice of the initialization the solution to the problems \eqref{eq:optimization_all_params} can be
 trivially parallelized; in addition, we can
{potentially} 
  cope with much coarser discretizations of $\mathcal{P}$.
\item
As discussed in section  \ref{sec:pMOR}, the weak-greedy algorithm requires multiple constructions of the ROM and might hence be expensive; in addition, it cannot be efficiently parallelized.
To address this issue, we initialize Algorithm \ref{alg:weak_greedy} with the parameters $\{\mu^{\star,i}\}_{i=1}^{n_0}$ obtained by applying the strong-greedy algorithm
  to the snapshot set generated for registration (cf. Line 3, Algorithm \ref{alg:offline_training}). Since the snapshot set is generated using the ROM, the strong-greedy algorithm can be applied to the generalized coordinates.
For completeness, we report the strong-greedy procedure in Appendix \ref{sec:appendix_greedy}.
\end{itemize}

\section{Numerical results}
\label{sec:numerics}
We present below extensive numerical investigations for the model problems introduced in section \ref{sec:model_problems}. Further numerical tests are provided in 
Appendix \ref{sec:appendix_transbump}. 
We assess performance based on $n_{\rm test}=20$ out-of-sample parameters $\mathcal{P}_{\rm test}=\{  \mu_{\rm test}^j \}_{j=1}^{n_{\rm test}}$ with 
$\mu_{\rm test}^1,\ldots,\mu_{\rm test}^{n_{\rm test}} \overset{\rm iid}{\sim} {\rm Uniform}(\mathcal{P})$; for each $\mu \in \mathcal{P}_{\rm test}$, we report the HF  $L^2$ error $E_{\mu}^{\rm hf}$, the sub-optimality index 
$\eta_{\mu}^{\rm hf}$ and the total enthalpy error
$E_{\mu}^{\infty}$ such that
\begin{equation}
\label{eq:L2error}
E_{\mu}^{\rm hf} = \frac{\|  q_{\mu}^{\rm hf} -  \widehat{q}_{\mu}^{\rm hf}  \|_{L^2(\Omega_{\mu})}     }{\|  q_{\mu}^{\rm hf}    \|_{L^2(\Omega_{\mu})} },
\end{equation} 
\begin{equation}
\label{eq:suboptimality_index}
\eta_{\mu}^{\rm hf} = \frac{\|  q_{\mu}^{\rm hf} -  \widehat{q}_{\mu}^{\rm hf}  \|_{L^2(\Omega_{\mu})}     }{
\min_{\zeta \in \mathcal{Z}_n}
\|  q_{\mu}^{\rm hf}  - \zeta \circ \Phi_{\mu}^{-1}   \|_{L^2(\Omega_{\mu})} },
\end{equation} 
\begin{equation}
\label{eq:total_enthalpy_error}
E_{\mu}^{\infty} = \frac{\|  H_{\rm tot, \mu}^{\rm true} -  \widehat{H}_{\rm tot, \mu}  \|_{L^2(\Omega_{\mu})}     }{
\|  H_{\rm tot, \mu}^{\rm \infty}    \|_{L^2(\Omega_{\mu})} }.
\end{equation} 
The relative error $E_{\mu}^{\rm hf}$ measures the accuracy of the reduced-order estimate with respect to the HF model employed  for training  --- it is hence a measure of the overall ability of the MOR procedure to approximate the truth model of the PDE. 
The suboptimality index 
$\eta_{\mu}^{\rm hf} $ measures the extent to which the LSPG projection scheme is suboptimal compared to the best-fit error: it hence allows to directly evaluate the effectiveness of the ROM, which encompasses the choice of the test space, initialization, and hyper-reduction. 
Finally, the total enthalpy error \eqref{eq:total_enthalpy_error}  measures the accuracy of the state estimate with respect to the exact solution to the PDE, in terms of enthalpy preservation.
Simulations are
performed in Matlab 2022a
\cite{MATLAB:2022} based on an in-house code, and executed over a commodity Linux  workstation (RAM 32 GB, Intel i7 CPU 3.20 GHz x 12).

\subsection{Inviscid flow through a nozzle}
\label{sec:nozzle_numerics}
We perform $N_{\rm it}=3$ iterations of Algorithm \ref{alg:offline_training} without acceleration. We initialize the algorithm using an uniform HF grid with $N_{\rm e}=60$ triangles and quadratic ($\texttt{p}=2$) polynomials; then,  we increase the size of the mesh by a factor $1.5$ at each iteration: this implies that the generated HF meshes have $N_{\rm e}=60$, $N_{\rm e}=90$ and $N_{\rm e}=135$ elements at iterations one, two and three, respectively.
We consider a regular $15 \times 15$ grid of parameters
$\mathcal{P}_{\rm train}$
 for registration and a regular $10 \times 10$ grid of parameters
$\mathcal{P}_{\rm train, gr}$
 in  Algorithm \ref{alg:weak_greedy}. We rely on the HF solver to generate the dataset of simulations at iteration one, while we rely on the ROM from previous iterations to generate the snapshot set (cf. Line 3, Algorithm \ref{alg:offline_training}) for $k=2,\ldots,N_{\rm it}$. We consider the tolerance $\texttt{tol}=10^{-3}$ in 
Algorithm \ref{alg:weak_greedy} and we consider an initial regular $3\times 3$ grid of parameters to initialize the ROM: the algorithm generates ROBs of size $n=15$, $n=20$ and $n=10$.

\begin{figure}[H]
\centering
 \subfloat[] 
{  \includegraphics[width=0.4\textwidth]
 {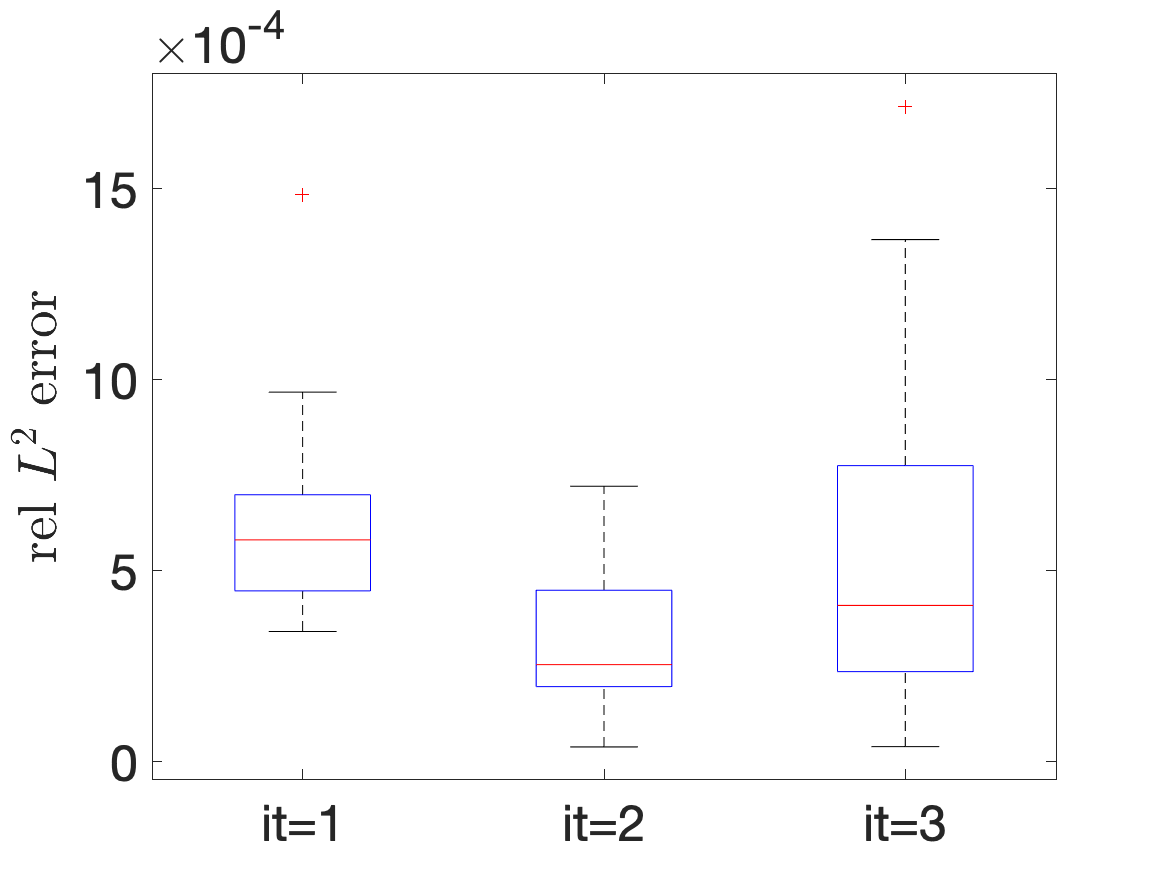}}
   ~~
 \subfloat[] 
{  \includegraphics[width=0.4\textwidth]
 {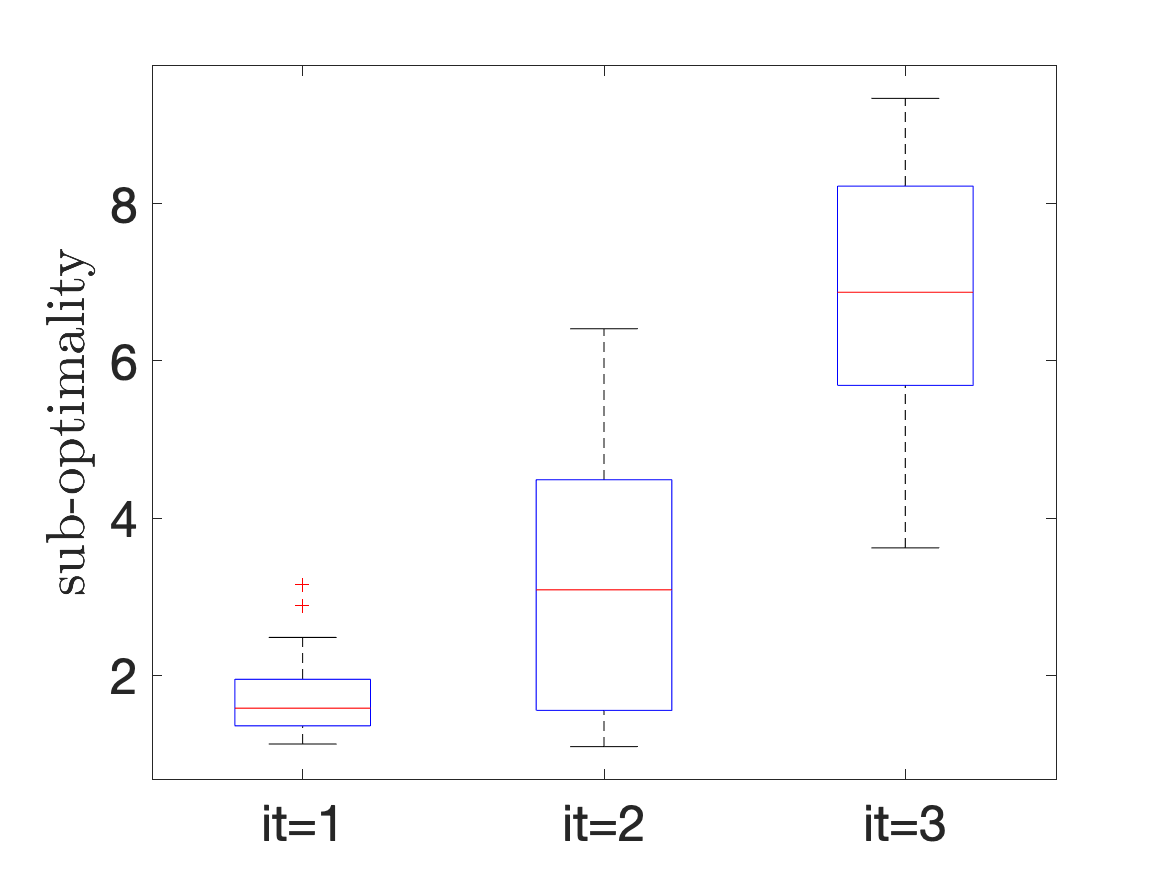}}

 \subfloat[] 
{  \includegraphics[width=0.4\textwidth]
 {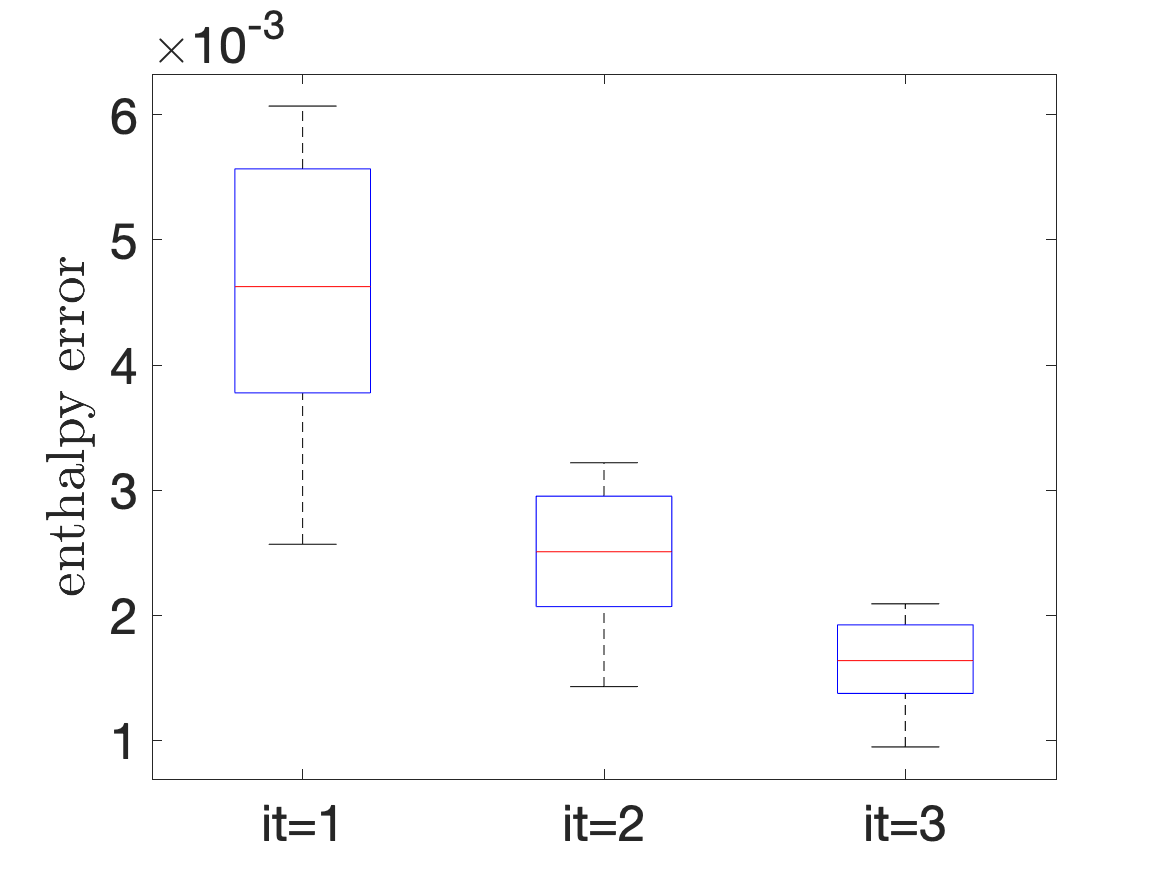}}
 ~~
  \subfloat[] 
{  \includegraphics[width=0.4\textwidth]
 {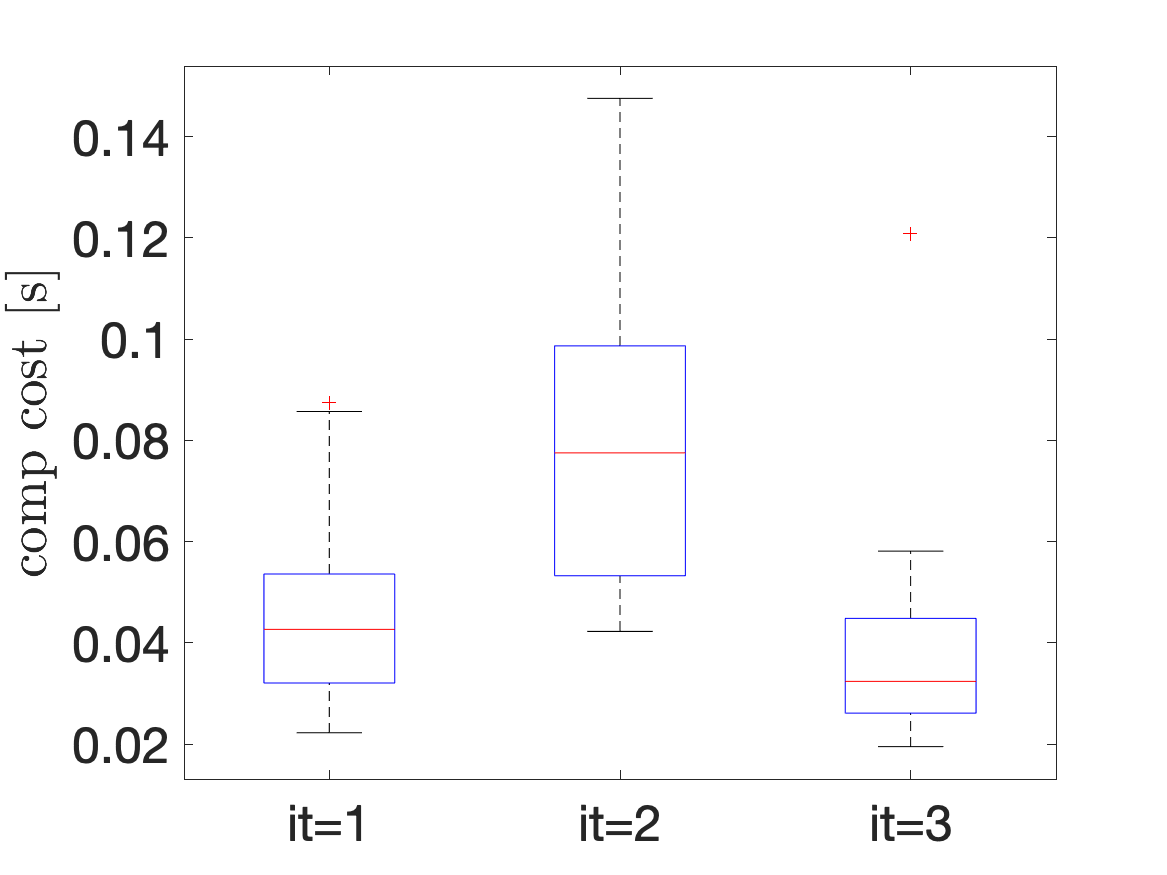}}
 
 \caption{nozzle flow. Performance of the ROM for three iterations of the adaptive (basic) procedure.
}
 \label{fig:nozzle_basic_rom}
 \end{figure}

Figure \ref{fig:nozzle_basic_rom} shows the performance of the ROM.
Figure \ref{fig:nozzle_basic_rom}(a) shows the relative error over the test set, which mildly depends on the size of the mesh.
Figure \ref{fig:nozzle_basic_rom}(b) shows the suboptimality index: interestingly, we observe that the performance of the projection scheme deteriorates as we increase the size of the mesh: we plan to investigate this behavior in a subsequent work; nevertheless, we observe that $\eta_{\mu}^{\rm hf} \lesssim 10$ for all numerical experiments.
Figure \ref{fig:nozzle_basic_rom}(c) shows the total enthalpy error: as expected, the error  decreases as we increase the size of the mesh. 
Figure \ref{fig:nozzle_basic_rom}(d) shows the wall-clock online cost: thanks to hyper-reduction, results do not depend on the size of the underlying mesh but they clearly depend on the size $n$ of the ROB. 

%nozzle_basic_rom

\begin{figure}[H]
\centering
 \subfloat[$it=1$] 
{  \includegraphics[width=0.33\textwidth]
 {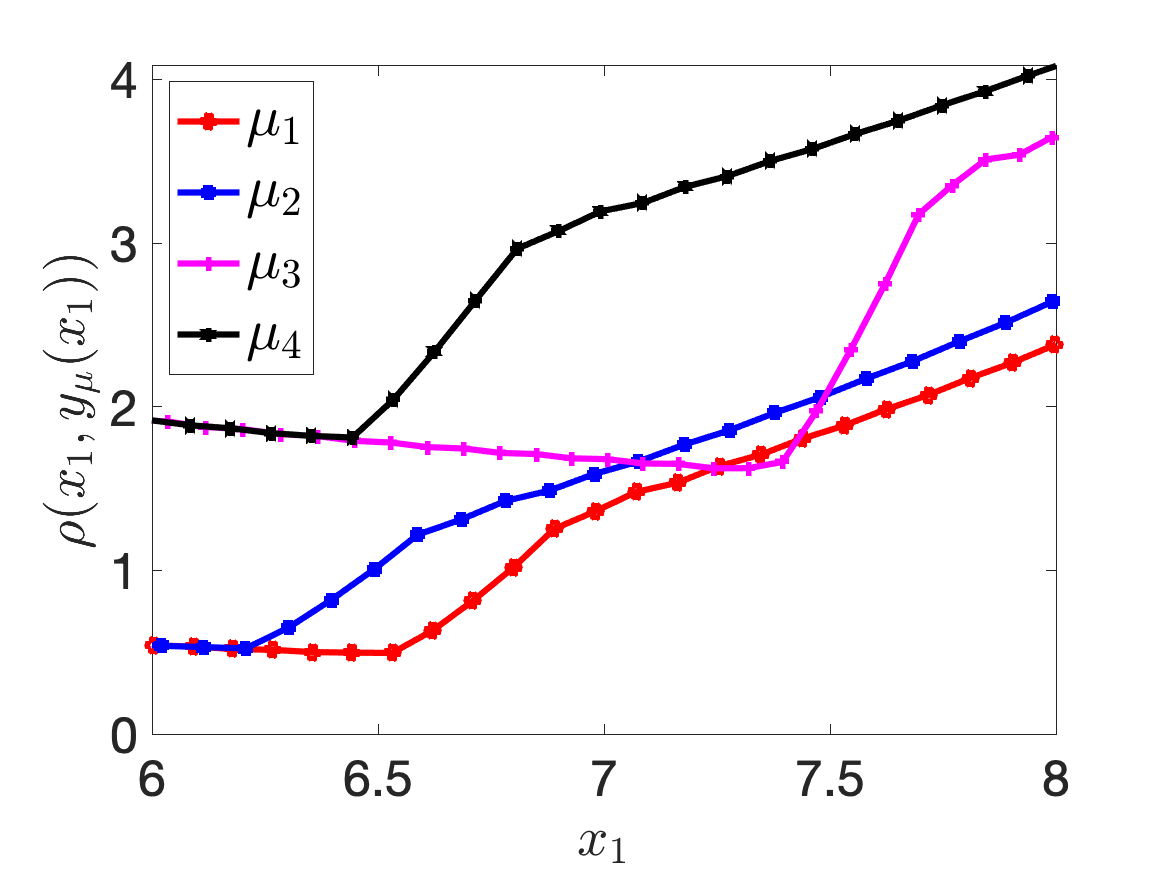}}
   ~~
 \subfloat[$it=2$] 
{  \includegraphics[width=0.33\textwidth]
 {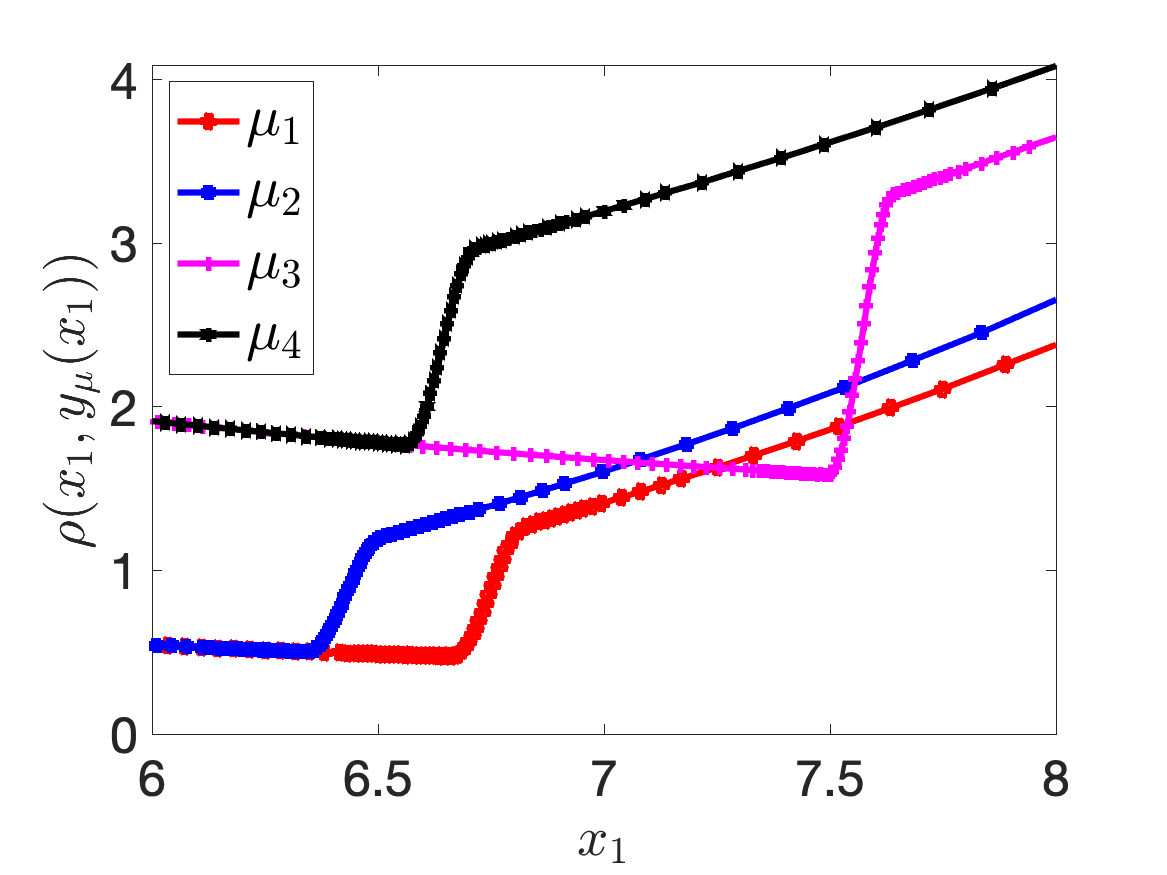}}
~~
 \subfloat[$it=3$] 
{  \includegraphics[width=0.33\textwidth]
 {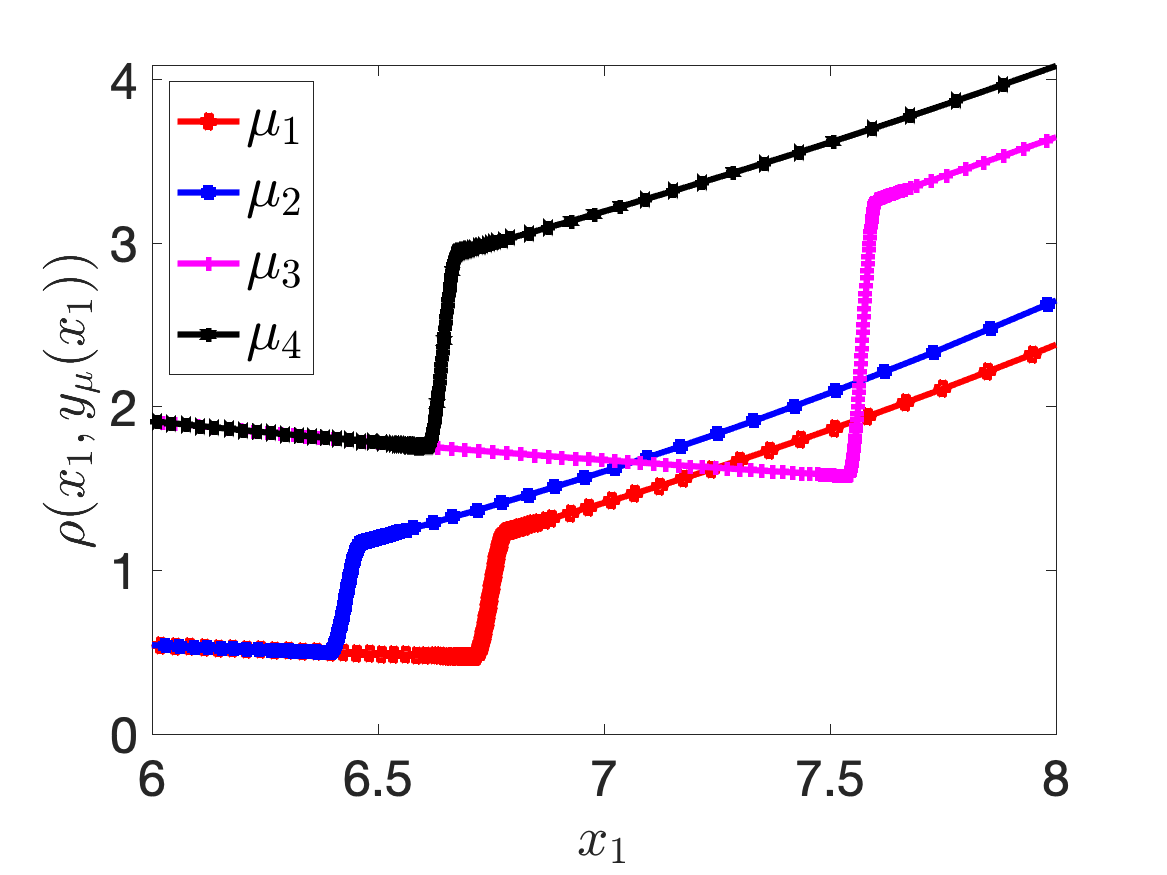}}
 
 \subfloat[$it=1$] 
{  \includegraphics[width=0.33\textwidth]
 {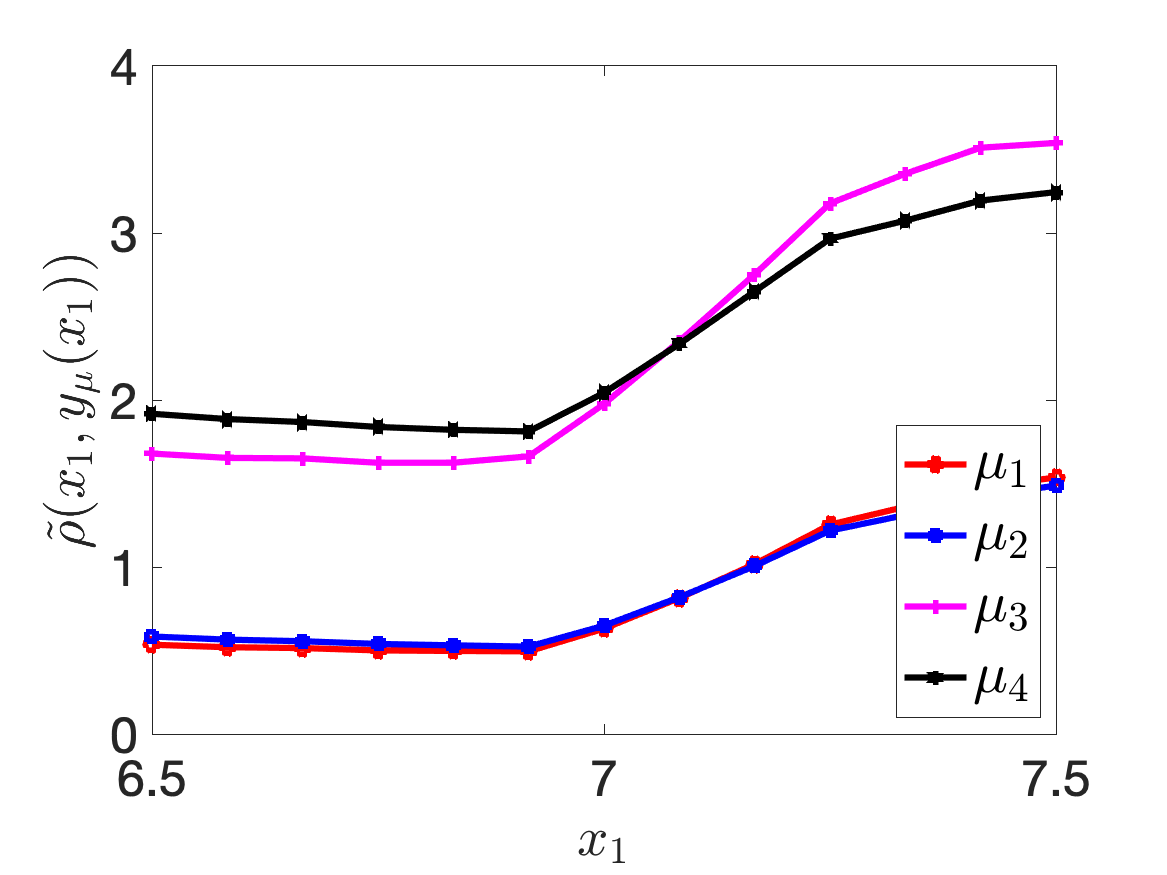}}
   ~~
 \subfloat[$it=2$] 
{  \includegraphics[width=0.33\textwidth]
 {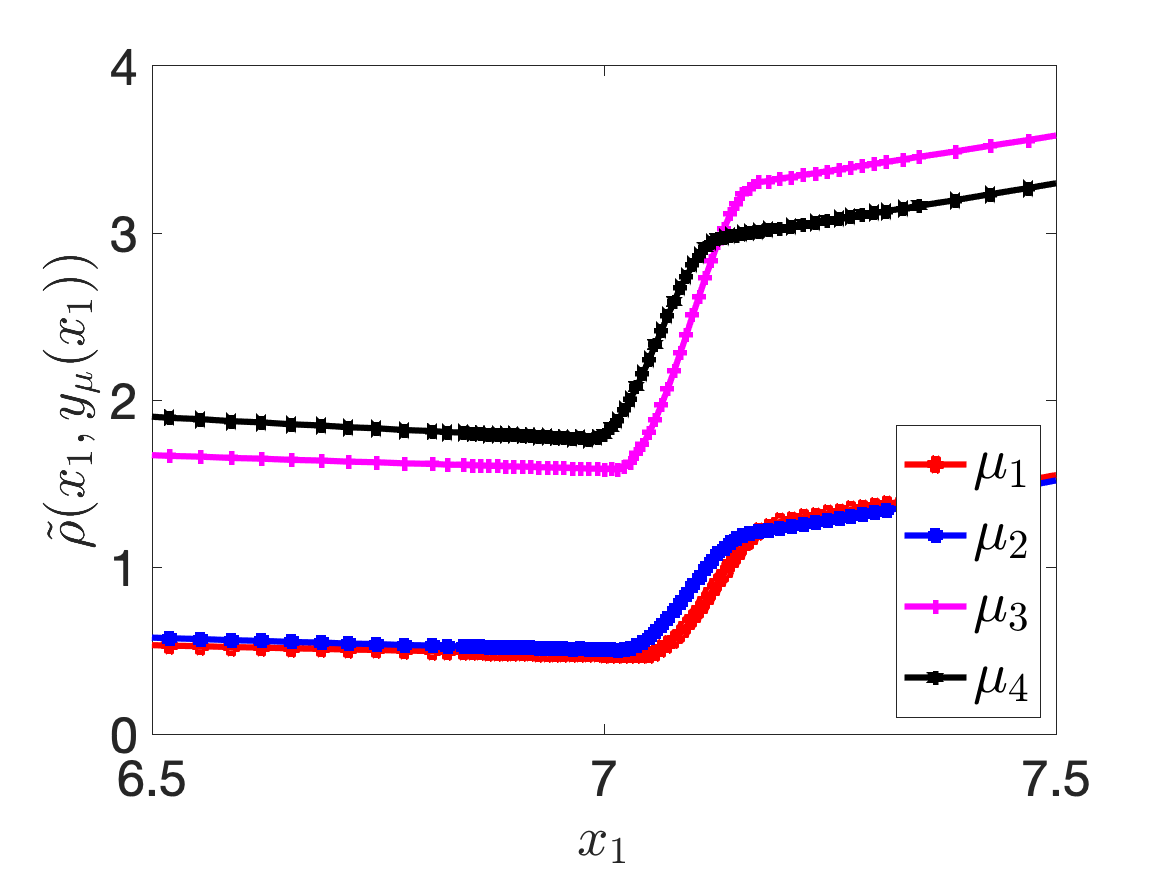}}
~~
 \subfloat[$it=3$] 
{  \includegraphics[width=0.33\textwidth]
 {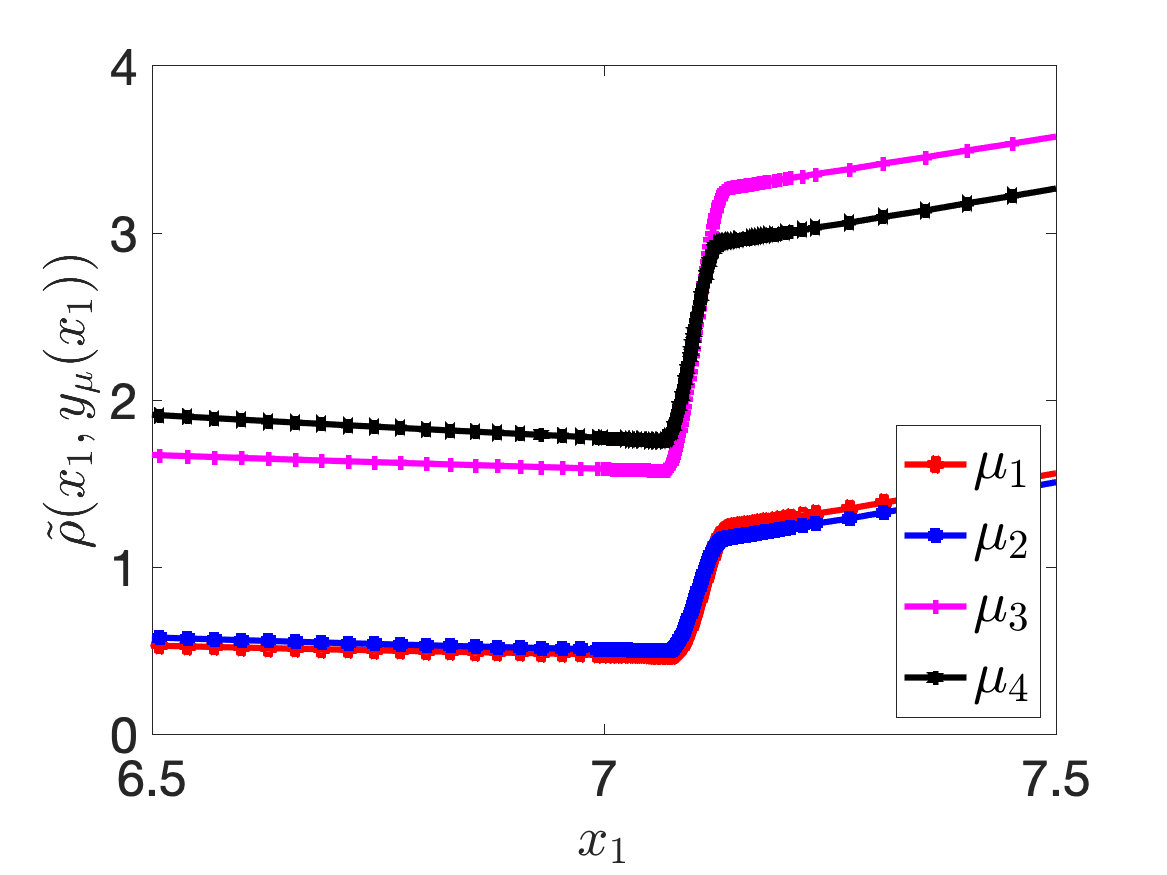}}
 
 \caption{nozzle flow. Behavior of the (modified) density field  in physical (cf. (a)-(b)-(c)) and reference 
(cf. (d)-(e)-(f))   configuration for four values of the parameter and three iterations of the adaptive algorithm (basic version).
}
 \label{fig:nozzle_basic_registration}
 \end{figure}

Figure \ref{fig:nozzle_basic_registration} shows the behavior of the modified density $\bar{\rho}:= A \rho$ in the proximity of the shock for four parameter values, for three iterations of the algorithm, and for both physical and reference configurations. We observe that registration is effective to track the position of the discontinuity. 

\begin{figure}[H]
\centering
 \subfloat[] 
{  \includegraphics[width=0.33\textwidth]
 {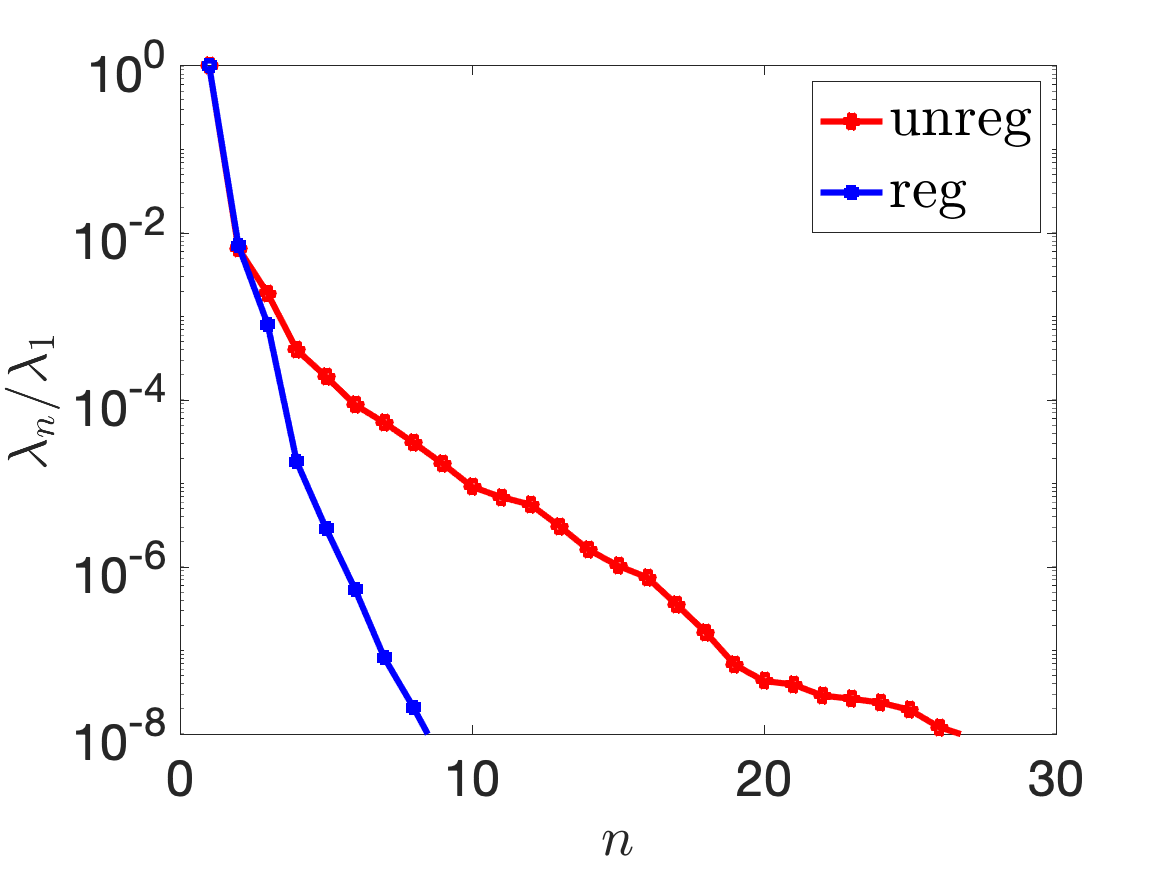}}
   ~~
 \subfloat[] 
{  \includegraphics[width=0.33\textwidth]
 {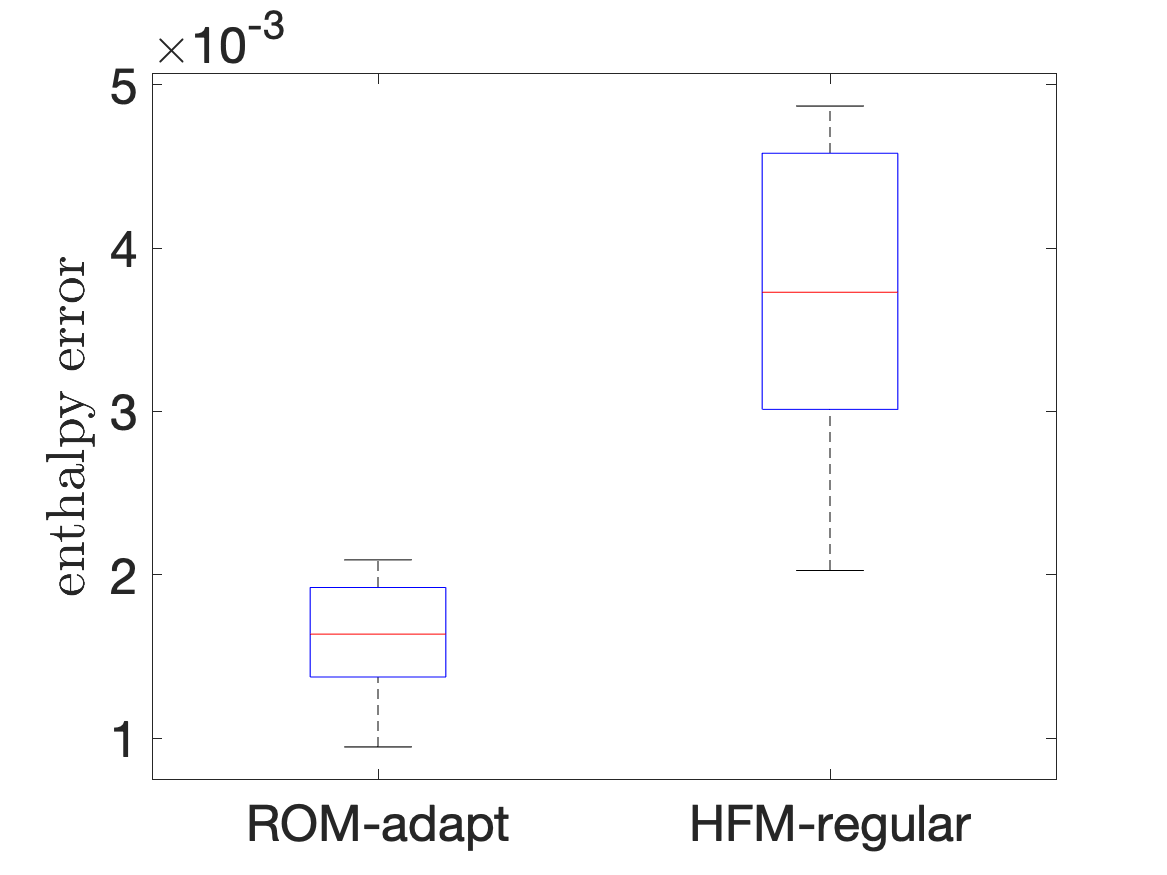}}
~~
 \subfloat[] 
{  \includegraphics[width=0.33\textwidth]
 {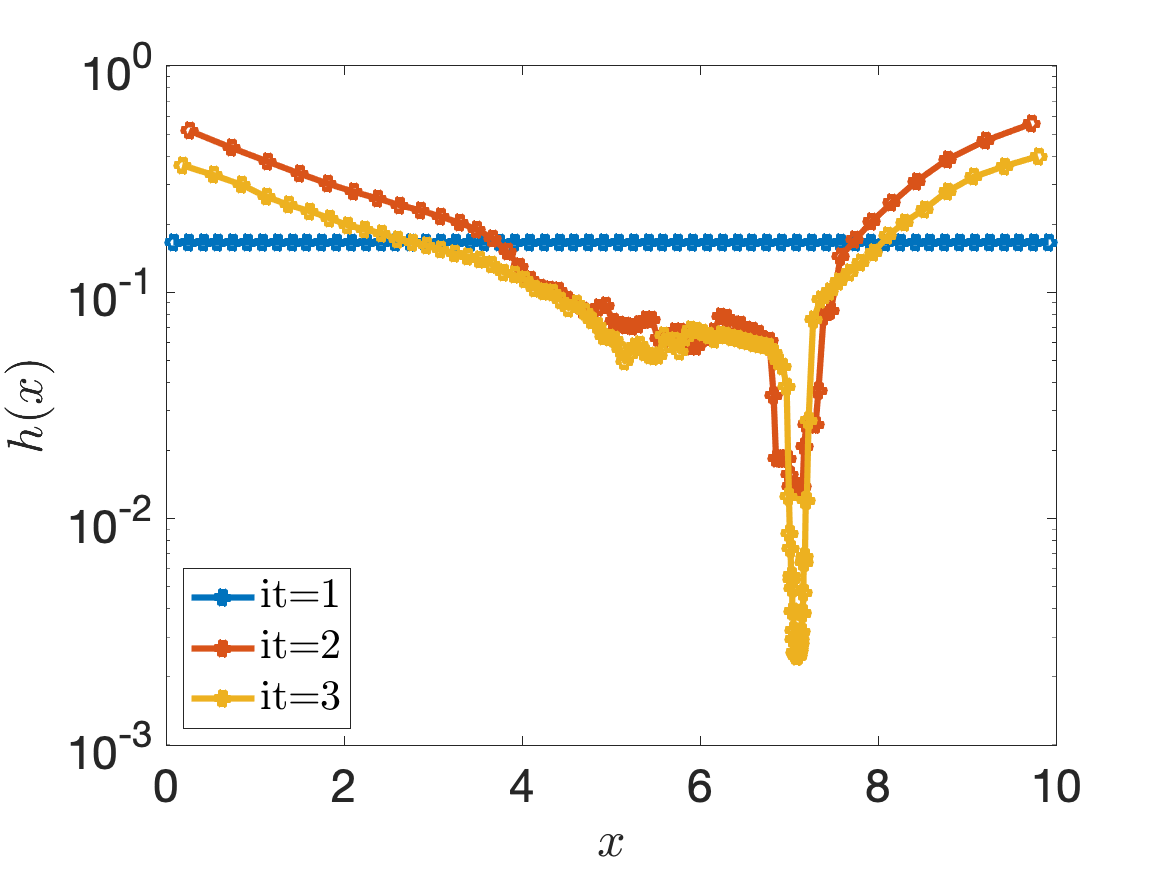}}
 
\caption{nozzle flow. Effect of registration on compressibility and mesh adaptation.
(a) POD eigenvalues in reference and physical configurations.
(b) total enthalpy error for registered ROM and unregistered HFM  on regular mesh with  $N_{\rm e}=135$ elements.
(c) mesh density 
$h:\Omega\to \mathbb{R}_+$, 
$h|_{\texttt{D}_k}=| \texttt{D}_k|$ for the sequence of considered meshes.
}
 \label{fig:nozzle_advanced_plot}
 \end{figure}

In Figure \ref{fig:nozzle_advanced_plot}, we investigate the effect of registration on solution manifold compressibility and mesh adaptation. 
Towards this end, we consider the adaptive reduced-order and   HF models 
associated with the third iteration of Algorithm \ref{alg:offline_training}, and a 
  HF model defined over a Cartesian  ``static''
mesh with the same number of elements, 
$N_{\rm e}=135$.
First, in Figure \ref{fig:nozzle_advanced_plot}(a), 
we compare the behavior of the normalized POD eigenvalues associated with the snapshot set in physical (``unreg'') and reference (``reg'')  configurations.
 We observe that registration significantly improves the convergence of the POD eigenvalues that can be regarded as a ``proxy'' of the linear complexity of the corresponding solution manifold.
Figure \ref{fig:nozzle_advanced_plot}(b)  shows the behavior of the error in total enthalpy for the final registered ROM and the static HF model based on an uniform mesh: we clearly notice that the HF model
--- which has the same number of degrees of freedom as the HF model used to generate the ROM ---
 is significantly less accurate than the adapted ROM.
Finally, Figure \ref{fig:nozzle_advanced_plot}(c) shows the behavior of the mesh density $h:\Omega\to \mathbb{R}_+ $ such that $h|_{\texttt{D}_k} = |\texttt{D}_k|$ for the sequence of meshes generated during Algorithm \ref{alg:offline_training}: we observe that registration allows us to refine the mesh over a very narrow portion of the computational domain and hence enables significant computational savings.

%nozzle_registration

% \label{fig:nozzle_advanced_plot}

\subsection{Inviscid  flow over a Gaussian bump}
\label{sec:transonic_bump_numerics}

We perform $N_{\rm it}=3$ iterations of Algorithm \ref{alg:offline_training} without  and with acceleration; we consider both isotropic and anisotropic mesh adaptation based on the software \texttt{mmg2d} and on the metrics introduced in 
section \ref{sec:mesh_adaptation}.
As in the previous case, we rely on a regular 
$15 \times 15$ grid of parameters
$\mathcal{P}_{\rm train}$
 for registration and a regular $10 \times 10$ grid of parameters $\mathcal{P}_{\rm train, gr}$
 in  Algorithm \ref{alg:weak_greedy}. 
 We set $\texttt{tol}=10^{-3}$ in the termination condition of Algorithm \ref{alg:weak_greedy}.
 To reduce training costs of the first snapshot generation, we first perform a weak-greedy algorithm to generate a ROM that is later used to generate the snapshot set.
 In all our tests, we consider the initial grid depicted in Figure \ref{fig:transbump_basic_mesh}(a) with $N_{\rm e}=3448$, and we rely on a quadratic approximation.
 We state upfront that the registration algorithm returns a low-rank mapping with $m=2$ modes for all runs considered.
 
\subsubsection{Basic approach}
\label{sec:transonic_bump_numerics_basic}

\begin{figure}[H]
\centering
 \subfloat[] 
{  \includegraphics[width=0.4\textwidth]
 {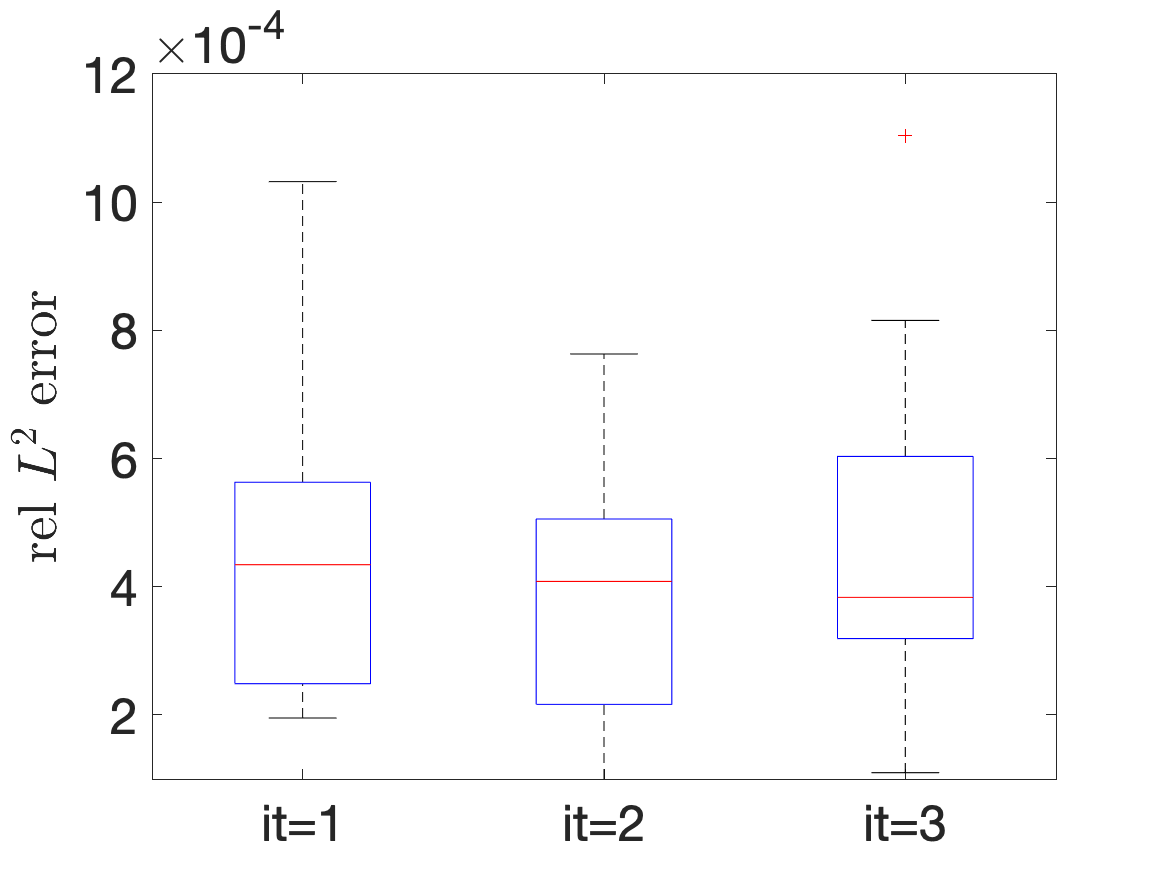}}
   ~~
 \subfloat[] 
{  \includegraphics[width=0.4\textwidth]
 {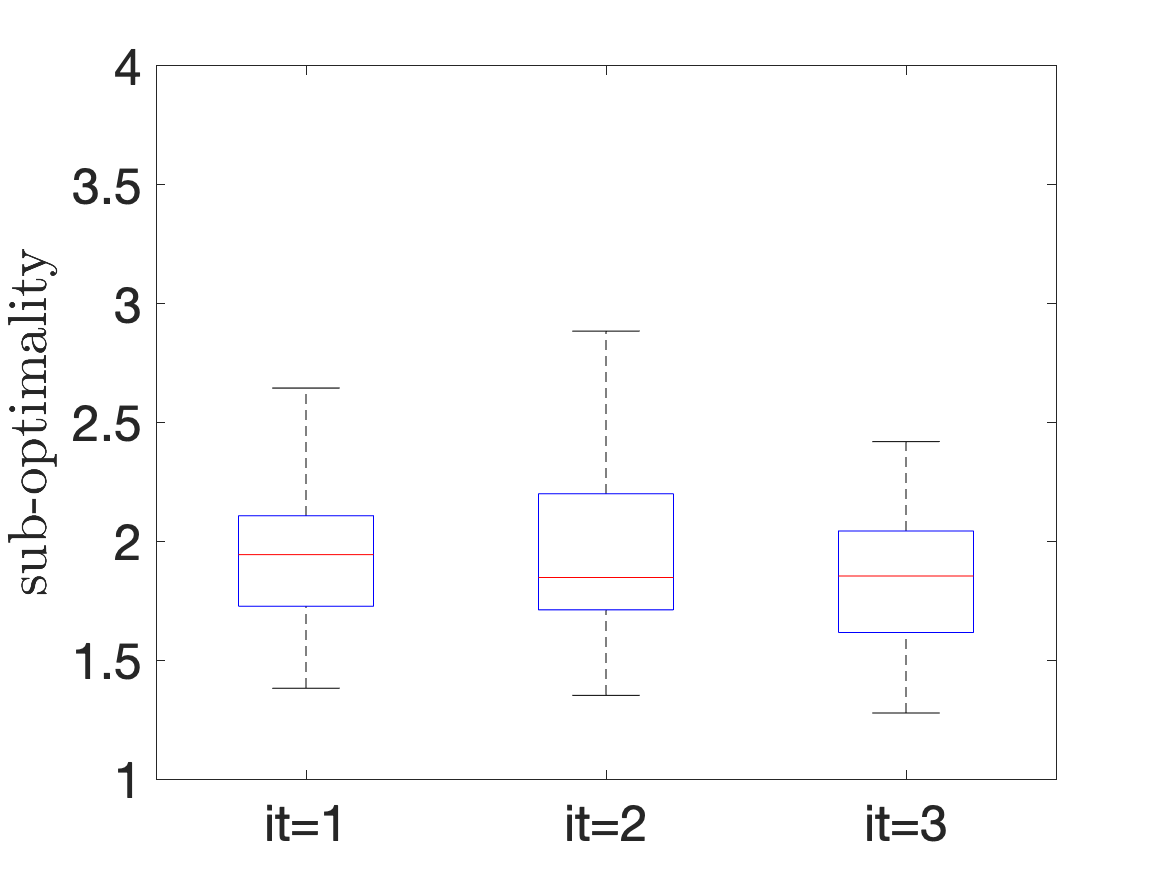}}

 \subfloat[] 
{  \includegraphics[width=0.4\textwidth]
 {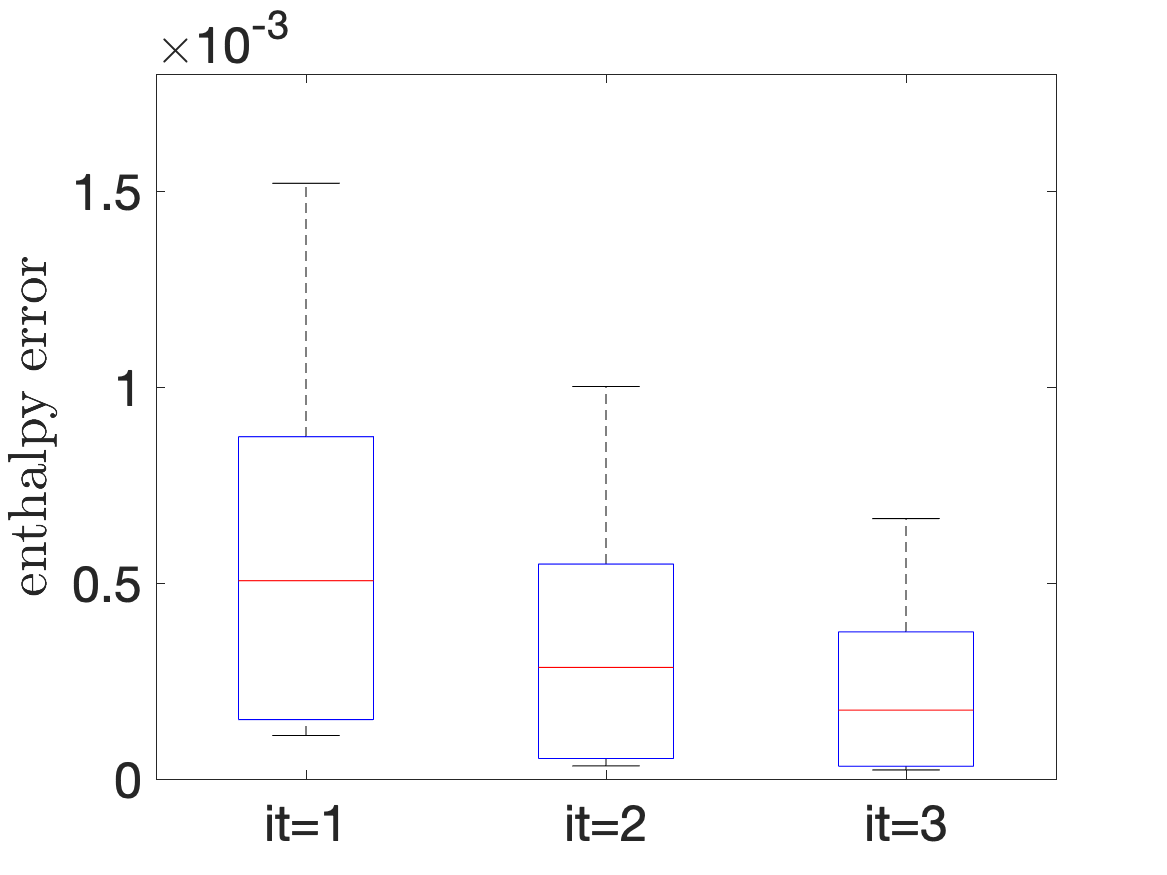}}
 ~~
  \subfloat[] 
{  \includegraphics[width=0.4\textwidth]
 {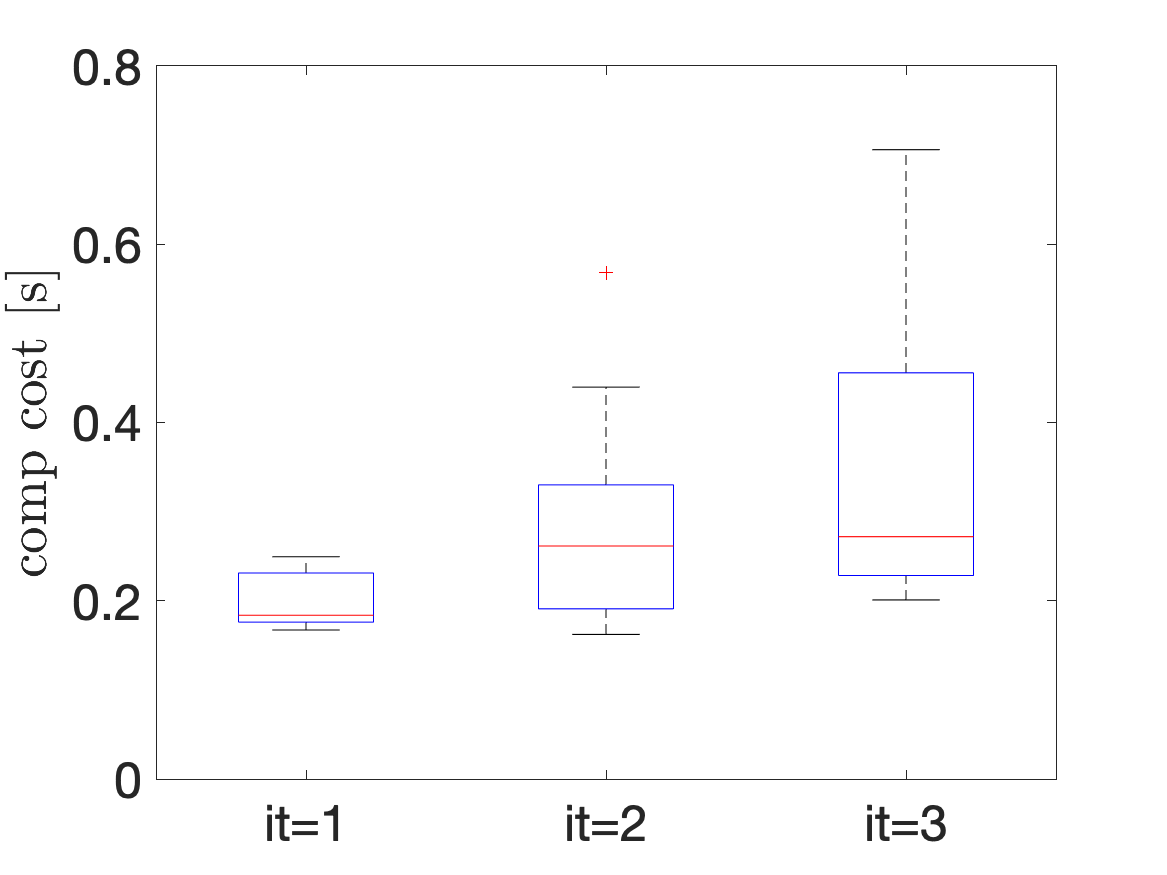}}
 
 \caption{transonic bump. Performance of the ROM for three iterations of the adaptive (basic) procedure.
}
 \label{fig:transbump_basic_rom}
 \end{figure}

We first study the performance of the standard (without acceleration) approach based on isotropic mesh adaptation.
Figure \ref{fig:transbump_basic_rom} replicates the results in Figure \ref{fig:nozzle_basic_rom} for the transonic bump problem.
We observe that the ROM achieves accurate performance over the test set with respect to the HF estimate for all three iterations: results are hence in good agreement with the selected tolerance ($\texttt{tol}=10^{-3}$)  of Algorithm \ref{alg:weak_greedy}.
The suboptimality index ranges from one to three for all experiments: this indicates that our projection scheme is extremely effective for this model problem.
The total enthalpy error decreases as we increase the size of the mesh, while the computational cost is nearly the same for all iterations.

%transbump_basic_rom

\begin{figure}[H]
\centering
 \subfloat[$it=1$] 
{  \includegraphics[width=0.33\textwidth]
 {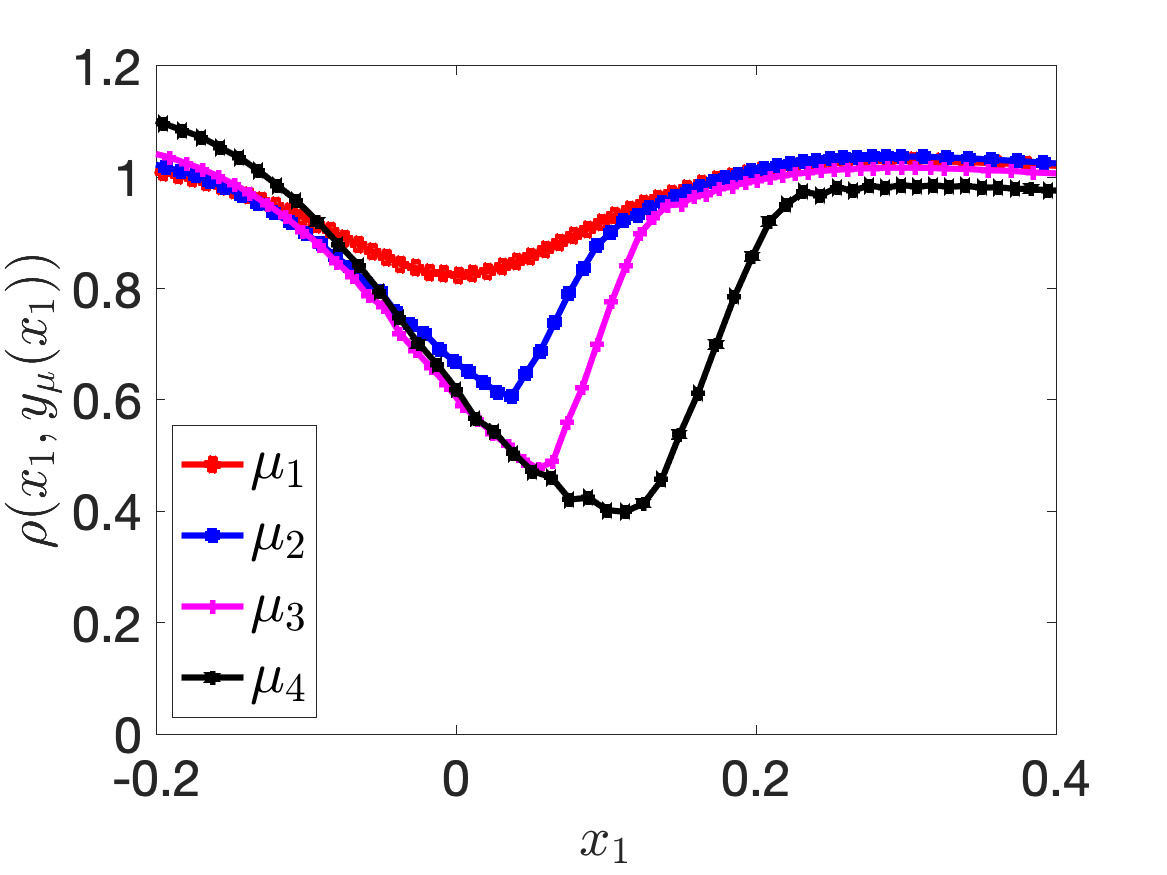}}
   ~~
 \subfloat[$it=2$] 
{  \includegraphics[width=0.33\textwidth]
 {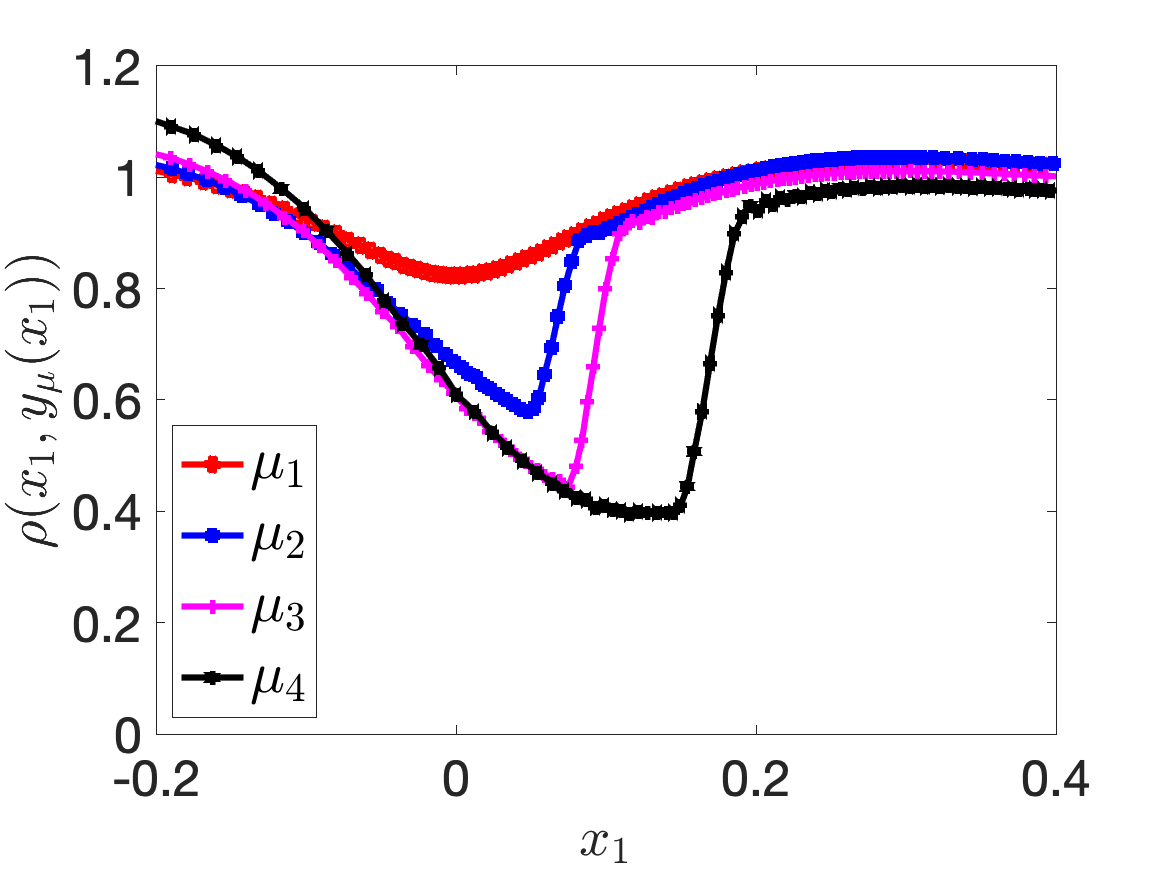}}
~~
 \subfloat[$it=3$] 
{  \includegraphics[width=0.33\textwidth]
 {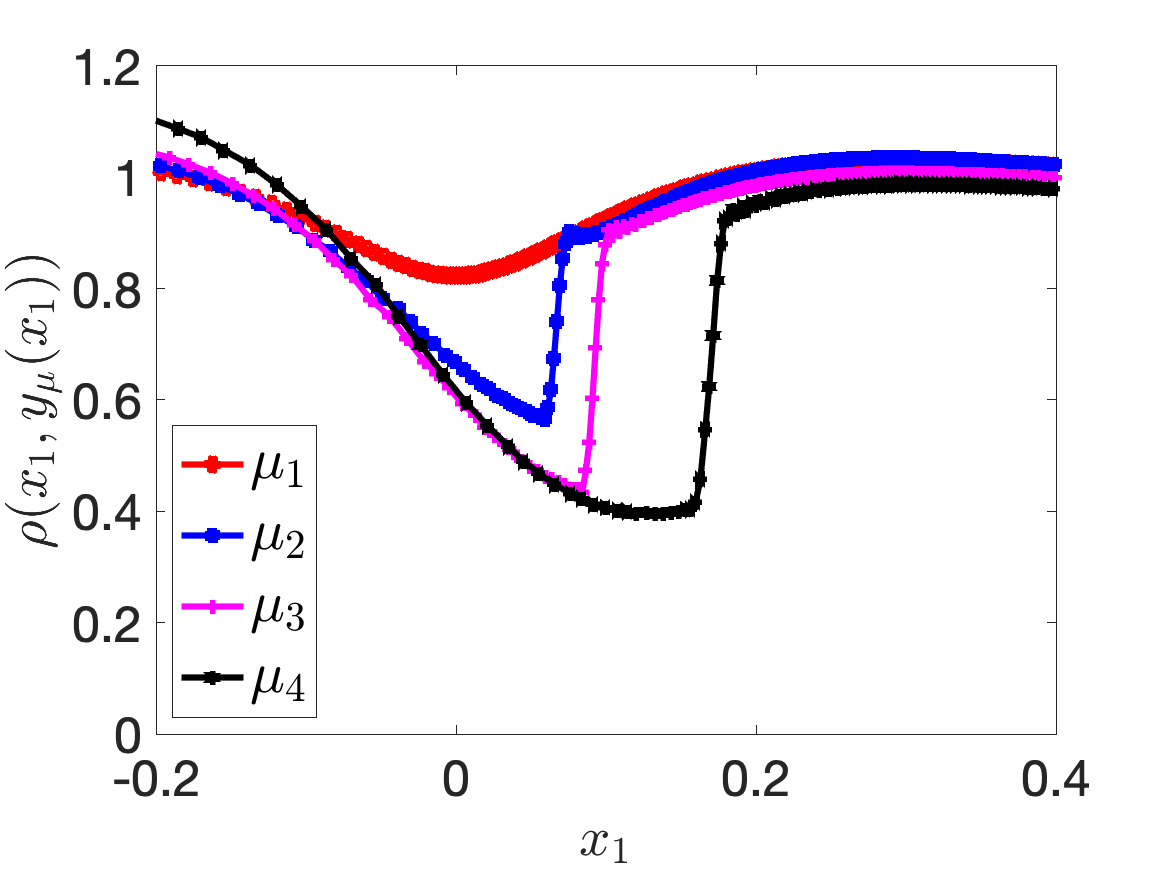}}
 
 \subfloat[$it=1$] 
{  \includegraphics[width=0.33\textwidth]
 {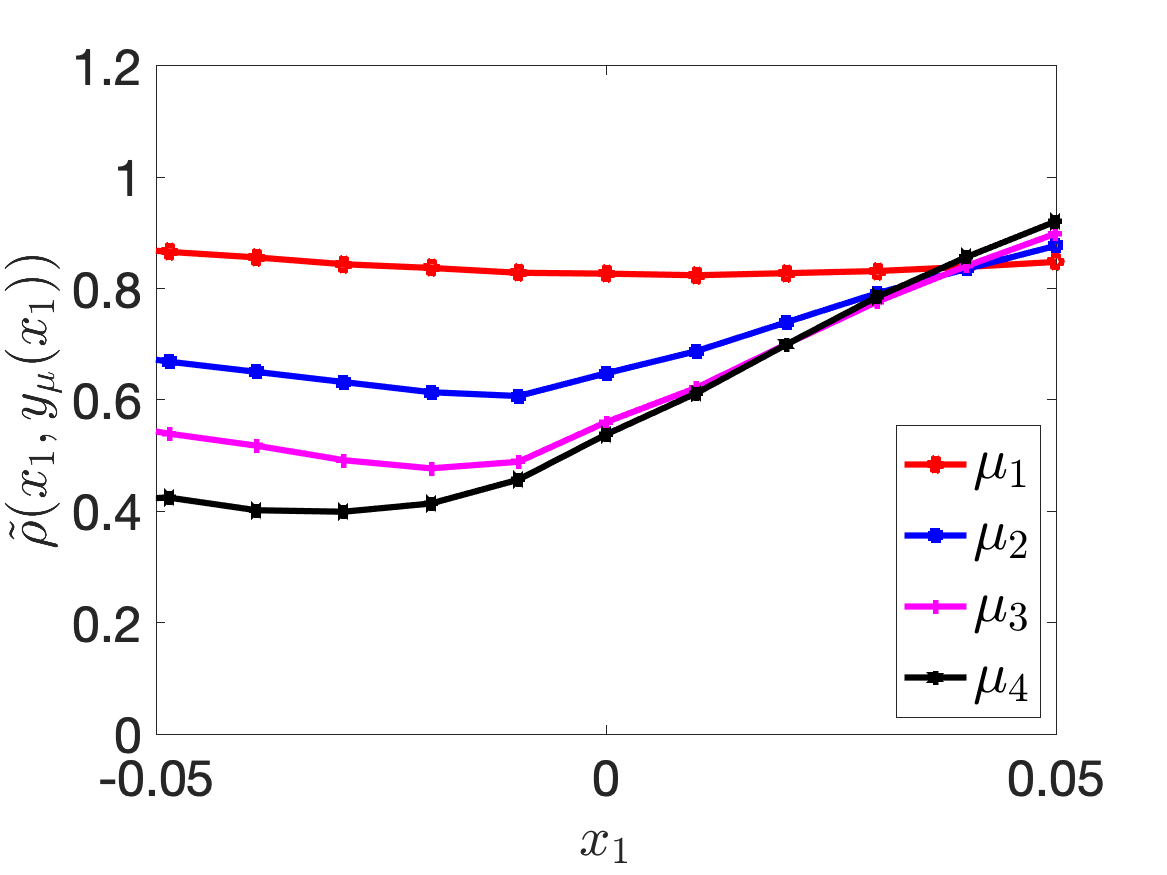}}
   ~~
 \subfloat[$it=2$] 
{  \includegraphics[width=0.33\textwidth]
 {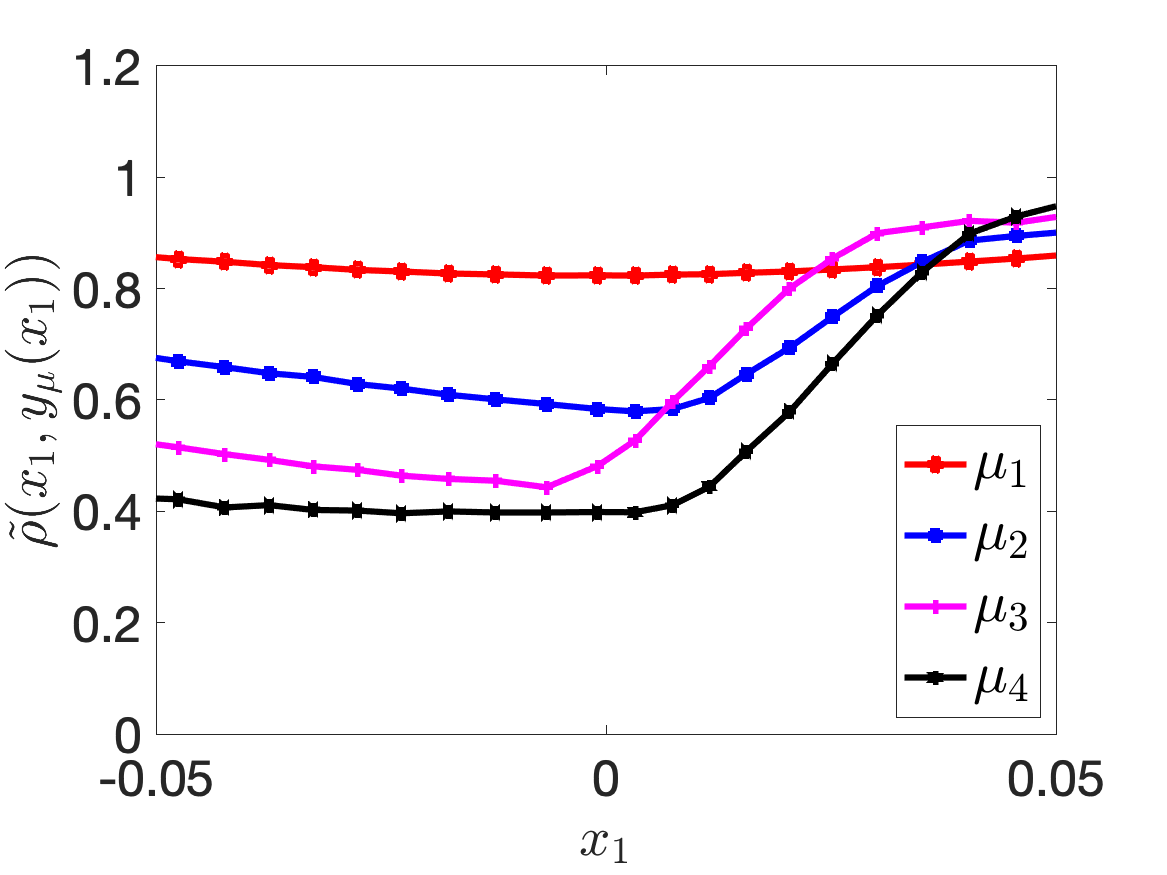}}
~~
 \subfloat[$it=3$] 
{  \includegraphics[width=0.33\textwidth]
 {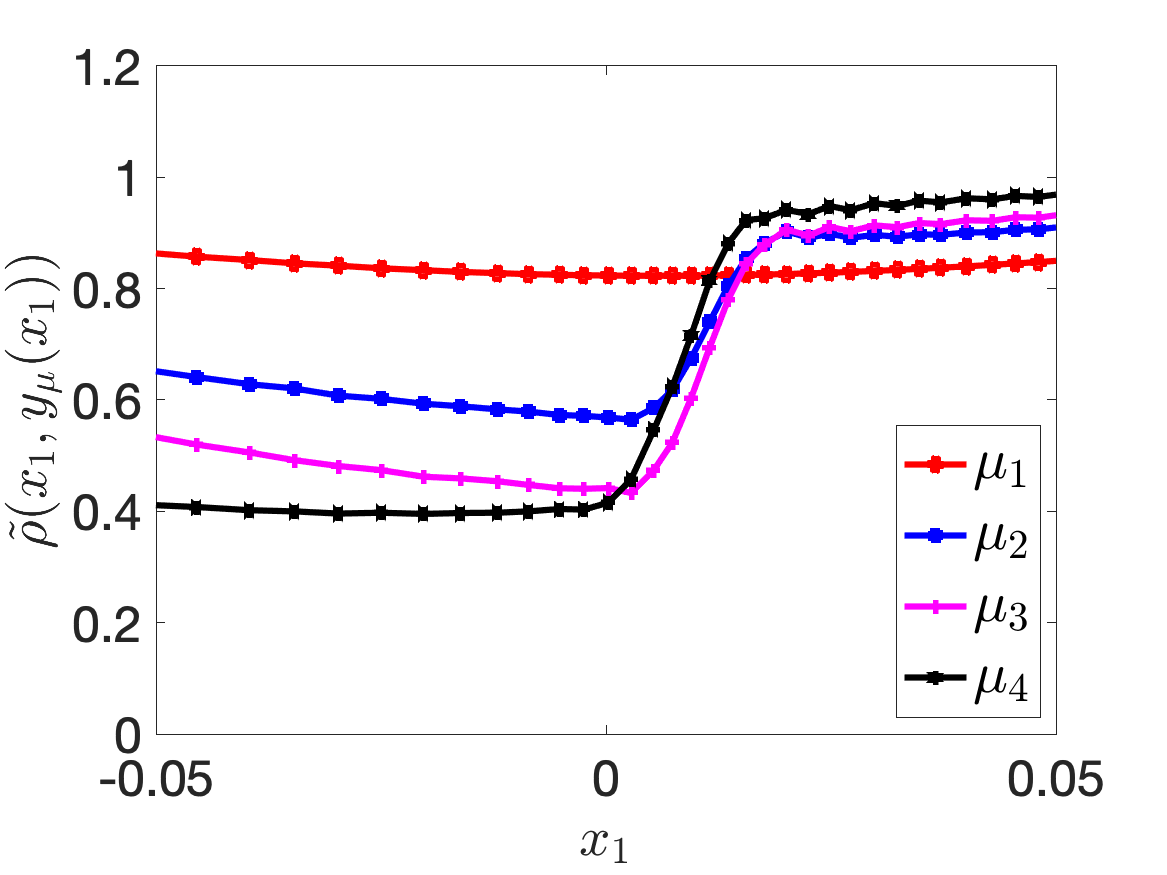}}
 
 \caption{transonic bump. 
Behavior of the (modified) density field  in physical (cf. (a)-(b)-(c)) and reference 
(cf. (d)-(e)-(f))   configuration for four values of the parameter and three iterations of the adaptive algorithm (basic version).
}
 \label{fig:transbump_basic_registration}
 \end{figure}

Figure \ref{fig:transbump_basic_registration} shows the behavior of the density field over the bump in physical and reference configuration for four values of the parameter and three iterations of the adaptive algorithm. We clearly notice the effect of the registration to nearly ``freeze'' the position of the shock --- when present --- in the reference configuration. We also notice that mesh adaptation is effective to sharpen the approximation of the shock as we increase the size of the mesh.

%transbump_registration

\begin{figure}[H]
\centering
 \subfloat[$it=1$] 
{  \includegraphics[width=0.33\textwidth]
 {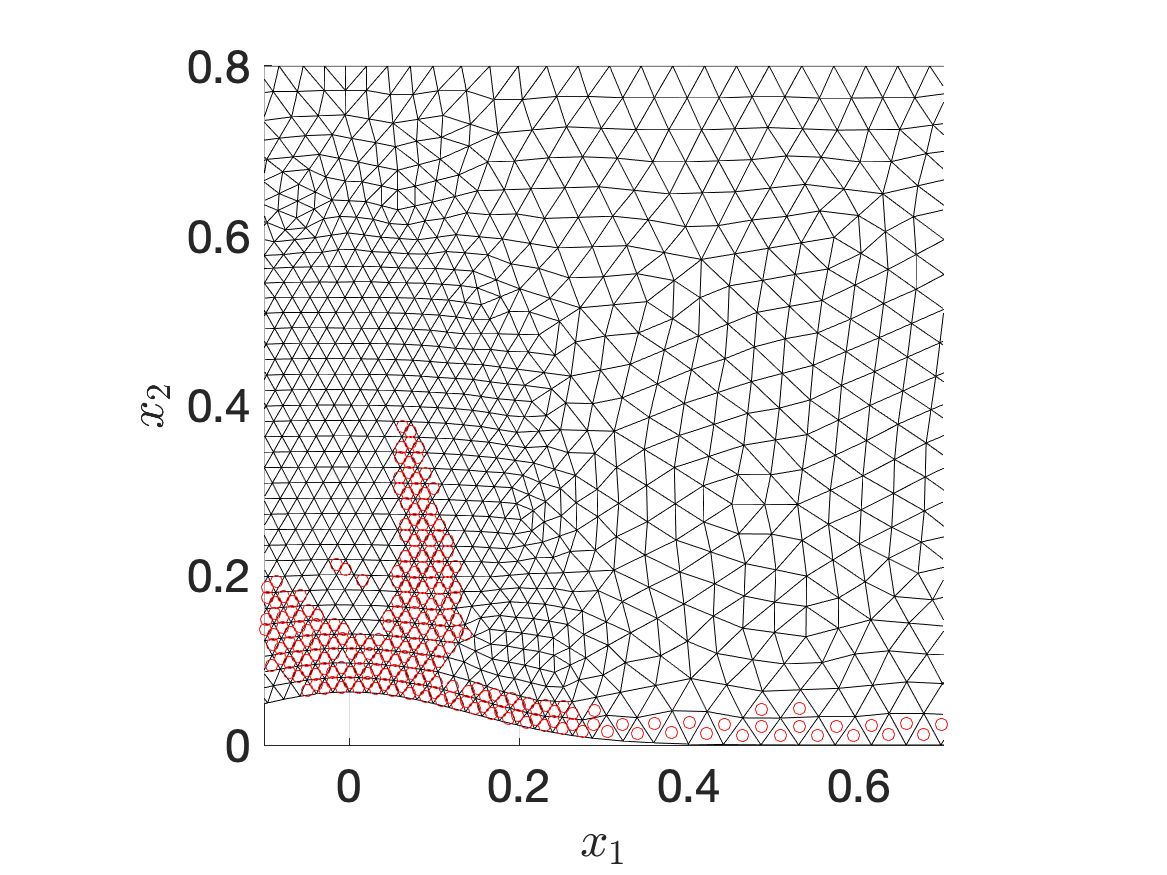}}
   ~~
 \subfloat[$it=2$] 
{  \includegraphics[width=0.33\textwidth]
 {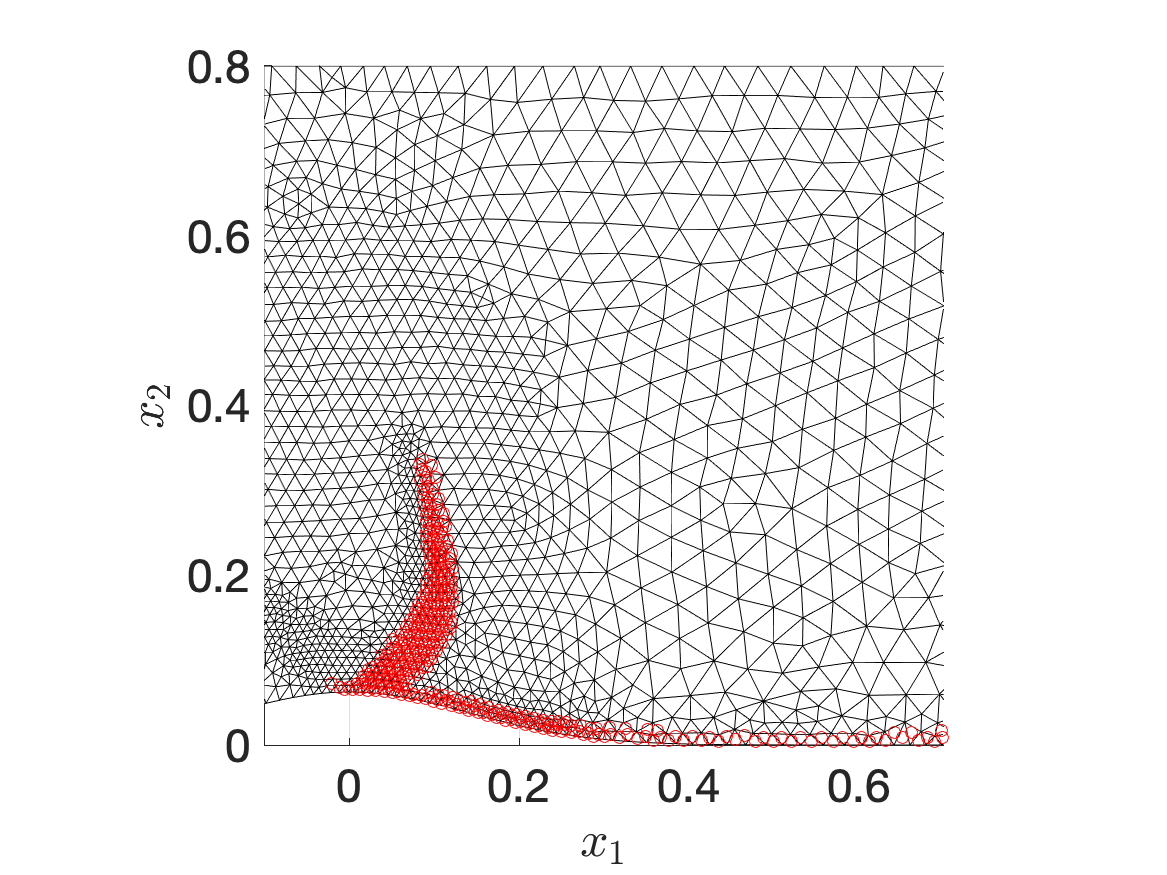}}
~~
 \subfloat[$it=3$] 
{  \includegraphics[width=0.33\textwidth]
 {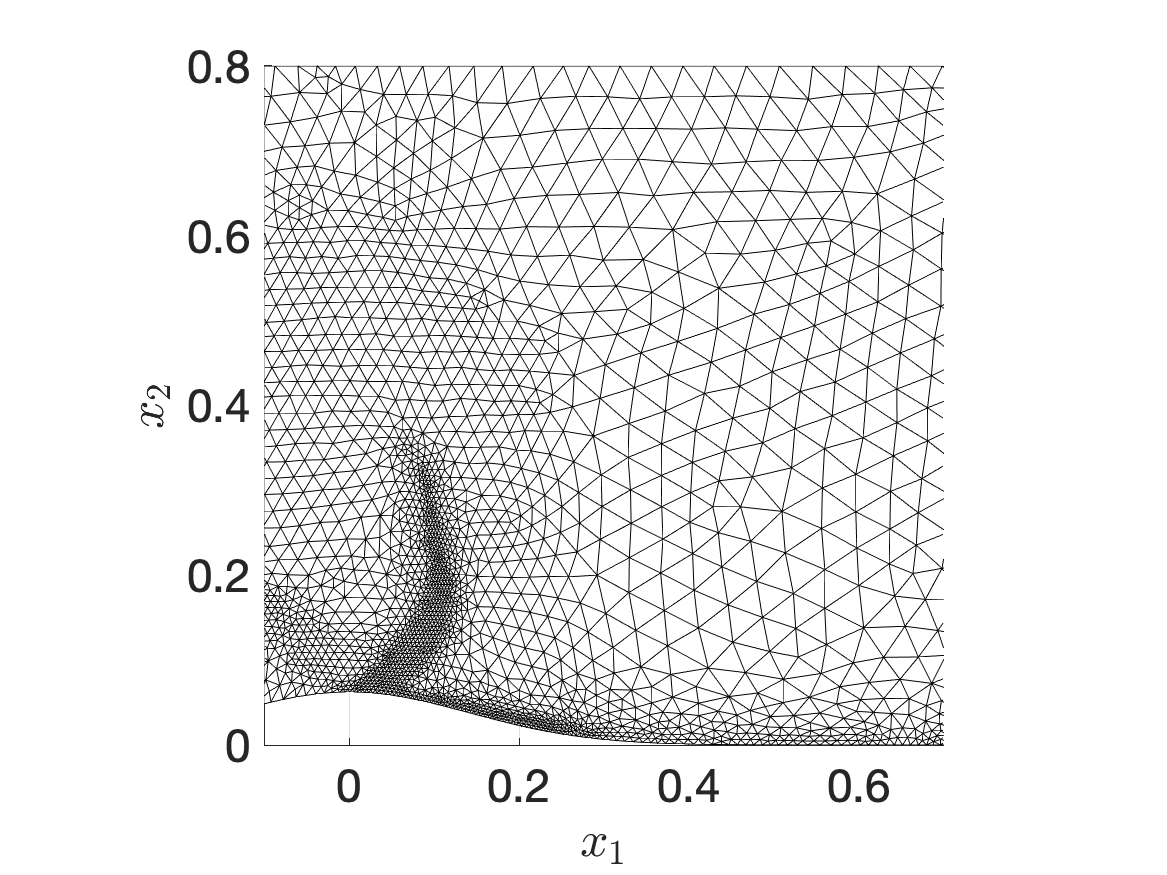}}
 
 \caption{transonic bump. Visualization of the reference mesh in the proximity of the bump, for three iterations of the adaptive algorithm (basic). 
 Red dots indicate the centers of the marked elements.
}
 \label{fig:transbump_basic_mesh}
 \end{figure}

Figure \ref{fig:transbump_basic_mesh} shows the reference mesh   in the proximity of the bump, for three iterations of the adaptive algorithm; 
red dots indicate the centers of the marked elements at iterations $it=2$ and $it=3$. Interestingly, we observe that the mesh is adapted in the proximity of the shock and in the proximity of the lower wall, in the area  downstream of the bump: as for the previous example, registration facilitates the task of parametric mesh adaptation by 
``freezing''
the coherent flow structure in the reference domain.

%transbump_meshadapt

Table \ref{tab:offline_costs_transbump} provides an overview of the offline costs. We notice that the costs are dominated by  snapshot generation at iteration $it=1$ --- which involves the construction of the ROM --- and by the weak-greedy algorithms.
We also observe that in our implementation the overhead costs of the greedy method are significant, while the cost of mesh adaptation is completely negligible.
For the registration algorithm, we distinguish between the cost to estimate the sensors 
$\{ s_{\mu}^{\rm hf} : \mu\in\mathcal{P}_{\rm train} \}$ 
in \eqref{eq:target_transonic_bump} (cf. Remark \ref{remark:implementation_registration}) and the cost of solving the parametric registration problem:  the former involves mesh interpolation over a curved HF mesh and is embarrassingly parallel, while the latter is dominated by the solution to the optimization problems \eqref{eq:optimization_all_params} for the first iteration.

\begin{table}[H]
\centering
\begin{tabular}{||c| c| c| c||} 
\hline 
\hline 
& it = 1 & it = 2 &  it = 3 \\[0.5ex]  \hline \hline 
ROB size: &  $14$ &  $17$ & $19$ \\[0.5ex] \hline 
mesh size: & $3448$ &  $4653$ & $6667$ \\[0.5ex] \hline 
snapshot generation: & $1459.37$ & $43.86$ & $59.59$ \\[0.5ex]   \hline 
registration (sensor def.): & $197.44$ & $202.13$ & $216.09$ \\[0.5ex]  \hline 
registration (optimization):& $402.93$ & $508.71$ & $770.13$ \\[0.5ex]  \hline 
mesh adaptation:& $0.00$ & $0.35$ & $0.57$ \\[0.5ex] \hline 
greedy alg (HF solves):& $762.15$ & $1539.81$ & $3961.10$ \\[0.5ex] \hline 
greedy alg (overhead): & $295.36$ & $632.19$ & $1200.36$ \\[0.5ex] \hline 
%PTC iterations (avg): & $18.64$ & $20.53$ & $29.00$ \\[0.5ex] %\hline 
\end{tabular}
\caption{transonic bump. Offline training costs (in seconds) of the adaptive (basic) approach.}
\label{tab:offline_costs_transbump}
\end{table}

\subsubsection{Acceleration of training through multi-fidelity strategies}
\label{sec:transonic_bump_numerics_accelerated}

We investigate the effect of the acceleration strategy discussed in section \ref{sec:adaptive_training}; to facilitate the comparison with the results of the previous section, we consider isotropic mesh adaptation. First, Table \ref{tab:initialization_PTC} investigates the effect of the initialization strategy on the convergence of the HF solver; the computational cost includes the interpolation cost.
We observe that our initialization strategy reduces the number of iterations required for convergence by roughly a factor three and computational costs by roughly a factor three for the final iteration.

%table comparison cost
\begin{table}%[H]
\centering
\begin{tabular}{|l|c|c|c|c|c|c|} 
\hline 
& \multicolumn{2}{|c|}{$it=1$} & \multicolumn{2}{|c|}{$it=2$} & \multicolumn{2}{|c|}{$it=3$} \\
\hline  
& basic & acc. & basic & acc. & basic & acc. \\[0.5ex]  
\hline 
avg nbr its & $18.64$  & $8.28$  & $20.53$  & $8.71$ & $29.00$  & $8.62$  \\ 
\hline  
avg cost & $54.44$  & $25.57$  & $90.58$  & $42.66$ & $208.48$  & $68.26$  \\ 
\hline 
\end{tabular}
\medskip
\caption{performance  of the HF solver for two initialization strategies.}
\label{tab:initialization_PTC}
\end{table}

Table \ref{tab:adaptive_overview} provides an overview of the two approaches in terms of the two metrics 
 \eqref{eq:L2error} and \eqref{eq:total_enthalpy_error}.
 Note that the acceleration strategy reduces offline costs by roughly $30\%$ mostly due to the reduction of the cost of the HF solves; it also slightly reduces the online costs by providing a more accurate initialization for GNM. Further numerical investigations
are provided in Appendix \ref{sec:appendix_transbump}.
We insist that the current implementation does not exploit parallel computing: since the acceleration strategy enables a much more efficient parallelization (cf. section \ref{sec:adaptive_training}), we expect more significant gains for the accelerated procedure when  combined with  parallel computing.

\begin{table}[H]
\centering
\begin{tabular}{|l || c|c|c|| c|c|c||} \hline &  \multicolumn{3}{c||}{$L^2$ error (avg)} &  \multicolumn{3}{|c||}{enthalpy error (avg)} \\ \cline{2-7}  
& 1 & 2 &  3 & 1 & 2 & 3  \\ \hline  
Basic & $0.53 \cdot 10^{-3}$ & $0.44 \cdot 10^{-3}$ & $0.48 \cdot 10^{-3}$ & $0.58 \cdot 10^{-3}$ & $0.35 \cdot 10^{-3}$ &  $0.23 \cdot 10^{-3}$  \\ \hline
Accelerated & $0.29 \cdot 10^{-3}$ & $0.31 \cdot 10^{-3}$ & $0.49 \cdot 10^{-3}$ & $0.59 \cdot 10^{-3}$ & $0.37 \cdot 10^{-3}$ &  $0.23 \cdot 10^{-3}$  \\ \hline
\end{tabular} 

\medskip 
\begin{tabular}{|l || c|c|c|| c|c|c|| c|} \hline &  \multicolumn{3}{c||}{ROB size} &  \multicolumn{3}{|c||}{online cost (avg)} &  offline cost  \\ \cline{2-8}  
& 1 & 2 &  3 & 1 & 2 & 3  & \\ \hline 
Basic & $14$    & $17$   & $19$  & $0.20$  & $0.28$ & $0.36$& 03:24:13 \\ \hline 
Accelerated & $18$    & $17$   & $21$  & $0.25$  & $0.25$ & $0.36$& 02:11:47 \\ \hline 
\end{tabular} 

%\medskip
\caption{Comparison of the performance of the basic and accelerated adaptive procedures.}
\label{tab:adaptive_overview}
\end{table}

\begin{figure}[H]
\centering
\includegraphics[width=0.33\textwidth]
 {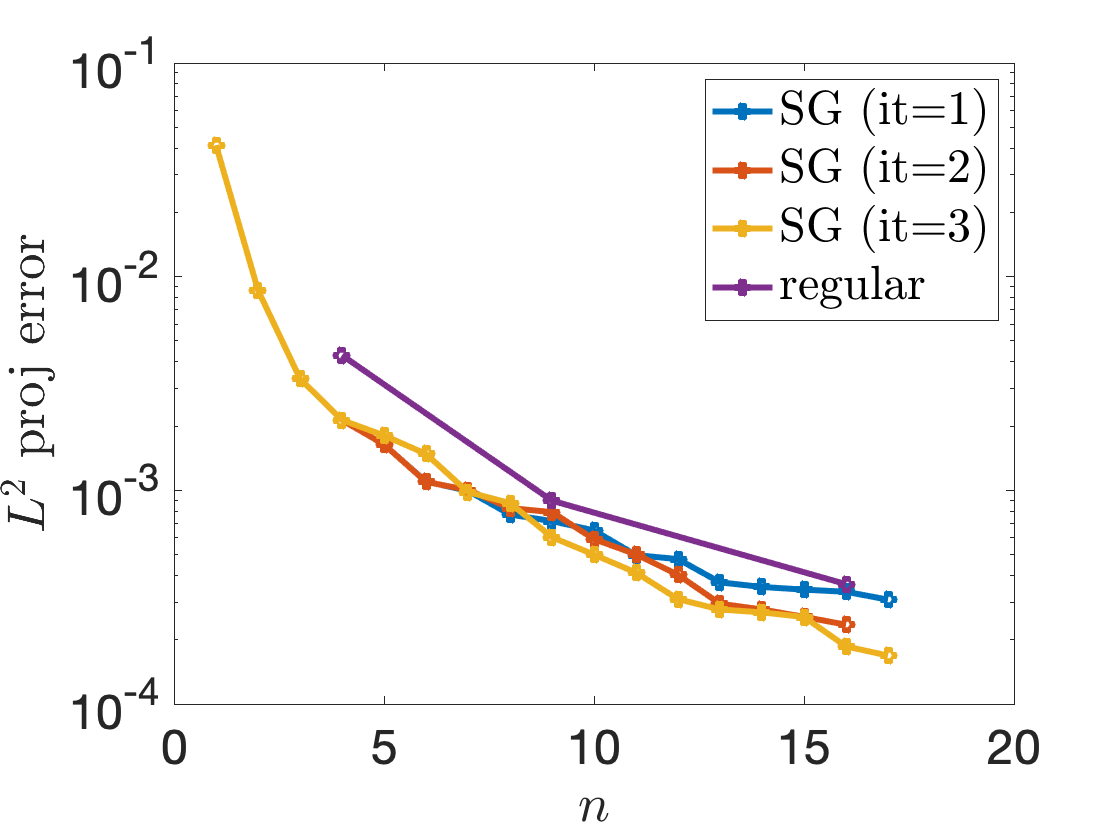} 
  
 \caption{transonic bump. Performance of greedy sampling based on datasets of different fidelity over the test set.
}
 \label{fig:transbump_greedy_sampling}
 \end{figure}

 Figure \ref{fig:transbump_greedy_sampling} investigates the performance of the greedy strategy.
 We perform the strong-greedy algorithm on the snapshot sets generated at iterations one, two and three to identify the ``optimal'' parameters
 $\mathcal{P}_{\star}^{it,n} =\{\mu^{\star, it,i} \}_{i=1}^n$ for $it=1,2,3$. Then, we compute the projection error
 $$
 E_{\rm proj, \mu}
 :=\frac{ \min_{\zeta \in   \mathcal{Z}_{\star}^{it,n}  }  \| \widetilde{q}_{\mu}^{\rm hf} - \zeta   \|    }{
\| \widetilde{q}_{\mu}^{\rm hf}   \| 
 },
 \quad
 {\rm where} \;\;
  \mathcal{Z}_{\star}^{it,n}  = 
  {\rm span} \left\{ \widetilde{q}_{\mu}^{\rm hf}  \; : \;
  \mu \in \mathcal{P}_{\star}^{it,n}
      \right\}
 $$
over the test set of $n_{\rm test}=20$ simulations;
here, $ \widetilde{q}_{\mu}^{\rm hf}$ refers to the HF estimate obtained using the DG model at the final ($it=3$)  iteration.
To provide a concrete reference, we compare the results obtained using  regular grids  ($2\times 2, 3\times 3, 4\times 4$) of parameters. We observe that the strong greedy algorithm based on iteration $it=2$ provides results that are nearly as good as the results obtained based on the snapshot set of  iteration $it=3$. This empirical finding suggests that the application of the strong-greedy method to a lower-fidelity snapshot set might provide an inexpensive yet effective  sampling strategy for model reduction.

\subsubsection{Accelerated training with anisotropic mesh adaptation}
\begin{figure}[H]
\centering
 \subfloat[$it=1$] 
{  \includegraphics[width=0.33\textwidth]
 {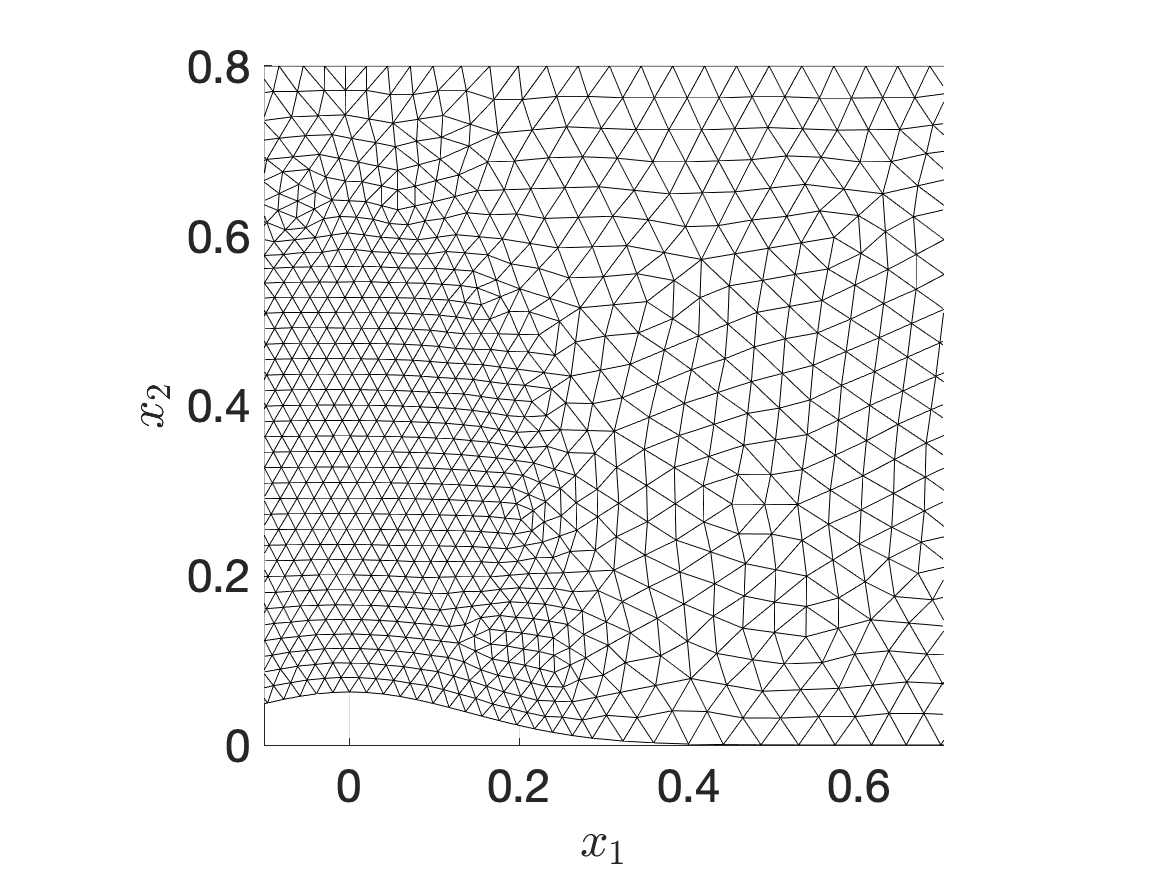}}
   ~~
 \subfloat[$it=2$] 
{  \includegraphics[width=0.33\textwidth]
 {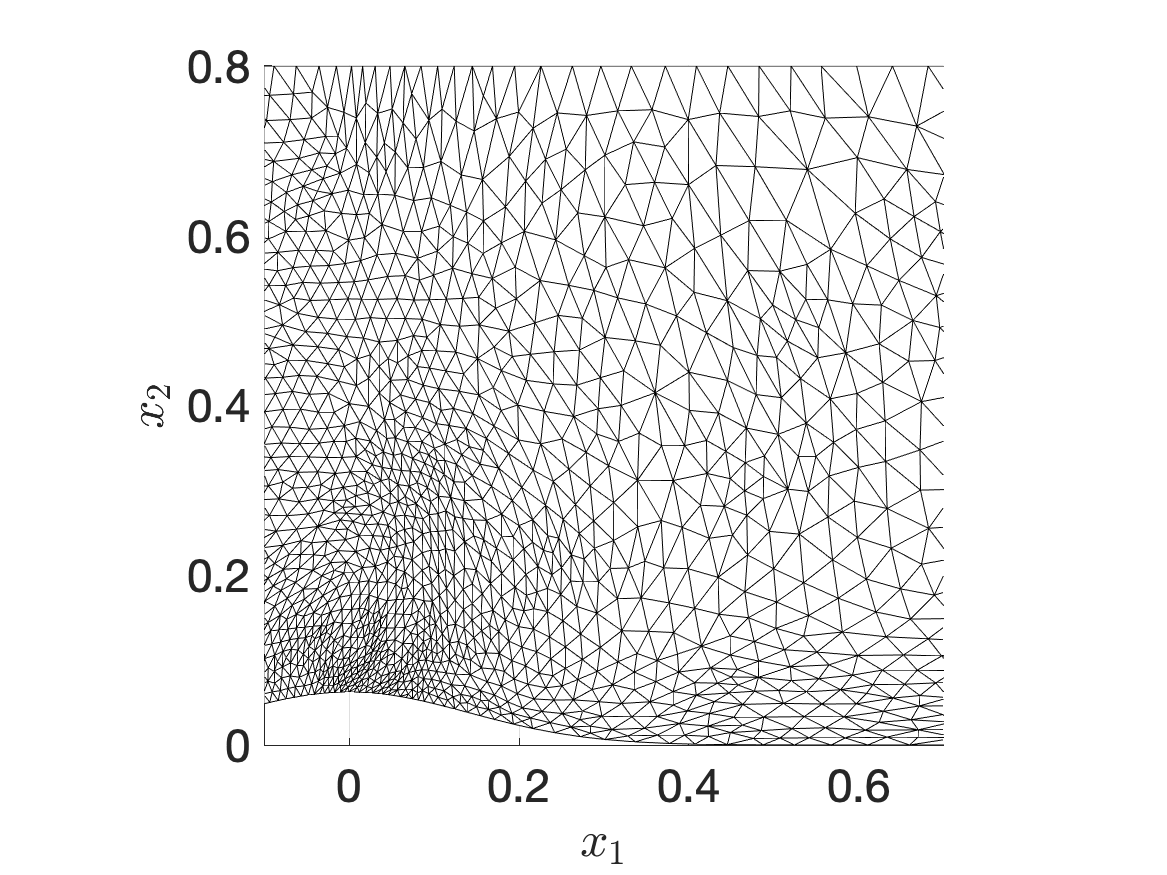}}
~~
 \subfloat[$it=3$] 
{  \includegraphics[width=0.33\textwidth]
 {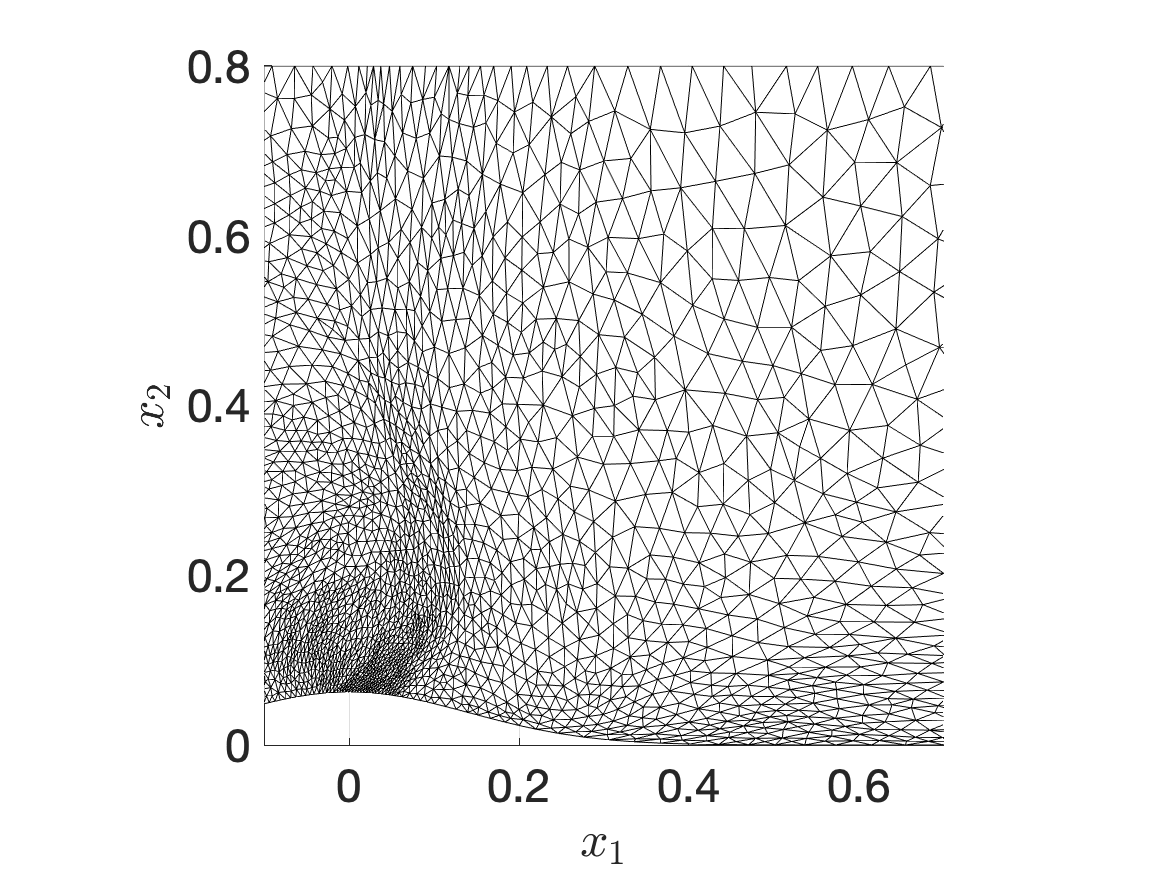}}
 
 \caption{transonic bump. Visualization of the reference mesh in the proximity of the bump, for three iterations of the adaptive algorithm (accelerated) with anisotropic mesh adaptation
 ($N_{\rm e}=3448$,   $N_{\rm e}=4440$, $N_{\rm e}=6389$).
}
 \label{fig:transbump_aniso_mesh}
 \end{figure}

We execute three iterations of Algorithm \ref{alg:offline_training} with anisotropic mesh adaptation (cf. section \ref{sec:mesh_adaptation}). We initially set the parameter $N$ in \eqref{eq:hessian_based_metric} equal to $750$ and we increase it at each iteration by a factor $1.5$. 
Figure \ref{fig:transbump_aniso_mesh} shows the sequence of meshes generated by Algorithm \ref{alg:offline_training}. We observe that the meshes are nearly isotropic in the proximity of the shock while they exhibit elongated  elements in the downstream region (the minimum  radius ratio is roughly $0.05$).
We notice 
that the adapted mesh for a single field is significantly more anisotropic in the proximity of the shock, but it becomes less and less anisotropic as we combine metrics associated with different parameters.
This is likely due to the fact that the shock is not sharply tracked in the reference configuration.

%transbump_meshadapt

Table \ref{tab:adaptive_overview_aniso} compares performance of the accelerated training strategy based on isotropic and anisotropic MA: we observe that the two approaches lead to comparable performance for this model problem.
We notice that the HF model requires slightly more PTC iterations to converge for anisotropic meshes: the difference is much more significant when we initialize the solver with  the free-stream flow. This observation shows the importance of exploiting prior information to  properly initialize the HF solver.

\begin{table}[H]
\centering
\begin{tabular}{|l||c|c|c|| c|c|c||}
\hline
& \multicolumn{3}{c||}{$L^2$ error (avg)} &
 \multicolumn{3}{c||}{enthalpy error (avg)}
 \\ \cline{2-7}
 &1 & 2 &  3 & 1 & 2 & 3 \\ \hline
 Isotropic MA & $0.29 \cdot 10^{-3}$ & $0.31 \cdot 10^{-3}$ & $0.49 \cdot 10^{-3}$
& $0.59 \cdot 10^{-3}$ & $0.37 \cdot 10^{-3}$ &  $0.23 \cdot 10^{-3}$  \\ \hline
 Anisotropic MA & $0.29 \cdot 10^{-3}$ & $0.36 \cdot 10^{-3}$ & $0.50 \cdot 10^{-3}$
& $0.59 \cdot 10^{-3}$ & $0.36 \cdot 10^{-3}$ &  $0.26 \cdot 10^{-3}$
  \\ \hline
\end{tabular}

\medskip 
\begin{tabular}{|l||c|c|c|| c|c|c|| c|}
\hline
& \multicolumn{3}{c||}{ROB size} &
 \multicolumn{3}{c||}{online cost (avg)}
 &
 offline cost 
 \\ \cline{2-8}
  &1 & 2 &  3 & 1 & 2 & 3  & \\ \hline
Isotropic MA & $18$    & $17$   & $21$  & $0.25$  & $0.25$ & $0.36$& 02:11:47 \\ \hline
Anisotropic MA & $18$    & $17$   & $16$  & $0.26$  & $0.26$ & $0.24$& 01:58:48 \\ \hline 
  \end{tabular}

%\begin{tabular}{|l || c|c|c|| c|c|c|| c|} \hline &  \multicolumn{3}{c||}{ROB size} &  \multicolumn{3}{|c||}{online cost (avg)} &  offline cost  \\ \cline{2-8}  
%

%\end{tabular} 

\caption{comparison of the performance  of the  accelerated adaptive procedures with isotropic and anisotropic mesh adaptation.}
\label{tab:adaptive_overview_aniso}
\end{table}

\section{Summary and discussion}
\label{sec:conclusions}
We developed and numerically validated an adaptive strategy for the simultaneous construction of high-fidelity and reduced-order approximations for parametric problems with discontinuous solutions.
The approach relies on registration to track moving features of the solution field, metric-based mesh adaptation to devise an accurate mesh for the solution over a range of parameters, and projection-based model reduction to effectively estimate the (mapped) solution field. We show that registration is key to improve the compressibility of the solution manifold
(cf. Figure \ref{fig:nozzle_advanced_plot})
and enables parsimonious yet accurate HF approximations by complementing parameter-independent $h$-adaptation with parameter-dependent $r$-adaptation
(cf. Figures \ref{fig:nozzle_advanced_plot} and 
\ref{fig:transbump_basic_mesh}).
We also show that our adaptive training 
strategy provides increasingly more accurate approximations of the solution field (cf. Figures \ref{fig:nozzle_basic_registration} and \ref{fig:transbump_basic_registration}) and can  be significantly accelerated by exploiting information from previous iterations (cf. section \ref{sec:adaptive_training} and Tables \ref{tab:adaptive_overview} and \ref{tab:adaptive_overview_aniso}).

We plan to extend our work in several directions.
First, we wish to apply our framework to a broad range of problems in nonlinear mechanics, viscous compressible flows, and hydraulics: towards this end, we should extend our approach to unsteady PDEs and we should devise effective mesh and registration sensors for a broad range of solution features of interest.
Second, we plan to leverage  clustering techniques
to  further improve the quality of the HF mesh:
even  if MA allows us to optimize the size of the mesh,  the HF meshes can still be large, in particular in the presence of  parameter-induced topology changes  that cannot be captured by a single parametric deformation; by resorting to  clustering techniques, we hence expect to better control the distribution of the degrees of freedom in the spatio-parametric space.
% Second, we plan to devise clustering techniques to deal with parameter-induced topology changes that cannot be captured by a single parametric deformation.

\section*{Acknowledgements}
The authors acknowledge the support  provided by  Inria through the exploratory action program (project title: \emph{Adaptive Meshes for Model Order Reduction}, AM2OR).
TT acknowledges the support by European Union’s Horizon 2020 research and innovation program  under the Marie Skłodowska-Curie Actions, grant agreement 872442 (ARIA).
%thank  Professor Masayuki Yano (University of Toronto) for fruitful discussions.

\appendix
 
\section{Further greedy procedures employed at training stage}
\label{sec:appendix_greedy}
In this section, we provide two greedy algorithms that are used during the execution of 
Algorithm \ref{alg:offline_training}. 
Algorithm \ref{alg:registration} summarizes the parametric registration  procedure for the transonic bump test case;
on the other hand, Algorithm \ref{alg:strong_greedy} outlines the strong greedy procedure employed to initialize Algorithm \ref{alg:weak_greedy} and to select a subset of relevant solutions for mesh adaptation (cf. section \ref{sec:mesh_adaptation}).

We use notation
$$
\left[\widehat{\mathbf{a}}_{\mu}, 
\mathfrak{f}_{\mu}^{\star}  \right]
\, = \,
\texttt{registration} \left(
s_{\mu}^{\rm hf}, \,   \mathcal{S}_n, \,  \texttt{W}_{\rm p}, \, \mathcal{T}_{\rm hf}, \, 
\Psi_{\mu}^{\rm geo}, \, 
\mathbf{a}_{\mu}^{0}
\right)
$$
to refer to the function that takes as inputs 
(i) the target sensor $s_{\mu}^{\rm hf}:\Omega_{\rm p}\to \mathbb{R}$, 
(ii) the template space $\mathcal{S}_n$,
(iii) the ROB $\texttt{W}_m$ associated with the 
mapping space $\mathcal{U}_{\rm p} \subset \mathcal{U}_{\rm hf,p}$, 
(iv) the HF mesh $\mathcal{T}_{\rm hf}$, 
(v) the geometric mapping $\Psi_{\mu}^{\rm geo}:\Omega_{\rm p}\to \Omega_{\mu}$ and 
(vi) the initial guess 
$ \mathbf{a}_{\mu}^{0} \in \mathbb{R}^m$
for the optimizer,
and returns 
(I) the mapping coefficients $\widehat{\mathbf{a}}_{\mu}$ associated with a  local minimum of the problem \eqref{eq:optreg_statement}, and 
(II) the value of the target function
$\mathfrak{f}_{\mu}^{\star}  = 
\mathfrak{f}_{\mu}^{\rm tg}(\widehat{\mathbf{a}}_{\mu})
$.
We also introduce the function
$$
[  \texttt{W}_m, \; \{  
\mathbf{a}_{\mu}^{\rm proj}
  \}_{\mu\in \mathcal{P}_{\rm train}} ]  =
\texttt{POD} \left( 
\{ 
\widehat{\mathbf{a}}_{\mu}
  \}_{\mu\in \mathcal{P}_{\rm train}}, 
tol_{\rm pod}  ,
(\cdot, \cdot)_2
 \right),
$$
which implements POD
based on the method of snapshots with Euclidean inner product $(\cdot,\cdot)_2$:
the tolerance 
$tol_{\rm pod}>0$ drives 
the selection of the number of modes $m$ based on the energy criterion 
\begin{equation}
\label{eq:POD_cardinality_selection}
m := \min \left\{
m': \, \sum_{j=1}^{m'} \lambda_{j} \geq  \left(1 - tol_{\rm pod} \right) 
\sum_{i=1}^{n_{\rm train}} \lambda_i
\right\},
\end{equation} 
where $ \lambda_1\geq \ldots \geq \lambda_{n_{\rm train}}\geq 0$ are the eigenvalues of the Gramian matrix $\mathbf{C}\in \mathbb{R}^{n_{\rm train}\times n_{\rm train}}$ such that
$(\mathbf{C})_{k,k'}
= 
\mathbf{a}_{\mu^k}^{\star} \cdot 
\mathbf{a}_{\mu^{k'}}^{\star} 
$. 
The function \texttt{POD} returns also the mapping coefficients associated with the 
projected displacements
$\mathbf{a}_{\mu}^{\rm proj}$ onto the POD space;
the latter  are used to initialize the iterative method for the optimization problem in the subsequent iterations.

\begin{algorithm}[H]                      
\caption{: registration algorithm (\cite{taddei2021space}).}     
\label{alg:registration}     

 \small
\begin{flushleft}
\emph{Inputs:}  $\{  s_{\mu} : \mu\in 
\mathcal{P}_{\rm train} \}$ snapshot set, 
$\mathcal{S}_{n_0} = {\rm span} \{ 
 s_{\mu^{\star,(i)}}^{\rm hf} 
\}_{i=1}^{n_0}$ initial template space;
$\mathcal{T}_{\rm hf}$ mesh for HF computations.
\smallskip

\emph{Outputs:} 
${\mathcal{S}}_n $ template space, 
$\texttt{W}_m:\mathbb{R}^m \to    \mathcal{U}_{\rm p}
$ mapping ROB,
$\{  \varphi_{\rm p,\mu^k}^{\star} 
= \texttt{W}_m \mathbf{a}_{\mu^k}^{\star}
 \}_k$ optimal  mappings.
\end{flushleft}                      

 \normalsize 

\begin{algorithmic}[1]
\State
Initialization: 
$\mathcal{S}_{n=n_0} = \mathcal{S}_{n_0}$,
$\Xi_{\star} = \{\mu^{\star,(i)} \}_{i=1}^{n_0}$,
$\mathcal{U}_{\rm p} =\mathcal{U}_{\rm hf,p}$.
\vspace{3pt}

\For {$n=n_0, \ldots, n_{\rm max}-1$ }

\State
$
\left[\widehat{\mathbf{a}}_{\mu}, 
\mathfrak{f}_{\mu}^{\star}  \right]
\, = \,
\texttt{registration} \left(
s_{\mu}^{\rm hf}, \,   \mathcal{S}_n, \,  \texttt{W}_{\rm p}, \, \mathcal{T}_{\rm hf}, \, \mathbf{a}_{\mu}^{0}
\right)
$ for all $\mu \in \mathcal{P}_{\rm train}$,

\hfill
\emph{see sections \ref{sec:param_reg} and \ref{sec:adaptive_training} for definition of $\mathbf{a}_{\mu}^{0}$}
\vspace{3pt}

\State
$[\texttt{W}_m, \; \{  
\mathbf{a}_{\mu}^{\rm proj}
  \}_{\mu\in \mathcal{P}_{\rm train}} ]  =
\texttt{POD} \left( 
\{ \widehat{\mathbf{a}}_{\mu}  \}_{\mu\in \mathcal{P}_{\rm train}}, 
tol_{\rm pod}  ,
(\cdot, \cdot)_2
 \right),$
\vspace{3pt}

\If{$\max_{\mu\in \mathcal{P}_{\rm train}}  \mathfrak{f}_{\mu}^{\star}   < \texttt{tol}$}, \texttt{break}

\Else

\State
$\Xi_{\star} = 
\Xi_{\star} \cup 
\{ \mu^{\star,(n+1)} \}$
with $\mu^{\star, (n+1)}= {\rm arg} 
\max_{\mu\in \mathcal{P}_{\rm train}}  
\mathfrak{f}_{\mu}^{\star} 
$.
\vspace{3pt}

\State
$\mathcal{S}_{n+1} =
 {\rm span} 
\{ 
s_{\mu^{i,\star}}^{\rm hf}      
\circ   \Phi_{\rm p,\mu^{i,\star}}   
 \}_{i=1}^{n+1} 
 $.
 \vspace{3pt}
 
 \EndIf
\EndFor
\end{algorithmic}
\end{algorithm}

\begin{algorithm}[H]                      
\caption{: strong-greedy algorithm (see, e.g., \cite[section 7.3]{quarteroni2015reduced}). }     
\label{alg:strong_greedy}     

 \small
\begin{flushleft}
\emph{Inputs:}  $\{ \widehat{\boldsymbol{\alpha}}_{\mu} :   \mu \in  \mathcal{P} _{\rm train} \} \subset \mathbb{R}^n$ snapshot set, 
$n_0 \leq n$ size of the desired reduced space.
\smallskip

\emph{Outputs:} 
$\mathcal{P}_{\star} =  \{ \mu^{\star, i}  \}_{i=1}^{n_0}$ selected parameters.

\end{flushleft}                      

 \normalsize 

\begin{algorithmic}[1]
\State
Choose 
$\mathcal{Z} = \emptyset$, 
$\mathcal{P}_{\star} = \emptyset$.
\vspace{3pt}

\For {$i= 1, \ldots, n_0$}

\State
Compute $\mu^{\star, i} = {\rm arg} \max_{\mu\in \mathcal{P} _{\rm train}  } \min_{\boldsymbol{\alpha}\in \mathcal{Z}} \|  \boldsymbol{\alpha} - \widehat{\boldsymbol{\alpha}}_{\mu}  \|_2$
\vspace{3pt}

\State
Update
$ \mathcal{Z} =  \mathcal{Z} \cup {\rm span} \{  \widehat{\boldsymbol{\alpha}}_{\mu^{\star, i}}  \}$ and
$ \mathcal{P}_{\star} =  \mathcal{P}_{\star}  \cup  \{ \mu^{\star, i}  \}$
\vspace{3pt}
 
\EndFor
\end{algorithmic}
\bigskip

\end{algorithm}

\section{Further numerical results for the transonic bump}
\label{sec:appendix_transbump}
In this section, we provide detailed results of the accelerated iterative procedure discussed in section  
\ref{sec:adaptive_training}. 
We distinguish between results obtained using isotropic and anisotropic mesh adaptation.
 
\subsection{Acceleration with isotropic mesh adaptation} 
\begin{figure}[H]
\centering
 \subfloat[] 
{  \includegraphics[width=0.4\textwidth]
 {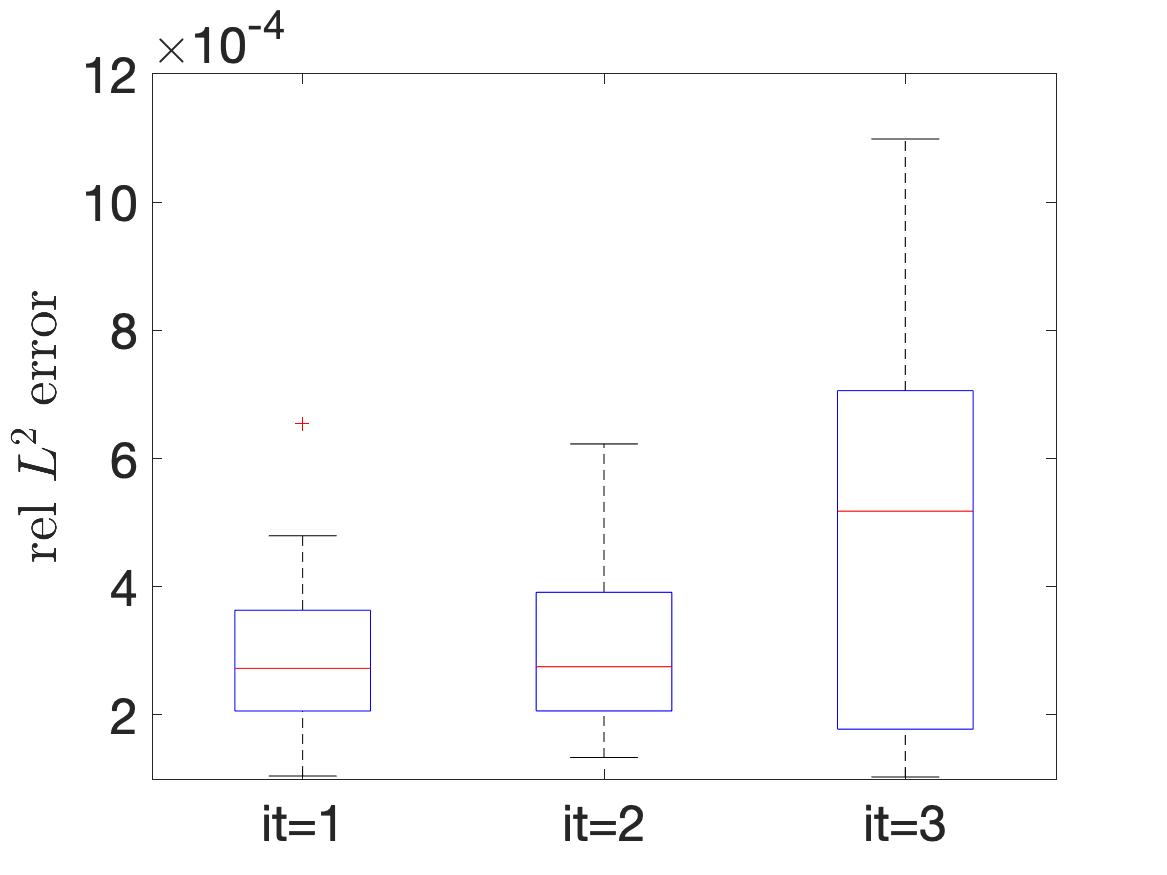}}
   ~~
 \subfloat[] 
{  \includegraphics[width=0.4\textwidth]
 {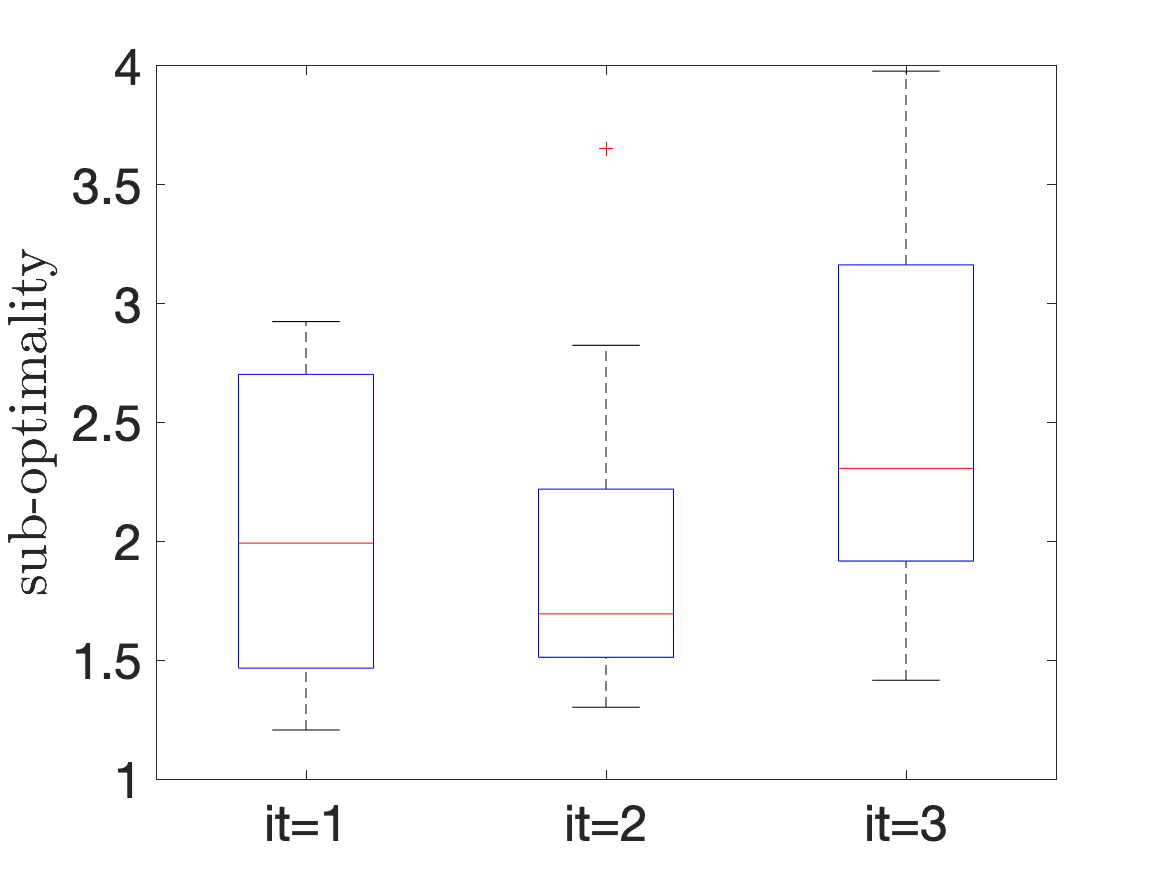}}

 \subfloat[] 
{  \includegraphics[width=0.4\textwidth]
 {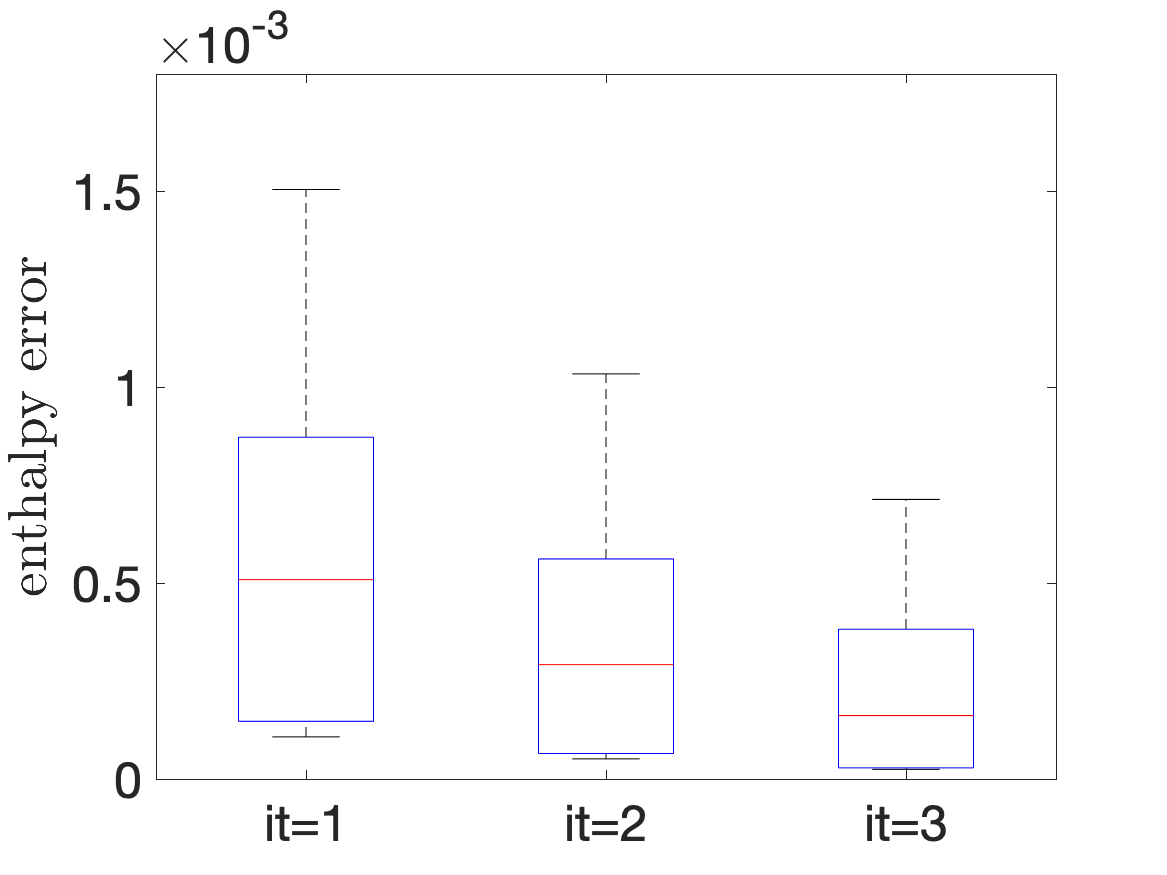}}
 ~~
  \subfloat[] 
{  \includegraphics[width=0.4\textwidth]
 {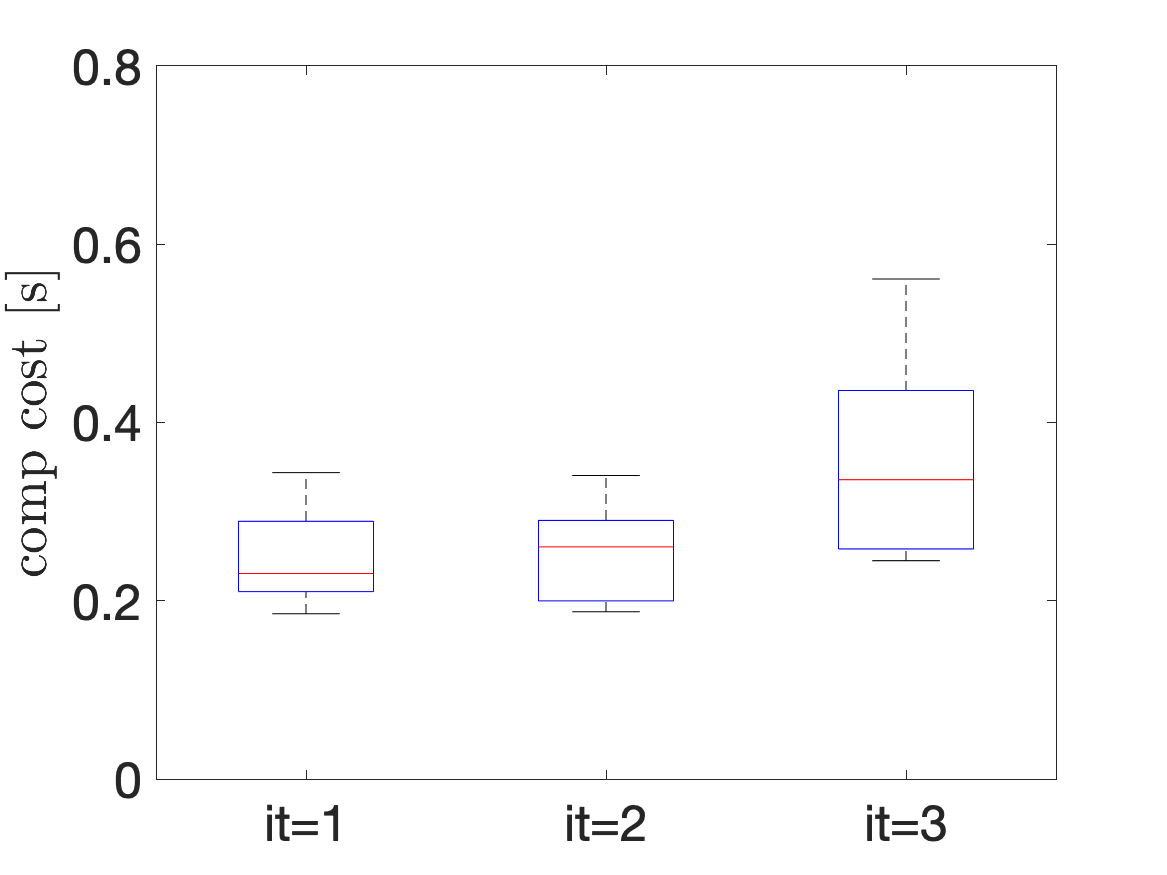}}
 
 \caption{transonic bump. Performance of the ROM for three iterations of the adaptive (accelerated) procedure with isotropic mesh adaptation.
}
 \label{fig:transbump_accelerated_rom}
 \end{figure}

\begin{figure}[H]
\centering
 \subfloat[$it=1$] 
{  \includegraphics[width=0.33\textwidth]
 {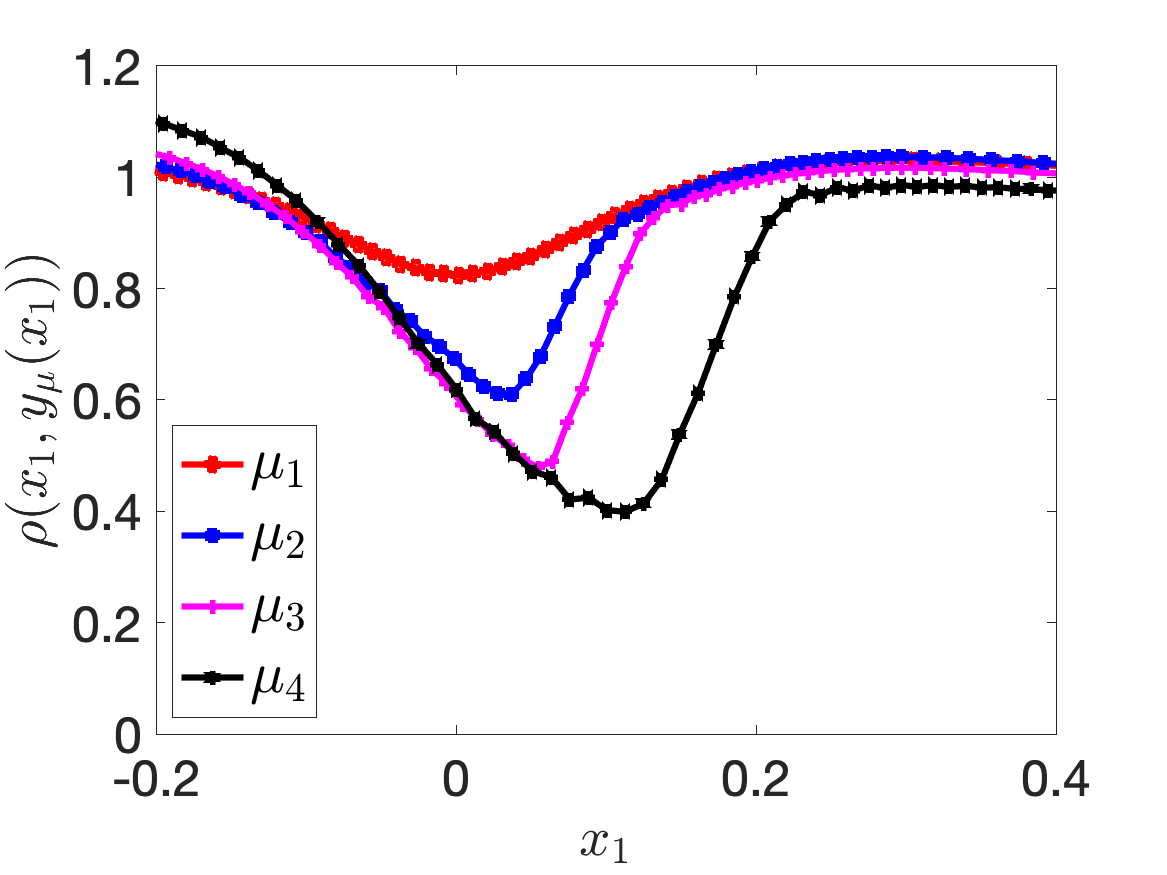}}
   ~~
 \subfloat[$it=2$] 
{  \includegraphics[width=0.33\textwidth]
 {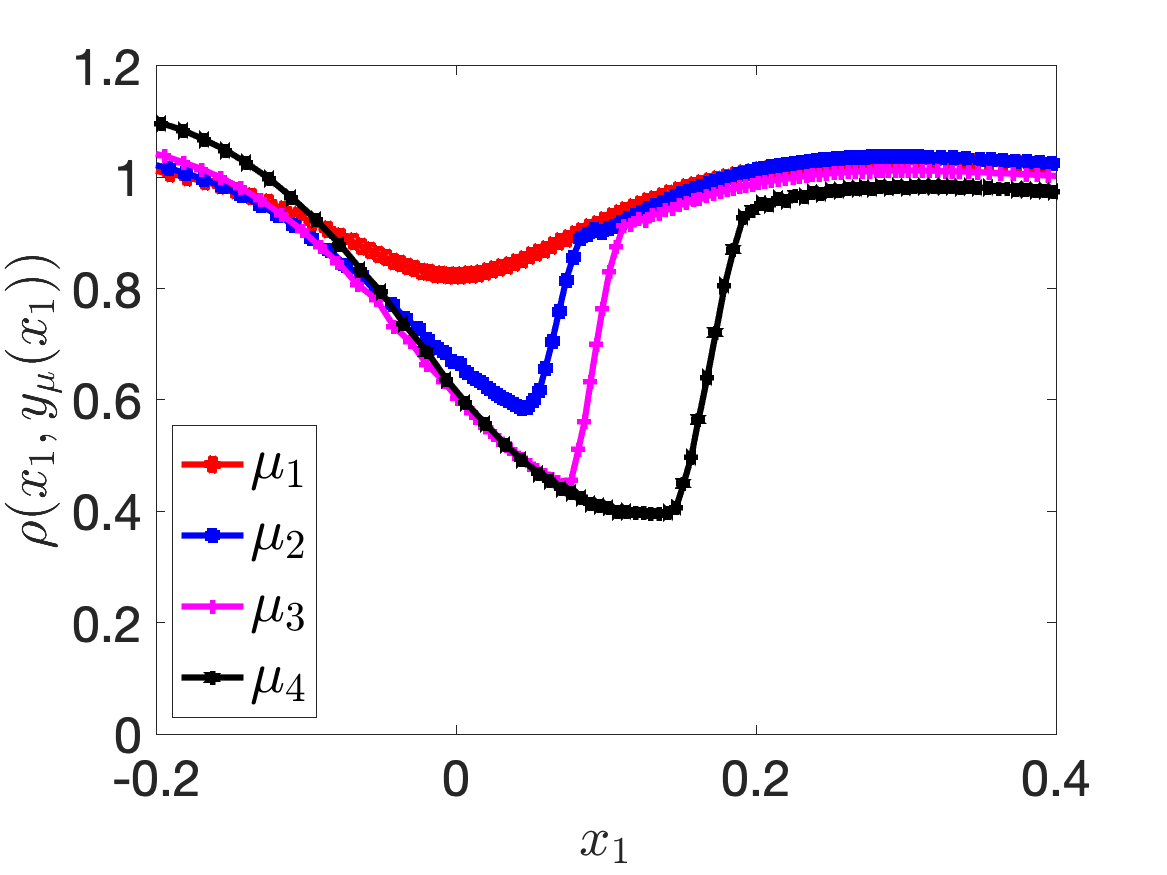}}
~~
 \subfloat[$it=3$] 
{  \includegraphics[width=0.33\textwidth]
 {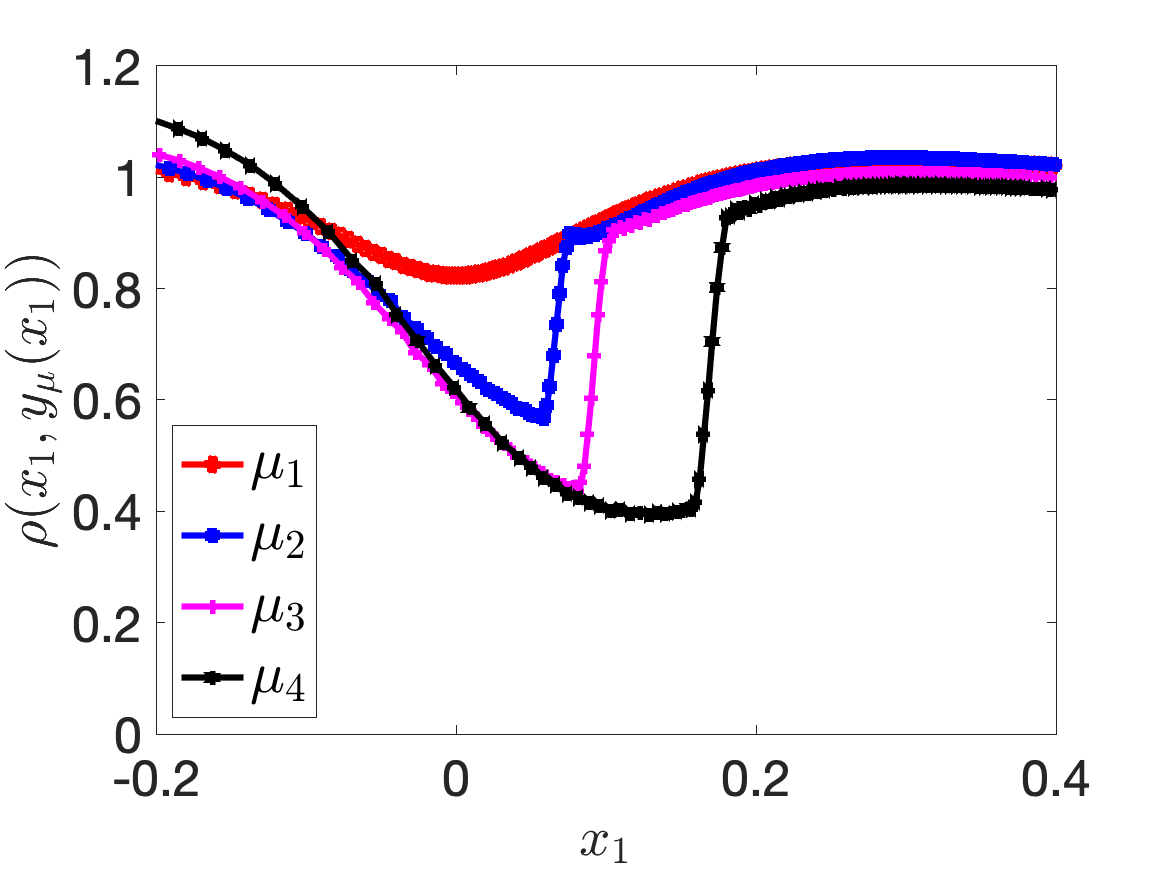}}
 
 \subfloat[$it=1$] 
{  \includegraphics[width=0.33\textwidth]
 {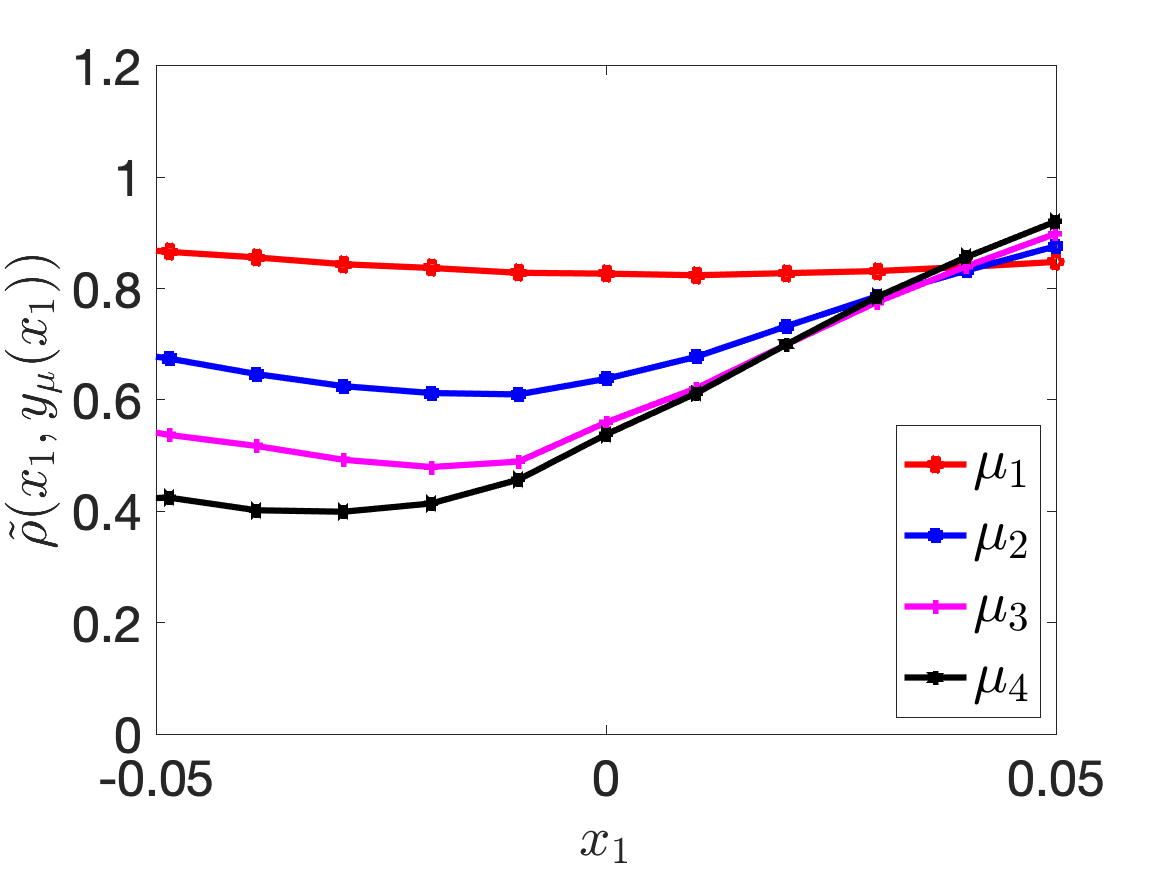}}
   ~~
 \subfloat[$it=2$] 
{  \includegraphics[width=0.33\textwidth]
 {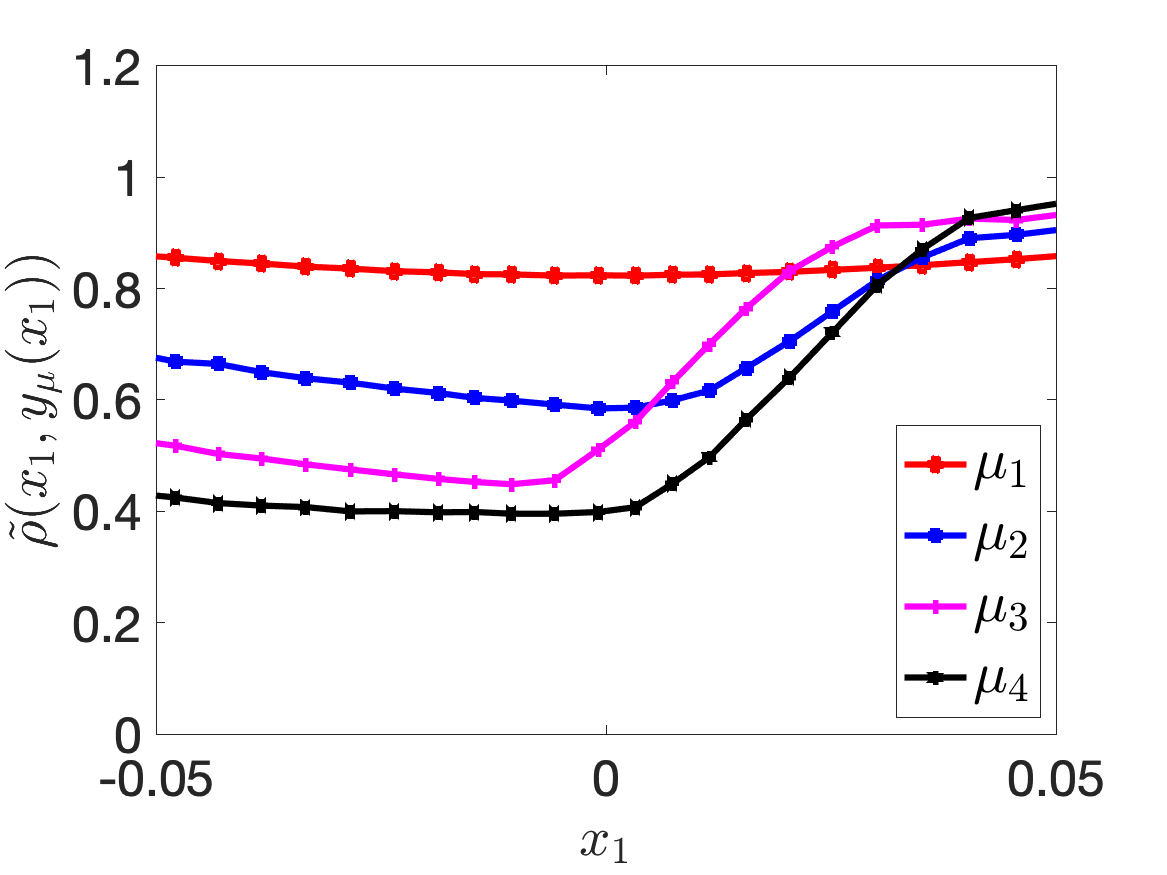}}
~~
 \subfloat[$it=3$] 
{  \includegraphics[width=0.33\textwidth]
 {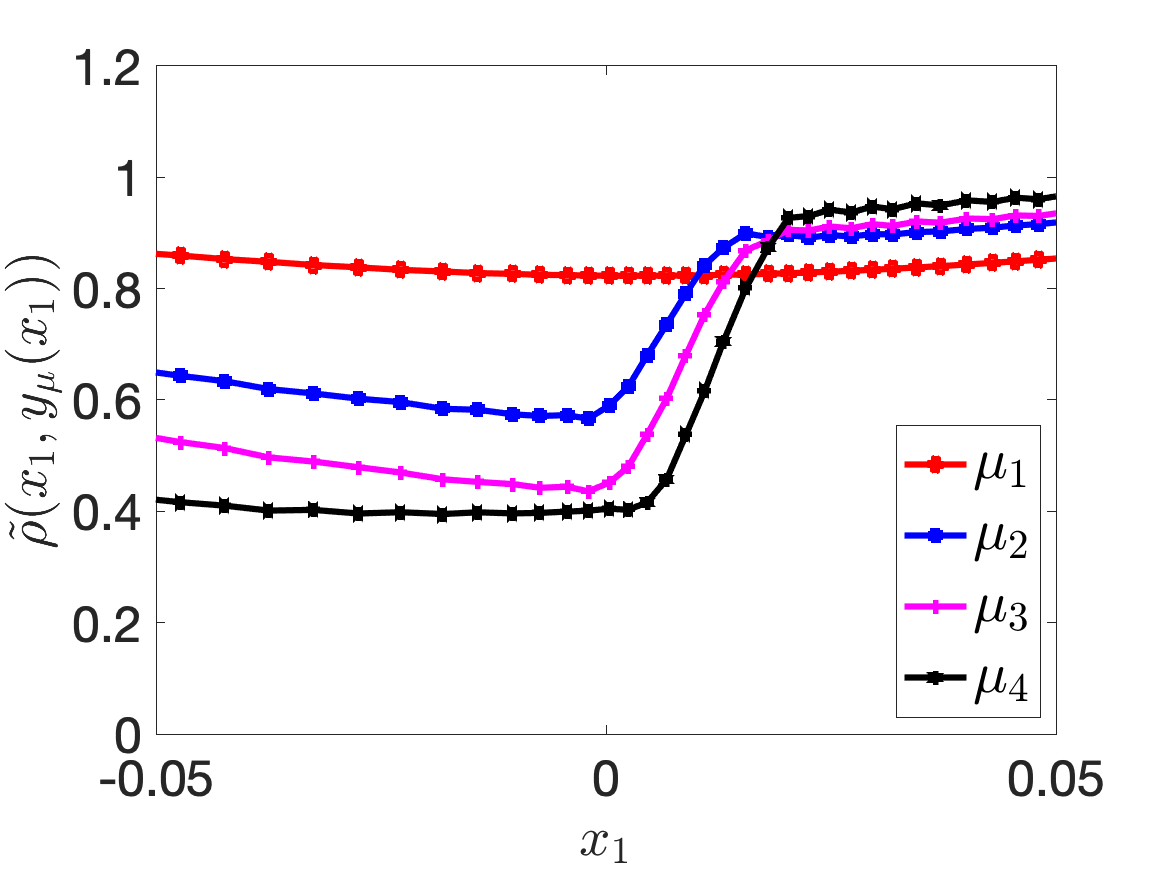}}
 
\caption{transonic bump. 
Behavior of the (modified) density field  in physical (cf. (a)-(b)-(c)) and reference 
(cf. (d)-(e)-(f))   configuration for four values of the parameter and three iterations of the adaptive algorithm (accelerated version) with isotropic mesh adaptation.
}
 \label{fig:transbump_accelerated_registration}
 \end{figure}

\begin{figure}[H]
\centering
 \subfloat[$it=1$] 
{  \includegraphics[width=0.33\textwidth]
 {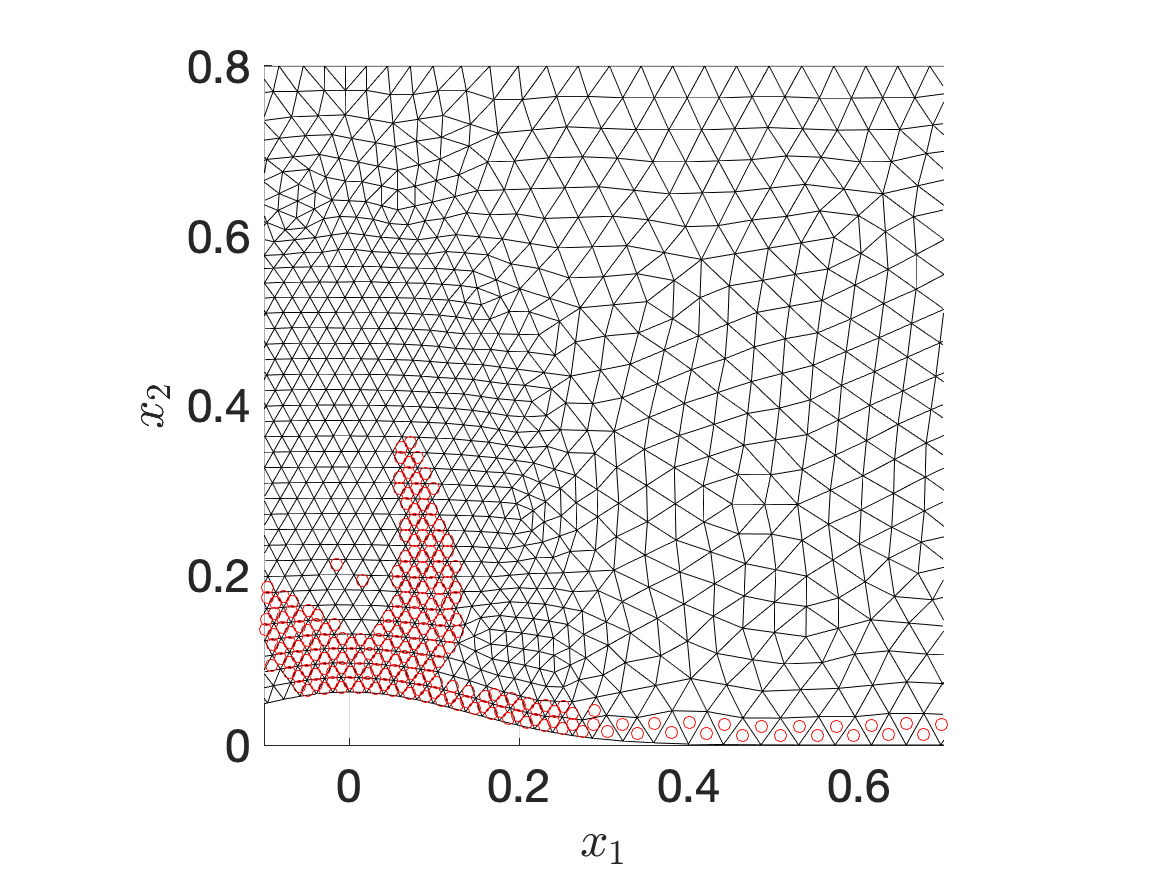}}
   ~~
 \subfloat[$it=2$] 
{  \includegraphics[width=0.33\textwidth]
 {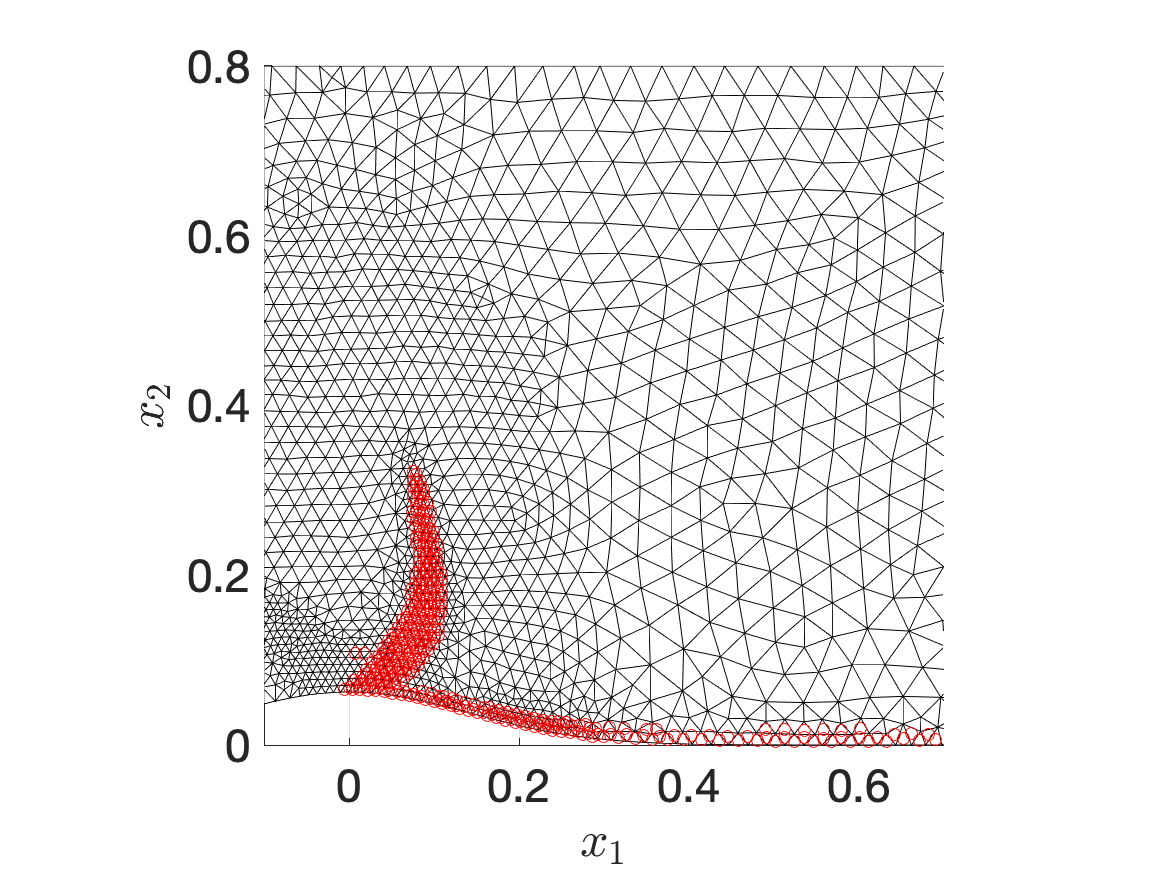}}
~~
 \subfloat[$it=3$] 
{  \includegraphics[width=0.33\textwidth]
 {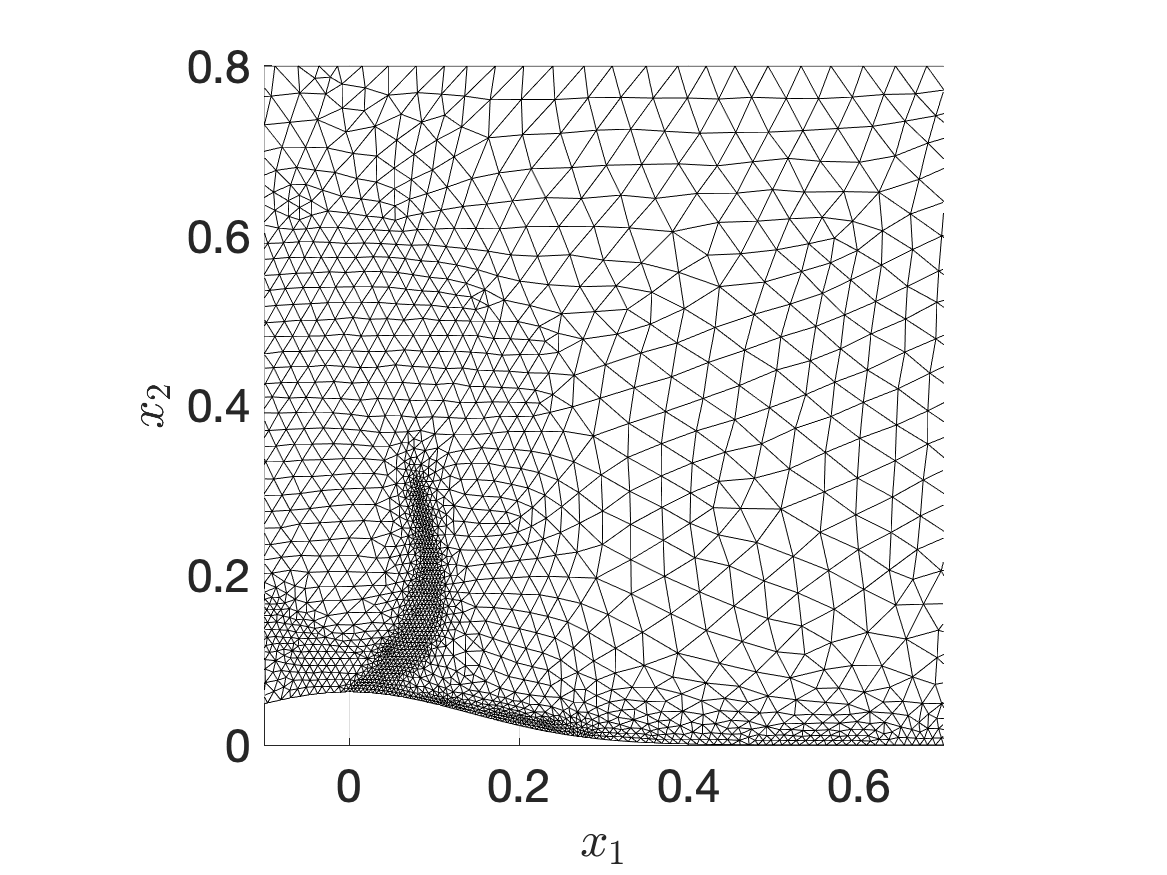}}
 
 \caption{transonic bump. Visualization of the reference mesh in the proximity of the bump, for three iterations of the adaptive algorithm (accelerated). 
 Red dots indicate the centers of the marked elements.
}
 \label{fig:transbump_accelerated_mesh}
 \end{figure}

 Figures 
\ref{fig:transbump_accelerated_rom}, 
\ref{fig:transbump_accelerated_registration} and
\ref{fig:transbump_accelerated_mesh}
  replicate the results of  Figures
 \ref{fig:transbump_basic_rom},
\ref{fig:transbump_basic_registration} and
\ref{fig:transbump_basic_mesh}: we observe that the results of the accelerated procedure are consistent with the ones obtained using the basic approach.
 
%transbump_basic_rom

%transbump_registration

%transbump_meshadapt

Table \ref{tab:offline_costs_transbump_acceleration} shows the details of the offline training costs. We observe that the vast majority of the computational gain is due to the reduction in the costs of the  HF solves and also in the overhead of the greedy algorithm.
We further remark that the acceleration strategy enables a much more efficient parallelization of the offline stage.
First, the computation of the initial set of HF solutions (cf. Line 1, Algorithm \ref{alg:weak_greedy}) is embarrassingly parallel; second, the solution to the registration problems \eqref{eq:optimization_all_params} based on the proposed initialization method is also parallel.

\begin{table}[H]
\centering
\begin{tabular}{||c| c| c| c||} 
\hline 
& it = 1 & it = 2 &  it = 3 \\[0.5ex]  \hline \hline 
ROB size: &  $18$ &  $17$ & $21$ \\[0.5ex] \hline 
mesh size: & $3448$ &  $4659$ & $6663$ \\[0.5ex] \hline 
snapshot generation: & $1455.19$ & $51.90$ & $54.51$ \\[0.5ex]   \hline 
registration (sensor def.): & $195.00$ & $187.80$ & $212.40$ \\[0.5ex]  \hline 
registration (optimization):& $390.30$ & $448.93$ & $728.93$ \\[0.5ex]  \hline 
mesh adaptation:& $0.00$ & $0.36$ & $0.61$ \\[0.5ex] \hline 
greedy alg (HF solves):& $460.26$ & $725.28$ & $1433.38$ \\[0.5ex] \hline 
greedy alg (overhead): & $282.49$ & $284.42$ & $994.63$ \\[0.5ex] \hline 
PTC iterations (avg): & $8.28$ & $8.71$ & $8.62$ \\[0.5ex] \hline 
\end{tabular}
\caption{transonic bump. Offline training costs (in seconds) of the adaptive (accelerated) approach.}
\label{tab:offline_costs_transbump_acceleration}
\end{table}

\subsection{Acceleration with anisotropic mesh adaptation}
\begin{figure}[H]
\centering
 \subfloat[] 
{  \includegraphics[width=0.4\textwidth]
 {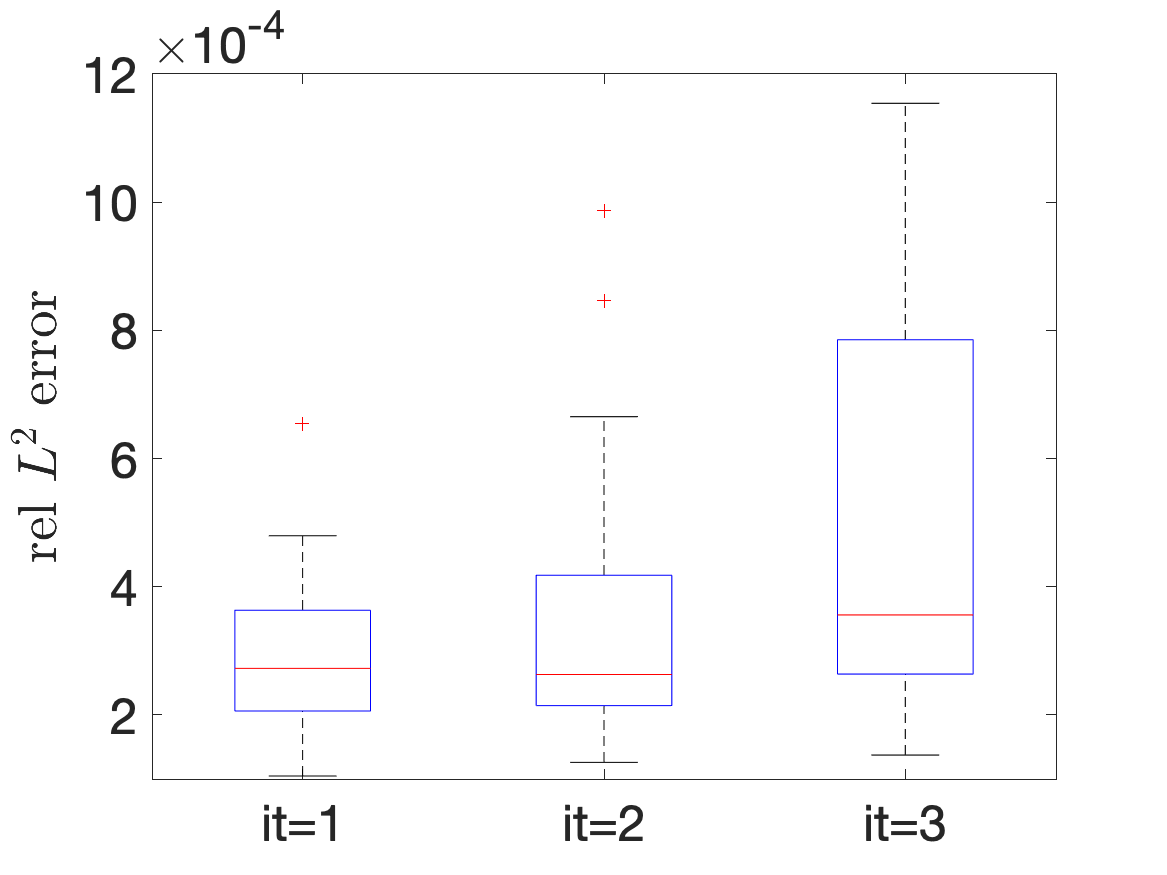}}
   ~~
 \subfloat[] 
{  \includegraphics[width=0.4\textwidth]
 {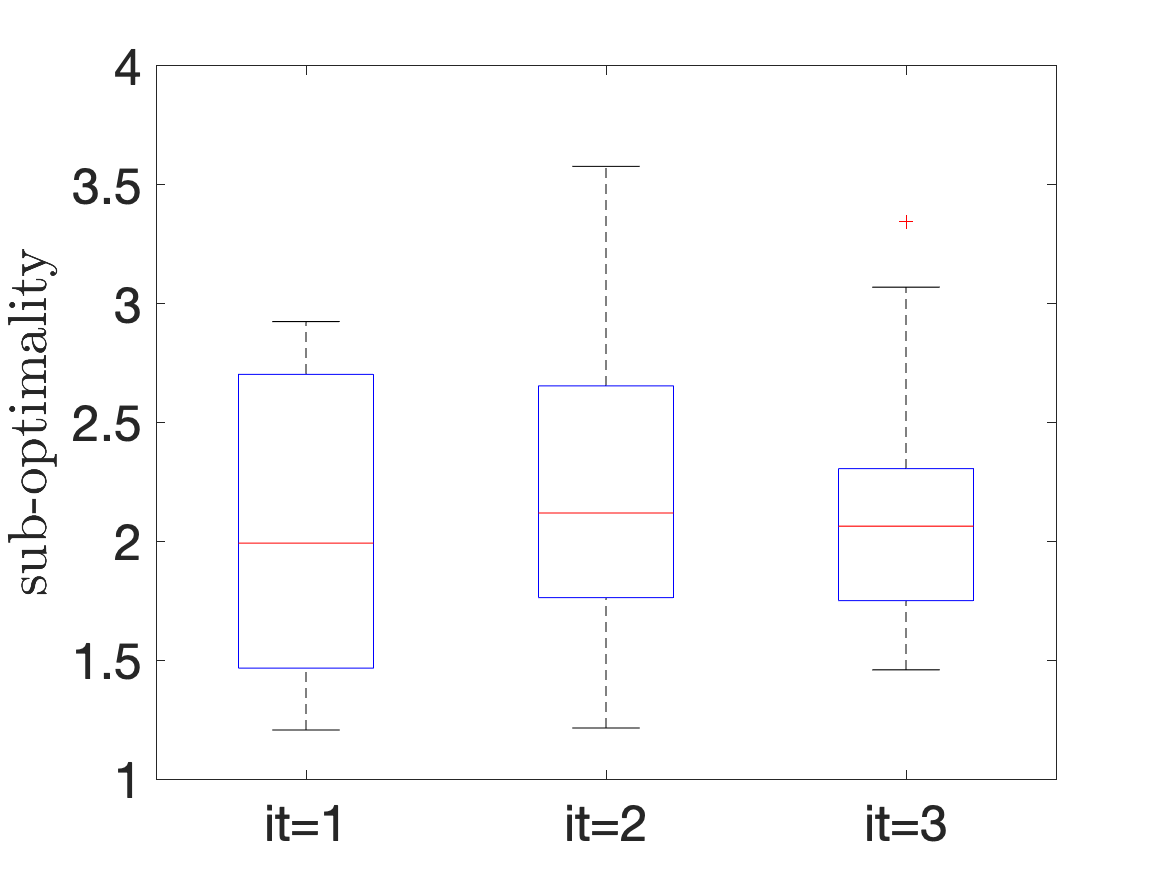}}

 \subfloat[] 
{  \includegraphics[width=0.4\textwidth]
 {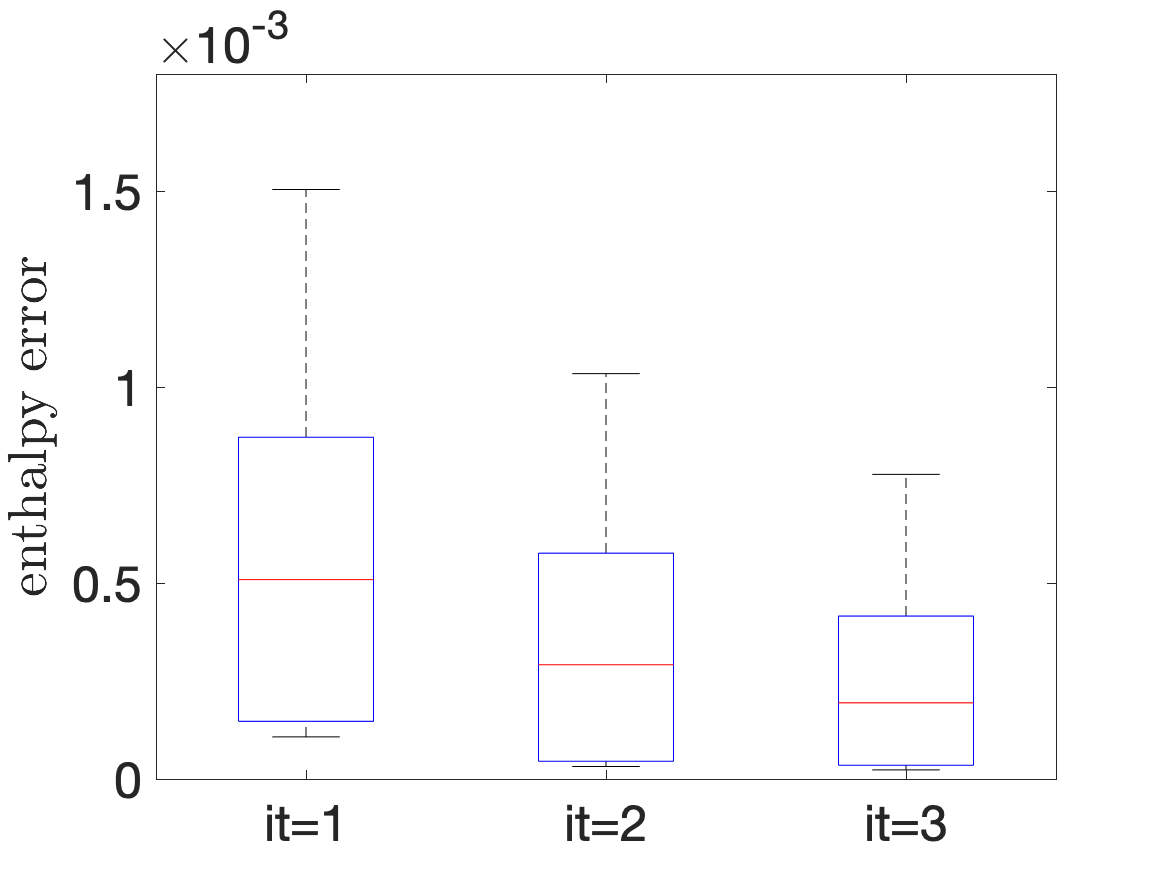}}
 ~~
  \subfloat[] 
{  \includegraphics[width=0.4\textwidth]
 {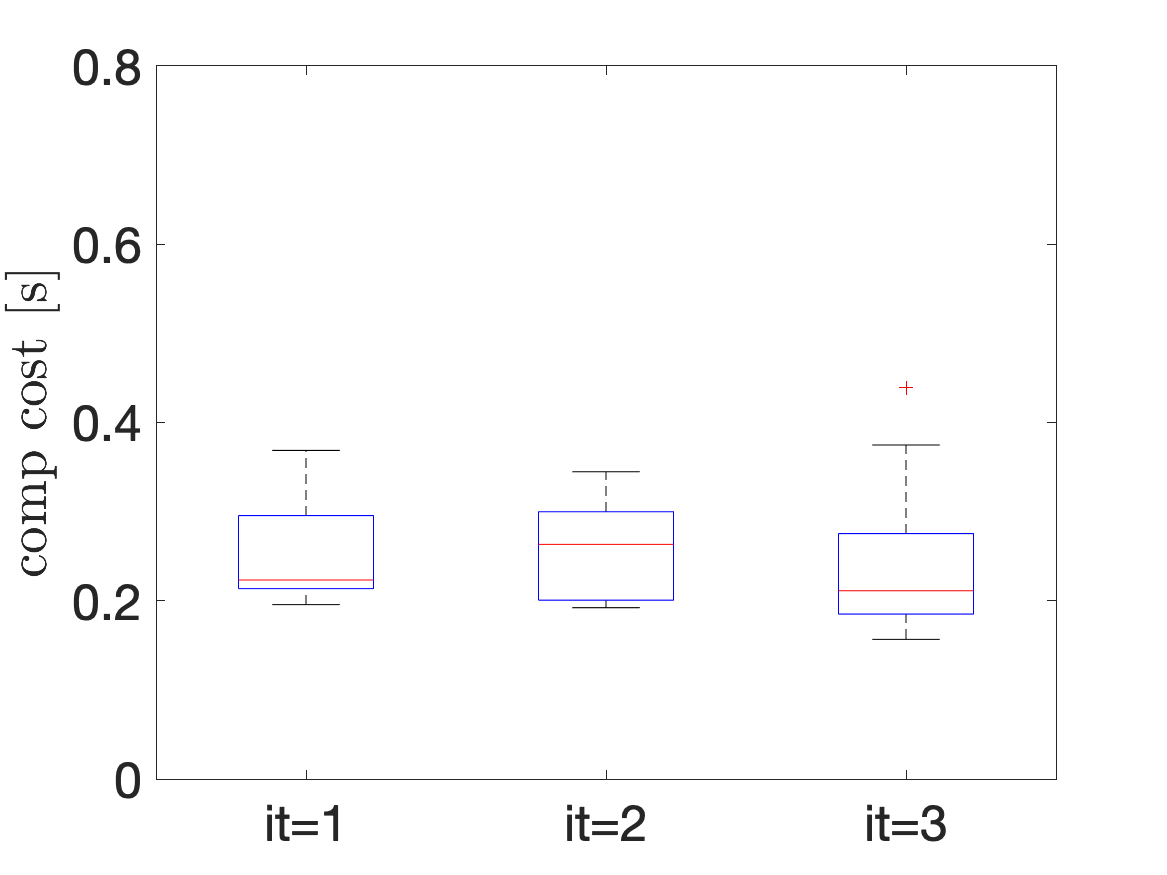}}
 
 \caption{transonic bump. Performance of the ROM for three iterations of the adaptive (accelerated) procedure with anisotropic mesh adaptation.
}
 \label{fig:transbump_accelerated_rom_aniso}
 \end{figure}

Figure \ref{fig:transbump_accelerated_rom_aniso} shows the performance of the ROM on the test set for the adaptive training procedure with anisotropic mesh adaptation: results are in good agreement with the results obtained using isotropic mesh adaptation. Similarly, Table 
\ref{tab:offline_costs_transbump_acceleration_aniso} details  the offline costs. As discussed in the main body of the paper,   anisotropic MA leads to a slight increase in the number of PTC iterations required for convergence. Note, however, that the number of iterations is still much lower than the one obtained with free-stream solution initialization.

\begin{table}[H]
\centering
\begin{tabular}{||c| c| c| c||} 
\hline 
& it = 1 & it = 2 &  it = 3 \\[0.5ex]  \hline \hline 
ROB size: &  $18$ &  $17$ & $16$ \\[0.5ex] \hline 
mesh size: & $3448$ &  $4440$ & $6389$ \\[0.5ex] \hline 
snapshot generation: & $1426.79$ & $51.47$ & $55.77$ \\[0.5ex]   \hline 
registration (sensor def.): & $177.64$ & $186.56$ & $191.43$ \\[0.5ex]  \hline 
registration (optimization):& $388.84$ & $442.69$ & $692.53$ \\[0.5ex]  \hline 
mesh adaptation:& $0.00$ & $1.01$ & $1.64$ \\[0.5ex] \hline 
greedy alg (HF solves):& $459.32$ & $910.38$ & $1357.74$ \\[0.5ex] \hline 
greedy alg (overhead): & $282.83$ & $270.21$ & $228.78$ \\[0.5ex] \hline 
PTC iterations (avg): & $8.28$ & $11.65$ & $10.19$ \\[0.5ex] \hline 
\end{tabular}
\caption{transonic bump. Offline training costs (in seconds) of the adaptive (accelerated) approach with anisotropic mesh adaptation.}
\label{tab:offline_costs_transbump_acceleration_aniso}
\end{table}
 
\bibliographystyle{abbrv}	
\bibliography{all_refs}

\end{document}